\documentclass[11pt]{article}%
\usepackage{amsfonts}
\usepackage{amssymb}
\usepackage{amsmath}
\usepackage{graphicx}
\usepackage{natbib}
\usepackage{epsfig}
\usepackage{subfigure}
\usepackage[left=2.6cm,top=2.4cm,right=2.6cm]{geometry}
\usepackage[doublespacing]{setspace}
\usepackage{fancyhdr}%
\newtheorem{theorem}{Theorem}
\newtheorem{corollary}{Corollary}
\newtheorem{definition}{Definition}

\newtheorem{proposition}{Proposition}
\newtheorem{remark}{Remark}
\newtheorem{example}{Example}
\newtheorem{identity}{Identity}
\newenvironment{proof}[1][Proof]{\noindent\textbf{#1.} }{$\hfill\square$\bigskip}

\newtheorem{defn0}{Definition}

\newtheorem{theorem0}{Theorem}
\newtheorem{proposition0}{Proposition}
\newtheorem{example0}{Example}
\newtheorem{remark0}{Remark}

\newtheorem{identity0}{Identity}
\renewenvironment{definition}{ \flushleft \begin{defn0} \rm}{\end{defn0}}
\renewenvironment{theorem}{\flushleft \begin{theorem0} \rm }{\end{theorem0}}
\renewenvironment{proposition}{ \flushleft \begin{proposition0} \rm }{\end{proposition0}}

\makeatletter
\renewcommand\section{\@startsection
{section}{1}{0pt}{3pt} {2pt} {\large\bf}}
\renewcommand\subsection{\@startsection
{subsection}{1}{0pt}{1pt} {1pt} {\large\bf}}
\renewcommand\subsubsection{\@startsection
{subsubsection}{1}{0pt}{-2.5ex plus-1ex minus-.2ex} {2.3ex plus.2ex}
{\normalfont\bf}}
\renewcommand{\thesection}{\arabic{section}}
\renewcommand{\thesubsection}{\arabic{section}.\arabic{subsection}}

\renewcommand{\cite}{\citet*}
\begin{document}
\bigskip%
\begin{center}
\large
{The Explicit Chaotic Representation of the powers of increments of L\'{e}%
vy Processes}
\vspace{0.2 in}

\large{Wing Yan Yip{$^{\dag}$}, David Stephens$^{\dag,\ddag}$, Sofia
Olhede$^{\dag}$}
\vspace{0.2 in}

\small{${\dag}$ Department of Mathematics, Imperial College London.}

\small{${\ddag}$ Department of Mathematics and Statistics, McGill
University.}

\small{wing.yip01@ic.ac.uk, dstephens@math.mcgill.ca, s.olhede@ic.ac.uk}

\small{\today}
\end{center}%

\rule{16cm}{0.03cm}

\noindent\textbf{Abstract}

An explicit formula for the chaotic representation of the powers of
increments, $\left(  X_{t+t_{0}}-X_{t_{0}}\right)  ^{n},$ of a L\'{e}vy
process is presented. \ There are two different chaos expansions of a square
integrable functional of a L\'{e}vy process: one with respect to the
compensated Poisson random measure and the other with respect to the
orthogonal compensated powers of the jumps of the L\'{e}vy process.
\ Computationally explicit formulae for both of these chaos expansions of
$\left(  X_{t+t_{0}}-X_{t_{0}}\right)  ^{n}$ are given in this paper.
\ Simulation results verify that the representation is satisfactory. \ The CRP
of a number of financial derivatives can be found by expressing them in terms
of $\left(  X_{t+t_{0}}-X_{t_{0}}\right)  ^{n}$ using Taylor's expansion.

\bigskip



\rule{16cm}{0.03cm}

\pagestyle{fancy}
\fancyhf{}
\fancyhead[LE,RO]{\textbf{\thepage}}
\renewcommand{\headrulewidth}{0pt}
\fancyhead[LO,RE]{\small{STATISTICS SECTION TECHNICAL REPORT TR-07-04}}

\newpage%

\section{Introduction\label{SectionIntroduction}\label{SectionReview}}

The chaotic representation of a square integrable functional of a L\'{e}vy
process is an expansion via its expectation plus a sum of iterated stochastic
integrals, see \cite{suj06} for a recent review of such representations.
\ There are two different types of chaos expansions: \cite{i56} proved a
Chaotic Representation Property (CRP) for any square integrable functional for
a general L\'{e}vy process. \ This representation is written using multiple
integrals with respect to a two-parameter random measure associated with the
L\'{e}vy process. \ \cite{ns00} proved the existence of a new version of the
CRP, which states that every square integrable L\'{e}vy functional can be
represented as its expectation plus an infinite sum of stochastic integrals
with respect to the orthogonalized compensated power jump processes of the
underlying L\'{e}vy process. \ \cite{bdlop03} and \cite{suj06} derived the
relationships between these two representations. \ However, these
representations are computationally intractable. \ For the powers of
increments, $\left(  X_{t+t_{0}}-X_{t_{0}}\right)  ^{n}$, of a L\'{e}vy
process, we instead derive computationally explicit formulae for the
integrands of these two chaotic expansions. \ Hence we have all the results
necessary to construct arbitrarily accurate computational formulae for the
L\'{e}vy functionals themselves.

Power jump processes are important in mathematical finance. \ \cite{bs06}
performed hypothesis tests on exchange data under the null of no jumps and
found that the tests were rejected frequently. \ In fact, at intraday scales,
prices move essentially by jumps and even at the scale of months, the
discontinuous behavior cannot be ignored in general. \ Only after
coarse-graining their behavior over longer time scales do we obtain something
similar to Brownian motion. \ Jumps can be understood both in terms of a
Poisson random measure, or equivalently, by using the Power jump processes.
\ Note that \cite[Proposition 2]{ns00} proved that all square integrable
random variables, adapted to the filtration generated by the L\'{e}vy process
denoted by $X,$ can be represented as a linear combination of powers of
increments of $X,$ see Section \ref{SectionLevy} below. \ In fact, for any
square integrable random variable, $F$, with derivatives of all order, we can
apply Taylor's Theorem to express $F$ in terms of a polynomial of powers of
increments of $X$. \ Thus, the CRP of a number of financial derivatives can be
found using this method, as is discussed further in Section
\ref{SectionGeneral}.

The derivation of an explicit formula for the CRP has been the focus of
considerable study, see for example \cite{ns01}, \cite{lsuv02}, \cite{l04} and
\cite{esv05}. \ All the explicit formulae for general L\'{e}vy functionals
derived in these papers use the Malliavin type derivatives to derive explicit
representations of stochastic processes for applications in finance. \ The
derivative operator $D$ is, in all of these cases, defined by its action on
the chaos expansions. \ In other words, the explicit chaos expansion must in
fact be known before $D$ can be applied to find the explicit form of the
predictable or chaotic representation, thus yielding a circular specification.
\ For example, \cite[Definition 1.7]{lsuv02} defined the derivative of $F$ in
the $l$-direction by%
\[
D_{t}^{\left(  l\right)  }F=\sum_{n=1}^{\infty}\sum_{i_{1},...,i_{n}}%
\sum_{k=1}^{n}1_{\left\{  i_{k}=l\right\}  }J_{n-1}^{\left(  i_{1}%
,...,\widehat{i_{k}},...,i_{n}\right)  }\left(  f_{i_{1},...,i_{n}}\left(
\cdots,t,\cdots\right)  1_{\Sigma_{n}^{\left(  k\right)  }\left(  t\right)
}\left(  \cdot\right)  \right)  ,
\]
and \cite[Section 3]{l04} defined the derivative operator by%
\[
D_{t,z}F=\sum_{n=1}^{\infty}nI_{n-1}\left(  f_{n}\left(  \cdot,t,z\right)
\right)  ,
\]
where%
\[
I_{n}\left(  f_{n}\right)  =\int_{\left[  0,T\right]  ^{n}\times\mathbb{R}%
_{0}^{n}}f_{n}\left(  t_{1},...,t_{n},z_{1},...,z_{n}\right)  \mathrm{d}%
\left(  \mu-\pi\right)  ^{\otimes n}.
\]
(Please refer to the corresponding papers for notation). \ Note that both of
these definitions require the knowledge of the functions $\left\{
f_{i_{1},...,i_{n}}\right\}  $'s or $f_{n}\left(  t_{1},...,t_{n}%
,z_{1},...,z_{n}\right)  $'s$,$ which are the integrands of the chaos
expansion of $F.$ \ 

\cite{j05} extended the CRP in \cite{ns00} to a large class of semimartingales
and derived the explicit representation of the power of a L\'{e}vy process
with respect to the corresponding non-compensated power jump processes, which
is discussed futher after Theorem \ref{newFormulaH} in this paper. \ Note that
L\'{e}vy processes are included in the class of semimartingales, see
\cite[Corollary 2.3.21, p.92]{kl01}. \ Our formula gives the explicit
representation with respect to the orthogonalized compensated power jump
processes. \ Our result is therefore complementary to Jamshidian's formula,
since our explicit formula gives the CRP with respect to the orthogonalized
processes, as defined by \cite{ns00}. \ 

In practical applications, it is often convenient to truncate the
representation given by the PRP. \ The truncated representation of a
stochastic process would yield a practically implementable approximation to
the stochastic process. \ This approximation would be used for simulating the
process, or with a finite number of traded higher order options, providing
pricing formulae. The truncation would be chosen with minimal variance
constraint. \ The advantage of expressing the sum in terms of stochastic
integrals with respect to the orthogonalized processes is that the error terms
omitted will be uncorrelated with the terms remained in the approximation.
\ Jamshidian's result holds for general semimartingales (a larger class than
ours) but our formula is designed for those with compensators equal to a
constant times $t$ only (which is satisfied by all L\'{e}vy processes). \ Our
results can be easily extended to semimartingales when the form of the
compensators is known. \ 

The rest of the paper is arranged as follow: Section \ref{SectionBackground}
gives the background information about the CRP for L\'{e}vy processes. \ We
give the explicit formulae for the CRP for $\left(  X_{t+t_{0}}-X_{t_{0}%
}\right)  ^{n}$ of a L\'{e}vy process in terms of power jump processes in
Section \ref{SectionPJ} and in terms of Poisson random measure in Section
\ref{Sectionprm}. \ We show that in the L\'{e}vy case, our formula complements
Jamshidian's formula. \ Section \ref{SectionGeneral} gives the representation
of a common kind of L\'{e}vy functionals with the use of Taylor's Theorem.
\ Simulation results for the explicit formulae are given in Section
\ref{SectionSimulation}. \ In Section \ref{SectionConclusion}, some concluding
remarks are provided. \ Proofs and plots are included as appendices at the end.

\section{Background\label{SectionBackground}}

\subsection{L\'{e}vy processes and their properties\label{SectionLevy}}

We give a brief account of L\'{e}vy processes and refer the reader to the work
by \cite{s99} for a more detailed account. \ A real-valued c\'{a}dl\'{a}g
stochastic process $X=\{X_{t},t\geq0\}$ defined in a complete probability
space $(\Omega,\mathcal{F},P)$ on $\mathbb{R}^{d}$ is called \textit{L\'{e}vy
process} if $X$ has stationary and independent increments with $X_{0}=0$,
where $\mathcal{F}$ is the filtration generated by $X:\mathcal{F}_{t}%
=\sigma\left\{  X_{t},0\leq s\leq t\right\}  .$ \ Denote the \textit{left
limit process} by $X_{t-}=\lim_{s\rightarrow t,s<t}X_{s},\ t>0,$ and the
\textit{jump size} at time $t$ by $\Delta X_{t}=X_{t}-X_{t-}$. \ 

A L\'{e}vy process is fully specified by its characteristic function. \ Let
$\phi_{X_{1}}\left(  u\right)  \,\ $be the \textit{characteristic function} of
the L\'{e}vy process at $t=1,$ $X_{1},$ that is, $\phi_{X_{1}}\left(
u\right)  =E\left[  e^{\mathrm{i}uX_{1}}\right]  .$ \ The characteristic
function of $X_{t}$ is then given by $\left(  \phi_{X_{1}}\left(  u\right)
\right)  ^{t}$ since the distribution of a L\'{e}vy process is infinitely
divisible, see \cite[chapter 2]{s99}. \ The cumulant characteristic function
$\psi\left(  u\right)  =\log\phi_{X_{1}}\left(  u\right)  $ is often called
the \textit{characteristic exponent}, which satisfies the
\textit{L\'{e}vy-Khintchine formula}:%
\begin{equation}
\psi\left(  u\right)  =\mathrm{i}\gamma u-\frac{1}{2}\sigma^{2}u^{2}%
+\int_{-\infty}^{+\infty}\left(  \exp\left(  \mathrm{i}ux\right)
-1-\mathrm{i}ux1_{\left\{  \left\vert x\right\vert <1\right\}  }\right)
\nu\left(  \mathrm{d}x\right)  ,\label{lk}%
\end{equation}
where $\gamma\in\mathbb{R},\ \sigma^{2}\geq0$ and $\nu$ is a measure on
$\mathbb{R}\backslash\left\{  0\right\}  $ with $\nu\left(  \left\{
0\right\}  \right)  =0$ and%
\[
\int_{-\infty}^{+\infty}\left(  1\wedge x^{2}\right)  \nu\left(
\mathrm{d}x\right)  <\infty.
\]
In general, a L\'{e}vy process consists of three independent components: a
linear deterministic component, a Brownian component and a pure jump
component. \ The L\'{e}vy measure $\nu\left(  \mathrm{d}x\right)  $ dictates
the jump process: jumps of sizes in set $A$ occur according to a Poisson
process with intensity parameter $\int_{A}\nu\left(  \mathrm{d}x\right)  $.
\ To model a generic L\'{e}vy process, only $\gamma,\ \sigma$ and a form for
$\nu\left(  \mathrm{d}x\right)  $ need to be specified.

In the rest of the paper, we assume that all L\'{e}vy measures concerned
satisfy, for some $\varepsilon>0$ and $\lambda>0$,
\begin{equation}
\int_{\left(  -\varepsilon,\varepsilon\right)  ^{c}}\exp\left(  \lambda
\left\vert x\right\vert \right)  \nu\left(  \mathrm{d}x\right)  <\infty
.\label{condition}%
\end{equation}
This condition implies that for $i\geq2,$ $\int_{-\infty}^{+\infty}\left\vert
x\right\vert ^{i}\nu\left(  \mathrm{d}x\right)  <\infty,$ and that the
characteristic function $E\left[  \exp\left(  \mathrm{i}uX_{t}\right)
\right]  $ is analytic in a neighborhood of 0.

Denote the $i$-th \textit{power jump process} by $X_{t}^{(i)}=\sum_{0<s\leq
t}(\Delta X_{s})^{i},\ i\geq2,\ $and for completeness let $X_{t}^{(1)}=X_{t}
$. \ In general, it is not true that $X_{t}=\sum_{0<s\leq t}\Delta X_{s}$;
this holds only in the bounded variation case, with $\sigma^{2}=0.$ \ By
definition, the quadratic variation of $X_{t},\ \left[  X,X\right]  _{t}%
=\sum_{0<s\leq t}(\Delta X_{s})^{2}=X_{t}^{(2)}$ when $\sigma^{2}=0$. \ The
power jump processes are also L\'{e}vy processes and jump at the same time as
$X_{t},$ but with jump sizes equal to the $i$-th powers of those of $X_{t}$,
see \cite{ns00}.

Clearly $E[X_{t}]=E[X_{t}^{(1)}]=m_{1}t$, where $m_{1}<\infty$ is a constant
and by \cite[p.32]{p04}, we have%
\begin{equation}
E[X_{t}^{(i)}]=E[\sum_{0<s\leq t}(\Delta X_{s})^{i}]=t\int_{-\infty}^{\infty
}x^{i}\nu(\mathrm{d}x)=m_{i}t<\infty,\ \ \ \mathrm{for}\ i\geq2,\label{mean}%
\end{equation}
thus defining $m_{i}$. \ \cite{ns00} introduced the \textit{compensated power
jump process }(or\textit{\ Teugels martingale}) of order $i$, $\left\{
Y_{t}^{\left(  i\right)  }\right\}  ,$ defined by
\begin{equation}
Y_{t}^{(i)}=X_{t}^{(i)}-E[X_{t}^{(i)}]=X_{t}^{(i)}-m_{i}t\ \ \ \mathrm{for}%
\ i=1,2,3,....\label{CompPower}%
\end{equation}
$Y_{t}^{\left(  i\right)  }$ is constructed to have a zero mean. \ It was
shown by \cite[Section 2]{ns00} that there exist constants $a_{i,1}%
,a_{i,2},...,a_{i,i-1}$ such that the processes defined by
\begin{equation}
H^{(i)}=Y^{(i)}+a_{i,i-1}Y^{(i-1)}+\cdots+a_{i,1}Y^{(1)},\label{H}%
\end{equation}
for$\ i\geq1$ are a set of pairwise strongly orthogonal martingales, and this
implies that for $i\neq j$, the process $H^{\left(  i\right)  }H^{\left(
j\right)  }$ is a martingale, see \cite{lsuv02}. \ For convenience, we define
$a_{i,i}=1$. \ \ \cite{ns00} proved that this strong orthogonality is
equivalent to the existence of an orthogonal family of polynomials with
respect to the measure%
\[
\mathrm{d}\eta\left(  x\right)  =\sigma^{2}\mathrm{d}\delta_{0}\left(
x\right)  +x^{2}\nu(\mathrm{d}x),
\]
where $\delta_{0}\left(  x\right)  =1$ when $x=0$ and zero otherwise, that is,
the polynomials $p_{n}$ defined by
\[
p_{n}\left(  x\right)  =\sum_{j=1}^{n}a_{n,j}x^{j-1}%
\]
are orthogonal with respect to the measure $\eta$:%
\[
\int_{\mathbb{R}}p_{n}\left(  x\right)  p_{m}\left(  x\right)  \mathrm{d}%
\eta\left(  x\right)  =0,\ \ \ n\neq m.
\]
We now state some key related results in the representation of stochastic
processes given in \cite{ns00}.

\vspace{-0.1 in}

\begin{itemize}
\item \textit{Denseness of polynomials} (\cite[Proposition 2]{ns00}): Let
$\mathcal{P}=\{X_{t_{1}}^{k_{1}}(X_{t_{2}}-X_{t_{1}})^{k_{2}}\cdots(X_{t_{n}%
}-X_{t_{n-1}})^{k_{n}}:n\geq0,\ 0\leq t_{1}<t_{2}<\cdots<t_{n},k_{1}%
,...,k_{n}\geq1\}$ be a family of stochastic processes. \ Then we have that
$\mathcal{P}$ is a total family in $L^{2}\left(  \Omega,\mathcal{F}%
_{T},P\right)  ,$ that is, the linear subspace spanned by $\mathcal{P}$ is
dense in $L^{2}\left(  \Omega,\mathcal{F}_{T},P\right)  $; each element in
$L^{2}(\Omega,\mathcal{F})$ can be represented as a linear combination of
elements in $\mathcal{P}$. \ 

\item \textit{Chaotic Representation Property} (CRP): Every random variable
$F$ in $L^{2}(\Omega,\mathcal{F})$ has a representation of the form
\begin{equation}
F=E(F)+\sum_{j=1}^{\infty}\sum_{i_{1},...,i_{j}\geq1}\int_{0}^{\infty}\int
_{0}^{{t_{1}}-}\cdots\int_{0}^{{t_{j-1}}-}f_{(i_{1},...,i_{j})}(t_{1}%
,...,t_{j})\mathrm{d}H_{t_{j}}^{(i_{j})}\cdots\mathrm{d}H_{t_{2}}^{(i_{2}%
)}\mathrm{d}H_{t_{1}}^{(i_{1})},\label{crp}%
\end{equation}

\noindent where the $f_{(i_{1},...,i_{j})}$'s are functions in $L^{2}%
(\mathbb{R}_{+}^{j})$. \ This result means that every random variable in
$L^{2}\left(  \Omega,\mathcal{F}_{T},P\right)  $ can be expressed as its
expectation plus an infinite sum of zero mean stochastic integrals with
respect to the orthogonalized compensated power jump processes of the
underlying L\'{e}vy process. \ Note that this representation does not
explicitly allow for calculation of the integrands.

\item \textit{Predictable Representation Property} (PRP) : From the CRP stated
above, we note that every random variable $F$ in $L^{2}\left(  \Omega
,\mathcal{F}_{T},P\right)  $ has a representation of the form%
\begin{equation}
F=E\left[  F\right]  +\sum_{i=1}^{\infty}\int_{0}^{\infty}\phi_{s}^{\left(
i\right)  }\ \mathrm{d}H_{s}^{\left(  i\right)  },\label{prp}%
\end{equation}
where $\phi_{s}^{\left(  i\right)  }$'s are predictable, that is, they are
$\mathcal{F}_{s-}$-measurable.
\end{itemize}

\subsection{Jamshidian's notation\label{SectionJam}}

In \cite{j05}, which extends the CRP to semimartingales, the power jump
processes and compensators were denoted and defined differently from
\cite{ns00}. \ The power jump processes were defined in \cite{j05} by \
\begin{equation}
\left[  X\right]  _{t}^{\left(  2\right)  }=\left[  X^{c}\right]  _{t}%
+\sum_{s\leq t}\left(  \Delta X_{s}\right)  ^{2}\text{ and }\left[  X\right]
_{t}^{\left(  n\right)  }=\sum_{s\leq t}\left(  \Delta X_{s}\right)
^{n}\text{ for }n=3,4,5,...,\label{JpowerJump}%
\end{equation}
where $\left[  X^{c}\right]  _{t}=\left[  X\right]  _{t}^{c}$ is the
continuous finite-variation (not martingale) part of $\left[  X\right]
_{t}^{\left(  2\right)  }$. \ Note that Jamshidian suppressed the time index
$t,$ but we add it here for clarification. \ Jamshidian denoted the
compensator of $\left[  X\right]  _{t}^{\left(  n\right)  }$ by $\left\langle
X\right\rangle _{t}^{\left(  n\right)  }.$ \ The compensator, $\left\langle
X\right\rangle _{t}^{\left(  n\right)  },$ is the predictable right-continuous
finite variation process such that $\left[  X\right]  _{t}^{\left(  n\right)
}-\left\langle X\right\rangle _{t}^{\left(  n\right)  }$ is a uniformly
integrable martingale. \ The compensated power jump process, denoted by
$X_{t}^{\left(  n\right)  }$, is thus defined by%
\begin{equation}
X_{t}^{\left(  n\right)  }=\left[  X\right]  _{t}^{\left(  n\right)
}-\left\langle X\right\rangle _{t}^{\left(  n\right)  }\text{ for
}n=2,3,4,....\label{Jcomp}%
\end{equation}

For L\'{e}vy processes, the compensators have the form $m_{i}t$, where
$m_{1}t=E[X_{t}]$ and $m_{i}t=t\int_{-\infty}^{\infty}x^{i}\nu(\mathrm{d}x) $
for $i=2,3,4,....$ \ However, for semimartingales, the general form of the
compensators is not known. \ 

\section{The chaotic representation with respect to power jump processes
\label{SectionPJ}}

In this section we firstly derive the explicit formulae for the CRP when $F$
in (\ref{crp}) is the power increment of a pure jump L\'{e}vy process and
extend it to a general L\'{e}vy process.\ 

\subsection{Pure jump case}

Let us first outline the form of the representation to introduce the reader to
the flavour of results in this paper. \ Suppose $t_{0}\geq0$ and let $G_{t}$
be a pure jump L\'{e}vy process with no Brownian part (that is, $\sigma^{2}%
=0$), $G_{t}^{\left(  i\right)  }$ be its $i$-th power jump process and
$\hat{G}_{t}^{\left(  i\right)  }$ be its $i$-th compensated power jump
process. \ Based on the structure of the expressions for $\left(  G_{t+t_{0}%
}-G_{t_{0}}\right)  ^{3}$ and $\left(  G_{t+t_{0}}-G_{t_{0}}\right)  ^{4}$
calculated using (\ref{Corrected1})-(\ref{Corrected})$,$ we desire to derive a
general formula for $\left(  G_{t+t_{0}}-G_{t_{0}}\right)  ^{k}%
,\ k=1,2,3,...,$ as this forms a starting point for the representation of
$X_{t} $. \ We notice that the number of stochastic integrals in each of the
above representation is less than the possible full representation specified
in the CRP by \cite{ns00}:%
\begin{align*}
\left(  X_{t+t_{0}}-X_{t_{0}}\right)  ^{k}  & =f^{\left(  k\right)  }\left(
t,t_{0}\right)  +\sum_{j=1}^{k}\sum_{\substack{\left(  i_{1},...,i_{j}\right)
\\\in\left\{  1,...,k\right\}  ^{j}}}\int_{t_{0}}^{t+t_{0}}\int_{t_{0}}%
^{t_{1}-}\cdots\int_{t_{0}}^{t_{j-1}-}f_{\left(  i_{1},...,i_{j}\right)
}^{\left(  k\right)  }\left(  t,t_{0},t_{1},...,t_{j}\right) \\
& \times\mathrm{d}Y_{t_{j}}^{\left(  i_{j}\right)  }\cdots\mathrm{d}Y_{t_{2}%
}^{\left(  i_{2}\right)  }\mathrm{d}Y_{t_{1}}^{\left(  i_{1}\right)  },
\end{align*}
where the $f_{\left(  i_{1},...,i_{j}\right)  }^{\left(  k\right)  }$'s are
deterministic functions in $L^{2}\left(  R_{+}^{j}\right)  $ and $Y$'s are
defined in (\ref{CompPower}). \ In $\left(  G_{t+t_{0}}-G_{t_{0}}\right)
^{2}$, we have $\int_{t_{0}}^{t+t_{0}}\int_{t_{0}}^{t_{1}-}$\textrm{d}$\hat
{G}_{t_{2}}^{\left(  1\right)  }$\textrm{d}$\hat{G}_{t_{1}}^{\left(  1\right)
},\ \int_{t_{0}}^{t+t_{0}}$\textrm{d}$\hat{G}_{t_{1}}^{\left(  1\right)  }$
and $\int_{t_{0}}^{t+t_{0}}$\textrm{d}$\hat{G}_{t_{1}}^{\left(  2\right)  }$,
that we shall represent via the list $\left\{  \left(  1,1\right)  ,\left(
1\right)  ,\left(  2\right)  \right\}  .$ \ We can do an equivalent
representation of $\left(  G_{t+t_{0}}-G_{t_{0}}\right)  ^{3}$ and $\left(
G_{t+t_{0}}-G_{t_{0}}\right)  ^{4}$ to get the following two lists:%
\begin{align*}
& \left\{  \left(  1,1,1\right)  ,\left(  1,1\right)  ,\left(  1,2\right)
,\left(  2,1\right)  ,\left(  1\right)  ,\left(  2\right)  ,\left(  3\right)
\right\}  .\\
& \left\{  \left(  1,1,1,1\right)  ,\left(  1,1,1\right)  ,\left(
1,1,2\right)  ,\left(  1,2,1\right)  ,\left(  2,1,1\right)  ,\right. \\
& \left.  \left(  1,1\right)  ,\left(  1,2\right)  ,\left(  2,1\right)
,\left(  2,2\right)  ,\left(  1,3\right)  ,\left(  3,1\right)  ,\left(
1\right)  ,\left(  2\right)  ,\left(  3\right)  ,\left(  4\right)  \right\}  .
\end{align*}
In general, the list of the orders of the compensated power jump processes of
the stochastic integrals in $\left(  G_{t+t_{0}}-G_{t_{0}}\right)  ^{k}$
depends on the collection of numbers%
\begin{equation}
\mathcal{I}_{k}=\left\{  \left(  i_{1},i_{2},...,i_{j}\right)  |\ j\in\left\{
1,2,...,k\right\}  ,\ i_{p}\in\left\{  1,2,...,k\right\}  \text{ and }%
\sum_{p=1}^{j}i_{p}\leq k\right\}  .\label{setI}%
\end{equation}
This construction is explained in the beginning of the proof of Theorem
\ref{Theorem3} (Appendix \ref{AppendixB}) using induction. \ A typical element
$\left(  i_{1},i_{2},...,i_{j}\right)  $ in $\mathcal{I}_{k}$ therefore
indexes a multiple stochastic integral $j$-times repeated with respect to the
power jump processes with powers $i_{1},i_{2},...,i_{j}$ and indexed
$t_{1},t_{2},...,t_{j}$. \ That is, $\left(  i_{1},i_{2},...,i_{j}\right)  $
indexes the integral
\[
\int_{t_{0}}^{t+t_{0}}\int_{t_{0}}^{t_{1}-}\cdots\int_{t_{0}}^{t_{j-1}%
}\mathrm{d}\hat{G}_{t_{j}}^{\left(  i_{1}\right)  }\cdots\mathrm{d}\hat
{G}_{t_{2}}^{\left(  i_{j-1}\right)  }\mathrm{d}\hat{G}_{t_{1}}^{\left(
i_{j}\right)  }.
\]

Next we consider the terms in the representation not involving any stochastic
integrals. \ That is, in $\left(  G_{t+t_{0}}-G_{t_{0}}\right)  ^{2} $,
$m_{1}^{2}t^{2}+m_{2}t$ is considered; in $\left(  G_{t+t_{0}}-G_{t_{0}%
}\right)  ^{3}$, $m_{1}^{3}t^{3}+3m_{1}m_{2}t^{2}+m_{3}t$ is considered, and
in $\left(  G_{t+t_{0}}-G_{t_{0}}\right)  ^{4}$, $m_{1}^{4}t^{4}+6m_{1}%
^{2}m_{2}t^{3}+\left(  4m_{1}m_{3}+3m_{2}^{2}\right)  t^{2}+m_{4}t$ is
considered. We use (\ref{Corrected1})-(\ref{Corrected}), given in Appendix
\ref{SectionNSTypo}, to derive the representation. \ This time the
representation can be simplified a great deal since we are not considering any
stochastic integrals. \ Denote the terms which do not contain any stochastic
integrals in $\left(  G_{t+t_{0}}-G_{t_{0}}\right)  ^{k} $ by $C_{t}^{\left(
k\right)  }$.

\begin{proposition}
\label{Theorem1} $C_{0}^{\left(  r\right)  }=0$ for all $r$, $C_{t}^{\left(
0\right)  }=1,\ C_{t}^{\left(  1\right)  }=m_{1}t,\ $and for $k=2,3,4,...,$%
\begin{equation}
C_{t}^{\left(  k\right)  }=\sum_{j=1}^{k-1}\binom{k}{j}m_{j}tC_{t}^{\left(
k-j\right)  }-\sum_{j=1}^{k-1}\binom{k}{j}m_{j}\int_{0}^{t}t_{1}%
\ \mathrm{d}C_{t_{1}}^{\left(  k-j\right)  }+m_{k}t.\label{C}%
\end{equation}

\end{proposition}

\qquad\qquad\qquad\qquad\qquad\qquad\qquad\qquad\qquad\qquad\qquad\qquad
\qquad\qquad\qquad\qquad\qquad\qquad\qquad\qquad\qquad\qquad

\begin{proof}
The results for $C_{0}^{\left(  r\right)  }$ and $C_{t}^{\left(  0\right)  }$
are trivial$.$ \ For $k=1$,$\ \left(  G_{t+t_{0}}-G_{t_{0}}\right)
=\int_{t_{0}}^{t+t_{0}}$\textrm{d}$\hat{G}_{t_{1}}^{\left(  1\right)  }%
+m_{1}t$ and hence $C_{t}^{\left(  1\right)  }=m_{1}t.$ \ For $k\geq2,$ the
terms in (\ref{Corrected1}) are equal to zero since $G_{t}$ has no Brownian
part ($\sigma^{2}=0$). \ The first term in (\ref{Corrected2}) contains a
stochastic integral and hence from the second term of (\ref{Corrected2}) and
(\ref{Corrected}), we have
\[
C_{t}^{\left(  k\right)  }=\sum_{j=1}^{k-1}\binom{k}{j}m_{j}\left(
t+t_{0}\right)  C_{t}^{\left(  k-j\right)  }-\sum_{j=1}^{k-1}\binom{k}{j}%
m_{j}\int_{t_{0}}^{t+t_{0}}t_{1}\ \mathrm{d}C_{t_{1}-t_{0}}^{\left(
k-j\right)  }+m_{k}t.
\]
Putting $u=t_{1}-t_{0}$ in the second term, we have%
\begin{align*}
C_{t}^{\left(  k\right)  }  & =\sum_{j=1}^{k-1}\binom{k}{j}m_{j}\left(
t+t_{0}\right)  C_{t}^{\left(  k-j\right)  }-\sum_{j=1}^{k-1}\binom{k}{j}%
m_{j}\int_{0}^{t}\left(  u+t_{0}\right)  \mathrm{d}C_{u}^{\left(  k-j\right)
}+m_{k}t\\
& =\sum_{j=1}^{k-1}\binom{k}{j}m_{j}tC_{t}^{\left(  k-j\right)  }-\sum
_{j=1}^{k-1}\binom{k}{j}m_{j}\int_{0}^{t}t_{1}\ \mathrm{d}C_{t_{1}}^{\left(
k-j\right)  }+m_{k}t.
\end{align*}
Note that $C_{t}^{\left(  k\right)  }$ is independent of $t_{0}.$ \ \ 
\end{proof}

Thus, given Proposition \ref{Theorem1}, $C_{t}^{\left(  k\right)  }$ can be
expressed in terms of $m_{i}$'s for any given $k$ and easily coded. \ We will
show in the followings that in the calculation of $\left(  G_{t+t_{0}%
}-G_{t_{0}}\right)  ^{k}$, all the $C_{t}^{\left(  j\right)  }$%
's$,\ j=0,1,...,k$ are required. \ In fact the coefficients of the stochastic
integrals in the representation depend only on $C_{t}^{\left(  j\right)  }%
$'s$,\ j=0,1,...,k$, as stated in Theorem \ref{Theorem3} below.

The next proposition gives the representation for $C_{t}^{\left(  k\right)  }$
in a non-recursive form. \ Let
\begin{equation}
\mathcal{L}_{k}=\left\{  \left(  i_{1},i_{2},...,i_{l}\right)  |l\in\left\{
1,2,...,k\right\}  ,i_{q}\in\left\{  1,2,...,k\right\}  ,i_{1}\geq i_{2}%
\geq\cdots\geq i_{l}\text{ and }\sum_{q=1}^{l}i_{q}=k\right\}  .\label{Lk}%
\end{equation}
The number of distinct values in a tuple $\phi_{k}=\left(  i_{1}^{(k)}%
,i_{2}^{(k)},...,i_{l}^{(k)}\right)  $ in $\mathcal{L}_{k}$ is less than or
equal to $l.$ \ When it is less than $l,$ it means some of the value(s) in the
tuple are repeated. \ Let the number of times $r\in\left\{
1,2,3,..,k\right\}  $ appears in the tuple $\phi_{k}=\left(  i_{1}^{(k)}%
,i_{2}^{(k)},...,i_{l}^{(k)}\right)  $ be $p_{r}^{\phi_{k}}.$

\begin{proposition}
\label{TheoremC2}%
\begin{equation}
C_{t}^{\left(  k\right)  }=\sum_{\phi_{k}=\left(  i_{1}^{(k)},i_{2}%
^{(k)},...,i_{l}^{(k)}\right)  \in\mathcal{L}_{k}}\frac{1}{l!}\left(
i_{1}^{(k)},i_{2}^{(k)},...,i_{l}^{(k)}\right)  !\left(  p_{1}^{\phi_{k}%
},p_{2}^{\phi_{k}},...,p_{k}^{\phi_{k}}\right)  !\left[  \prod\limits_{q\in
\phi_{k}}m_{q}\right]  t^{l}\label{C2_crp}%
\end{equation}
where $i_{1}^{(k)},...,i_{l}^{(k)}$ are the elements of $\phi_{k}$,
$p_{j}^{\phi_{k}}$'s are defined above and $\left(  i_{1}^{(k)},i_{2}%
^{(k)},...,i_{l}^{(k)}\right)  !$ is the multinomial coefficient: $\left(
i_{1}^{(k)},i_{2}^{(k)},...,i_{l}^{(k)}\right)  !=\frac{\left(  \sum_{j=1}%
^{l}i_{j}^{(k)}\right)  !}{i_{1}^{(k)}!i_{2}^{(k)}!\cdots i_{l}^{(k)}!}$
\end{proposition}

\begin{proof}
Proof is included in Appendix \ref{AppendixTheoremC2}. \ 
\end{proof}

Let $\Pi_{(i_{1},i_{2},...,i_{j}),t}^{\left(  k\right)  }$ be the coefficient
of $\int_{t_{0}}^{t+t_{0}}\int_{t_{0}}^{t_{1}-}\cdots\int_{t_{0}}^{t_{j-1}-}
$\textrm{d}$\hat{G}_{t_{j}}^{\left(  i_{1}\right)  }\cdots$\textrm{d}$\hat
{G}_{t_{2}}^{\left(  i_{j-1}\right)  }$\textrm{d}$\hat{G}_{t_{1}}^{\left(
i_{j}\right)  }$ in $\left(  G_{t+t_{0}}-G_{t_{0}}\right)  ^{k}.$ \ We then
have the following result. \ \ \qquad\qquad\qquad\qquad\qquad\qquad
\qquad\qquad\qquad\qquad\qquad\qquad\qquad\qquad\qquad\qquad\qquad\qquad
\qquad\qquad\qquad\qquad\qquad\qquad\qquad\qquad\qquad\qquad\qquad\qquad

\begin{proposition}
\label{Theorem2}%
\begin{equation}
\Pi_{(i_{1},i_{2},...,i_{j}),t}^{\left(  k\right)  }=\left(  i_{1}%
,i_{2},...,i_{j},n\right)  !C_{t}^{\left(  n\right)  }\text{ where }%
n=k-\sum_{p=1}^{j}i_{p}\text{.}\label{PI}%
\end{equation}
$\ $
\end{proposition}

\begin{proof}
The proof of Proposition \ref{Theorem2} is contained in the proof of Theorem
\ref{Theorem3}. \ \ 
\end{proof}

For example, say we want to find the coefficient of $\int_{t_{0}}^{t+t_{0}%
}\int_{t_{0}}^{t_{1}-}$\textrm{d}$\hat{G}_{t_{2}}^{\left(  1\right)  }%
$\textrm{d}$\hat{G}_{t_{1}}^{\left(  1\right)  }$ in $\left(  G_{t+t_{0}%
}-G_{t_{0}}\right)  ^{4}$, that is, we want to find $\Pi_{\left(  1,1\right)
,t}^{\left(  4\right)  }.$ \ To derive this coefficient, we first note that
$n=2$ and so $\Pi_{\left(  1,1\right)  ,t}^{\left(  4\right)  }=\frac
{4!}{1!1!2!}C_{t}^{\left(  2\right)  }=12\left(  m_{2}t+m_{1}^{2}t^{2}\right)
,$ which can be easily verified by calculating $\left(  G_{t+t_{0}}-G_{t_{0}%
}\right)  ^{4}$ using (\ref{Corrected1})-(\ref{Corrected}). \ Now we put the
above results together to get a general formula for $\left(  G_{t+t_{0}%
}-G_{t_{0}}\right)  ^{k}.$

\begin{theorem}
\label{Theorem3} Let $G_{t}$ be a L\'{e}vy process with no Brownian part
satisfying condition (\ref{condition}). \ Then the power of its increment can
be expressed by:%
\begin{equation}
\left(  G_{t+t_{0}}-G_{t_{0}}\right)  ^{k}=\sum_{\theta_{k}\in\mathcal{I}_{k}%
}\Pi_{\theta_{k},t}^{\left(  k\right)  }\mathcal{S}_{\theta_{k},t,t_{0}}%
+C_{t}^{\left(  k\right)  },\label{Thm4}%
\end{equation}
where $\mathcal{I}_{k}$ is defined in (\ref{setI}), $\Pi_{\theta_{k}%
,t}^{\left(  k\right)  }$ is defined in Proposition \ref{Theorem2}, the
$C_{t}^{\left(  k\right)  }$ are constants defined in Proposition
\ref{TheoremC2} and $\mathcal{S}_{\left(  i_{1},i_{2},...,i_{j}\right)
,t,t_{0}}$ is defined as the integral:%
\[
\mathcal{S}_{\left(  i_{1},i_{2},...,i_{j}\right)  ,t,t_{0}}=\int_{t_{0}%
}^{t+t_{0}}\int_{t_{0}}^{t_{1}-}\cdots\int_{t_{0}}^{t_{j-1}-}\mathrm{d}\hat
{G}_{t_{j}}^{\left(  i_{1}\right)  }\cdots\mathrm{d}\hat{G}_{t_{2}}^{\left(
i_{j-1}\right)  }\mathrm{d}\hat{G}_{t_{1}}^{\left(  i_{j}\right)  }.
\]

\end{theorem}

\begin{proof}
Proof is included in Appendix \ref{AppendixB}. \ \ 
\end{proof}

To derive the explicit formula for the power of increment of a L\'{e}vy
process, $\left(  X_{t+t_{0}}-X_{t_{0}}\right)  ^{n},$ with respect to
orthogonalized compensated power jump processes, we need the following proposition.

\begin{proposition}
\label{PropositionYtoH} The $n$-th compensated power jump processes,
$Y^{\left(  n\right)  }$, of a general L\'{e}vy processes satisfying condition
(\ref{condition}), can be expressed in terms of the orthogonalized compensated
power jump processes, $H^{\left(  i\right)  }$ for $i=1,2,...,n$, by%
\[
Y^{\left(  n\right)  }=H^{\left(  n\right)  }+\sum_{k=1}^{n-1}b_{n,k}%
H^{\left(  k\right)  },
\]
where $b_{n,k}$ denotes the sum of the set $\mathcal{M}^{n,k},$ which is
defined by
\[
\mathcal{M}^{n,k}=\left\{  \left(  -1\right)  ^{j-1}a_{i_{1},i_{2}}%
a_{i_{2},i_{3}}\cdots a_{i_{j-1},i_{j}}:i_{1}=n,i_{j}=k,i_{p}>i_{q}\text{ if
}p<q,i_{p}\in\mathbb{N}\ \text{for all }p\right\}  ,
\]
and $\mathcal{M}^{n,n}=\left\{  1\right\}  .$
\end{proposition}

\begin{proof}
Proof is included in Appendix \ref{AppendixYtoH}.
\end{proof}

\begin{theorem}
\label{TheoremH} Let $G_{t}$ be a L\'{e}vy process with no Brownian part
satisfying condition (\ref{condition}). \ Then the power of its increment in
terms of stochastic integrals with respect to the orthogonal martingales, $H,
$ is given by the following equation:%
\begin{equation}
\left(  G_{t+t_{0}}-G_{t_{0}}\right)  ^{k}=\sum_{\theta_{k}\in\mathcal{I}_{k}%
}\Pi_{\theta_{k},t}^{\left(  k\right)  }\mathcal{S}_{\theta_{k},t,t_{0}%
}^{\left(  H\right)  }+C_{t}^{\left(  k\right)  },
\end{equation}
where $\mathcal{I}_{k}$ is defined in (\ref{setI}), $\Pi_{\theta_{k}%
,t}^{\left(  k\right)  }$ is defined in Proposition \ref{Theorem2},
$C_{t}^{\left(  k\right)  }$ is defined in Proposition \ref{TheoremC2} and
$\mathcal{S}_{\left(  i_{1},i_{2},...,i_{j}\right)  ,t,t_{0}}^{\left(
H\right)  } $ is defined as the integral:%
\begin{align*}
& \mathcal{S}_{\left(  i_{1},i_{2},...,i_{j}\right)  ,t,t_{0}}^{\left(
H\right)  }\\
& \ =\sum_{k_{1}=1}^{i_{1}}\cdots\sum_{k_{j-1}=1}^{i_{j-1}}\sum_{k_{j}%
=1}^{i_{j}}b_{i_{1},k_{1}}\cdots b_{i_{j-1},k_{j-1}}b_{i_{j},k_{j}}\int
_{t_{0}}^{t+t_{0}}\int_{t_{0}}^{t_{1}-}\cdots\int_{t_{0}}^{t_{j-1}-}%
\mathrm{d}H_{t_{j}}^{\left(  k_{1}\right)  }\cdots\mathrm{d}H_{t_{2}}^{\left(
k_{j-1}\right)  }\mathrm{d}H_{t_{1}}^{\left(  k_{j}\right)  },
\end{align*}
$b_{n,k}$ is defined in Proposition \ref{PropositionYtoH}.
\end{theorem}

\begin{proof}
From Proposition \ref{PropositionYtoH}, we have%
\begin{align*}
& \mathcal{S}_{\left(  i_{1},i_{2},...,i_{j}\right)  ,t,t_{0}}\\
& ~=\int_{t_{0}}^{t+t_{0}}\int_{t_{0}}^{t_{1}-}\cdots\int_{t_{0}}^{t_{j-1}%
-}\mathrm{d}\hat{G}_{t_{j}}^{\left(  i_{1}\right)  }\cdots\mathrm{d}\hat
{G}_{t_{2}}^{\left(  i_{j-1}\right)  }\mathrm{d}\hat{G}_{t_{1}}^{\left(
i_{j}\right)  }\\
& ~=\int_{t_{0}}^{t+t_{0}}\int_{t_{0}}^{t_{1}-}\cdots\int_{t_{0}}^{t_{j-1}%
-}\mathrm{d}\left[  \sum_{k_{1}=1}^{i_{1}}b_{i_{1},k_{1}}H_{t_{j}}^{\left(
k_{1}\right)  }\right]  \cdots\mathrm{d}\left[  \sum_{k_{j-1}=1}^{i_{j-1}%
}b_{i_{j-1},k_{j-1}}H_{t_{2}}^{\left(  k_{j-1}\right)  }\right]
\mathrm{d}\left[  \sum_{k_{j}=1}^{i_{j}}b_{i_{j},k_{j}}H_{t_{1}}^{\left(
k_{j}\right)  }\right] \\
& ~=\sum_{k_{1}=1}^{i_{1}}\cdots\sum_{k_{j-1}=1}^{i_{j-1}}\sum_{k_{j}%
=1}^{i_{j}}b_{i_{1},k_{1}}\cdots b_{i_{j-1},k_{j-1}}b_{i_{j},k_{j}}\int
_{t_{0}}^{t+t_{0}}\int_{t_{0}}^{t_{1}-}\cdots\int_{t_{0}}^{t_{j-1}-}%
\mathrm{d}H_{t_{j}}^{\left(  k_{1}\right)  }\cdots\mathrm{d}H_{t_{2}}^{\left(
k_{j-1}\right)  }\mathrm{d}H_{t_{1}}^{\left(  k_{j}\right)  }.
\end{align*}
Hence, by using Theorem \ref{Theorem3}, we finish the proof.
\end{proof}

\begin{corollary}
By Theorem \ref{Theorem3},
\begin{align*}
\left(  G_{t+t_{0}}-G_{t_{0}}\right)  ^{m}\left(  G_{t+t_{0}}-G_{t_{0}%
}\right)  ^{n}  & =\left(  \sum_{\theta_{m}\in\mathcal{I}_{m}}\Pi_{\theta
_{m},t}^{\left(  m\right)  }\mathcal{S}_{\theta_{m},t,t_{0}}^{\left(
H\right)  }+C_{t}^{\left(  m\right)  }\right)  \left(  \sum_{\theta_{n}%
\in\mathcal{I}_{n}}\Pi_{\theta_{n},t}^{\left(  n\right)  }\mathcal{S}%
_{\theta_{n},t,t_{0}}^{\left(  H\right)  }+C_{t}^{\left(  n\right)  }\right)
\\
& =\sum_{\theta_{m+n}\in\mathcal{I}_{m+n}}\Pi_{\theta_{m+n},t}^{\left(
m+n\right)  }\mathcal{S}_{\theta_{m+n},t,t_{0}}^{\left(  H\right)  }%
+C_{t}^{\left(  m+n\right)  }.
\end{align*}

\end{corollary}

Hence, we can find out how to express the product of two iterative stochastic
integrals of orders $m$ and $n$ as a weighted sum of iterative stochastic
integrals of order $m+n,\ m+n-1,...,2,1$.

Note in Theorems \ref{Theorem3} and \ref{TheoremH}, the integrands of the
stochastic integrals do \textbf{not} involve $t_{0}$ nor any of the
integrating variables $t_{1},t_{2},...,t_{j}.$ \ They are completely
characterized by $C_{t}^{\left(  p\right)  }$'s$,$ where $p=0,1,...,k.$
\ Hence to find the chaotic representation of $\left(  G_{t+t_{0}}-G_{t_{0}%
}\right)  ^{k},$ we only need to know the moments of $G_{t}$, $m_{1}t=E\left[
X_{t}\right]  $ and $m_{p}=\int_{-\infty}^{\infty}x^{p}\nu\left(
\mathrm{d}x\right)  $ for $p=2,...,k$. \ This result is intuitive as
$(G_{t+t_{0}}-G_{t_{0}})$ is a stationary process.

\subsection{General case}

Next we want to derive the formula for the power of the increments of L\'{e}vy
processes with $\sigma\neq0$. \ Recall $X=\left\{  X_{t},t\geq0\right\}  $
denotes a general L\'{e}vy process, $X_{t}^{\left(  i\right)  }$ denotes its
$i$-th power jump process and $Y_{t}^{\left(  i\right)  }$ denotes its $i$-th
compensated power jump process as defined in (\ref{CompPower}). \ We define
$A_{1}\left(  X_{t+t_{0}},X_{t_{0}};k\right)  $ and $A_{2}\left(  X_{t+t_{0}%
},X_{t_{0}};k\right)  $ such that $\left(  X_{t+t_{0}}-X_{t_{0}}\right)
^{k}=A_{1}\left(  X_{t+t_{0}},X_{t_{0}};k\right)  +A_{2}\left(  X_{t+t_{0}%
},X_{t_{0}};k\right)  ,$ where $A_{1}\left(  X_{t+t_{0}},X_{t_{0}};k\right)  $
comprises all the terms not containing $\sigma$ in $\left(  X_{t+t_{0}%
}-X_{t_{0}}\right)  ^{k}.$ \ \ From (\ref{Corrected1})-(\ref{Corrected}), it
may be noted:
\begin{align}
\left(  X_{t+t_{0}}-X_{t_{0}}\right)  ^{k}  & =\frac{\sigma^{2}}{2}k\left(
k-1\right)  \left(  \left(  X_{t+t_{0}}-X_{t_{0}}\right)  ^{k-2}t-\int_{t_{0}%
}^{t+t_{0}}\left(  s-t_{0}\right)  \mathrm{d}\left(  X_{s}-X_{t_{0}}\right)
^{k-2}\right) \nonumber\\
& +\sum_{j=1}^{k}\binom{k}{j}\int_{t_{0}}^{t+t_{0}}A_{2}\left(  X_{s-}%
,X_{t_{0}};k-j\right)  \mathrm{d}Y_{s}^{\left(  j\right)  }\nonumber\\
& +\sum_{j=1}^{k-1}\binom{k}{j}m_{j}\left(  t+t_{0}\right)  A_{2}\left(
X_{t+t_{0}},X_{t_{0}};k-j\right) \nonumber\\
& -\sum_{j=1}^{k-1}\binom{k}{j}m_{j}\int_{t_{0}}^{t+t_{0}}s\ \mathrm{d}\left[
A_{2}\left(  X_{s},X_{t_{0}};k-j\right)  \right]  +A_{1}\left(  X_{t+t_{0}%
},X_{t_{0}};k\right)  .\label{ito}%
\end{align}

\begin{proposition}
\label{Theorem4}For any L\'{e}vy process $X_{t}$ satisfying condition
(\ref{condition}),
\[
\left(  X_{t+t_{0}}-X_{t_{0}}\right)  ^{k}=A_{1}\left(  X_{t+t_{0}},X_{t_{0}%
};k\right)  +\sum_{n=1}^{\left\lfloor k/2\right\rfloor }\frac{k!}{\left(
k-2n\right)  !}\frac{1}{n!}\frac{1}{2^{n}}\sigma^{2n}A_{1}\left(  X_{t+t_{0}%
},X_{t_{0}};k-2n\right)  t^{n}.
\]

\end{proposition}

\begin{proof}
The proof uses the same techniques as in the proof of Theorem \ref{Theorem3}.
\ Note that $A_{1}\left(  X_{t+t_{0}},X_{t_{0}};p\right)  ,$ where
$p=1,2,...,k,$ are given by Theorem \ref{Theorem3}. \ \ 
\end{proof}

Proposition \ref{Theorem4} gives the formula of $\left(  X_{t+t_{0}}-X_{t_{0}%
}\right)  ^{k}$ in terms of a summation of $A_{1}$, where $\left\lfloor
k/2\right\rfloor +1$ calculations of $A_{1}$ are needed. \ \ The next theorem
gives the formula in an alternative form which requires $A_{1}\,$to be
computed once only. \ 

\begin{definition}
\label{DefinitionCPI}Let $C_{t,\sigma}^{\left(  k\right)  }$ be the terms
obtained by replacing $m_{2}$ with $m_{2}+\sigma^{2}$ in $C_{t}^{\left(
k\right)  }$ (Proposition \ref{TheoremC2}) and $\Pi_{\left(  i_{1}%
,i_{2},...,i_{j}\right)  ,t,\sigma}^{\left(  k\right)  }$ be the terms
obtained by replacing $C_{t}^{\left(  k\right)  }$ with $C_{t,\sigma}^{\left(
k\right)  }$ in $\Pi_{\left(  i_{1},i_{2},...,i_{j}\right)  ,t}^{\left(
k\right)  }$ (Proposition \ref{Theorem2})$.$ \ We then note the following theorem.
\end{definition}

\begin{theorem}
\label{newFormula}For any L\'{e}vy process $X_{t}$ with $\sigma^{2}\neq0$ and
satisfying condition (\ref{condition}), the representation of $\left(
X_{t+t_{0}}-X_{t_{0}}\right)  ^{n}$ is given by Theorem \ref{Theorem3} with
$m_{2}$ replaced by $\left(  m_{2}+\sigma^{2}\right)  ,$ i.e.%
\[
\left(  X_{t+t_{0}}-X_{t_{0}}\right)  ^{n}=\sum_{\theta_{n}\in\mathcal{I}_{n}%
}\Pi_{\theta_{n},t,\sigma}^{\left(  n\right)  }\mathcal{S}_{\theta_{n}%
,t,t_{0}}^{\prime}+C_{t,\sigma}^{\left(  n\right)  },
\]
where $\mathcal{I}_{n}$ is defined in (\ref{setI}), $\Pi_{\theta_{n},t,\sigma
}^{\left(  n\right)  }$ and $C_{t,\sigma}^{\left(  n\right)  }$ are defined
above and $\mathcal{S}_{\left(  i_{1},i_{2},...,i_{j}\right)  ,t,t_{0}%
}^{\prime}$ is defined to be the integral:%
\[
\mathcal{S}_{\left(  i_{1},i_{2},...,i_{j}\right)  ,t,t_{0}}^{\prime}%
=\int_{t_{0}}^{t+t_{0}}\int_{t_{0}}^{t_{1}-}\cdots\int_{t_{0}}^{t_{j-1}%
-}\mathrm{d}Y_{t_{j}}^{\left(  i_{1}\right)  }\cdots\mathrm{d}Y_{t_{2}%
}^{\left(  i_{j-1}\right)  }\mathrm{d}Y_{t_{1}}^{\left(  i_{j}\right)  }.
\]

\end{theorem}

\begin{proof}
We define a new class of power jump processes by%
\begin{align}
\widetilde{X}_{t}^{\left(  2\right)  }  & =X_{t}^{\left(  2\right)  }%
+\sigma^{2}t\text{,}\nonumber\\
\text{ }\widetilde{X}_{t}^{\left(  j\right)  }  & =X_{t}^{\left(  j\right)
}\text{ \ \ for }j=1\text{ and }j=3,4,5,....\label{newX}%
\end{align}
We also define a new class of compensators%
\begin{align*}
\tilde{m}_{2}t  & =\left(  m_{2}+\sigma^{2}\right)  t\text{,}\\
\text{ }\tilde{m}_{j}t  & =m_{j}t\text{ \ \ for }j=1\text{ and }j=3,4,5,....
\end{align*}
Hence, by definition, the compensated power jump processes, $\tilde{Y}%
_{t}^{\left(  i\right)  }=\tilde{X}_{t}^{\left(  i\right)  }-\tilde{m}%
_{i}t=X_{t}^{\left(  i\right)  }-m_{i}t$ $=Y_{t}^{\left(  i\right)  }$ for all
$i\geq1.$ \ Therefore, the representation of $\left(  X_{t+t_{0}}-X_{t_{0}%
}\right)  ^{k}$ in terms of the stochastic integrals with respect to
$Y_{t}^{\left(  i\right)  }$ is the same no matter we start from using
$X_{t}^{\left(  i\right)  }$ or $\widetilde{X}_{t}^{\left(  i\right)  }.$ \ To
calculate the expression using $\widetilde{X}_{t}^{\left(  i\right)  }$, we
use equation (2) in \cite{ns00}:
\begin{align*}
& \left(  X_{t+t_{0}}-X_{t_{0}}\right)  ^{k}\\
& ~=\sum_{j=1}^{k}\binom{k}{j}\int_{t_{0}}^{t+t_{0}}\left(  X_{s-}-X_{t_{0}%
}\right)  ^{k-j}\mathrm{d}X_{s}^{\left(  j\right)  }\\
& ~~~~~~+\frac{\sigma^{2}}{2}k\left(  k-1\right)  \left(  \left(  X_{t+t_{0}%
}-X_{t_{0}}\right)  ^{k-2}t-\int_{0}^{t}s\ \mathrm{d}\left(  X_{s+t_{0}%
}-X_{t_{0}}\right)  ^{k-2}\right) \\
& ~=\sum_{j=1}^{k}\binom{k}{j}\int_{t_{0}}^{t+t_{0}}\left(  X_{s-}-X_{t_{0}%
}\right)  ^{k-j}\mathrm{d}X_{s}^{\left(  j\right)  }+\frac{\sigma^{2}}%
{2}k\left(  k-1\right)  \int_{0}^{t}\left(  X_{\left(  s+t_{0}\right)
-}-X_{t_{0}}\right)  ^{k-2}\mathrm{d}s\\
& ~=\sum_{j=1}^{k}\binom{k}{j}\int_{t_{0}}^{t+t_{0}}\left(  X_{s-}-X_{t_{0}%
}\right)  ^{k-j}\mathrm{d}X_{s}^{\left(  j\right)  }+\frac{\sigma^{2}}%
{2}k\left(  k-1\right)  \int_{t_{0}}^{t+t_{0}}\left(  X_{u-}-X_{t_{0}}\right)
^{k-2}\mathrm{d}u\\
& ~=\sum_{j=1}^{k}\binom{k}{j}\int_{t_{0}}^{t+t_{0}}\left(  X_{s-}-X_{t_{0}%
}\right)  ^{k-j}\mathrm{d}X_{s}^{\left(  j\right)  }+\binom{k}{2}\int_{t_{0}%
}^{t+t_{0}}\left(  X_{s-}-X_{t_{0}}\right)  ^{k-2}\mathrm{d}\left(  \sigma
^{2}s\right)  .
\end{align*}
By (\ref{newX}), we have%
\[
\left(  X_{t+t_{0}}-X_{t_{0}}\right)  ^{k}=\sum_{j=1}^{k}\binom{k}{j}%
\int_{t_{0}}^{t+t_{0}}\left(  X_{s-}-X_{t_{0}}\right)  ^{k-j}\mathrm{d}%
\tilde{X}_{s}^{\left(  j\right)  }.
\]
Using exactly the same calculation as the one leading to (\ref{Corrected1}%
)-(\ref{Corrected}), we have%
\begin{align*}
\left(  X_{t+t_{0}}-X_{t_{0}}\right)  ^{k}  & =\sum_{j=1}^{k}\binom{k}{j}%
\int_{t_{0}}^{t+t_{0}}\left(  X_{s-}-X_{t_{0}}\right)  ^{k-j}\mathrm{d}%
Y_{s}^{\left(  j\right)  }+\sum_{j=1}^{k-1}\binom{k}{j}\tilde{m}_{j}\left(
t+t_{0}\right)  \left(  X_{t+t_{0}}-X_{t_{0}}\right)  ^{k-j}\\
& -\sum_{j=1}^{k-1}\binom{k}{j}\tilde{m}_{j}\int_{t_{0}}^{t+t_{0}%
}s\ \mathrm{d}\left(  X_{s}-X_{t_{0}}\right)  ^{k-j}+\tilde{m}_{k}t.
\end{align*}
This is exactly the equation (\ref{Corrected2})-(\ref{Corrected}) we based on
in the derivation of Theorem \ref{Theorem3}, except that $m_{j}$ is replaced
by $\tilde{m}_{j}.$ \ Hence we now have a simple formula for the
representation of $\left(  X_{t+t_{0}}-X_{t_{0}}\right)  ^{k}$ in terms of the
stochastic integrals with respect to $Y_{t}^{\left(  i\right)  }$ by replacing
$m_{j}$ with $\tilde{m}_{j}$ in the formula given by Theorem \ref{Theorem3}.
\ In other words, we have%
\[
\left(  X_{t+t_{0}}-X_{t_{0}}\right)  ^{n}=\sum_{\theta_{n}\in\mathcal{I}_{n}%
}\Pi_{\theta_{n},t,\sigma}^{\left(  n\right)  }\mathcal{S}_{\theta_{n}%
,t,t_{0}}^{\prime}+C_{t,\sigma}^{\left(  n\right)  },
\]
where $\Pi_{\theta_{n},t,\sigma}^{\left(  n\right)  }$ and $C_{t,\sigma
}^{\left(  n\right)  }$ are defined above. \ Note that this representation
does not depend on the power jump processes directly since it is in terms of
the compensated power jump processes, $Y_{t}^{\left(  j\right)  }.$ \ So it
does not matter if we change the definition of the power jump processes, as
long as we change the compensators accordingly, we will get the same
compensated power jump processes. \ 
\end{proof}

\begin{theorem}
\label{newFormulaH}For any L\'{e}vy process $X_{t}$ with $\sigma^{2}\neq0$ and
satisfying condition (\ref{condition}), the representation of $\left(
X_{t+t_{0}}-X_{t_{0}}\right)  ^{n}$ is given by Theorem \ref{TheoremH} with
$m_{2}$ replaced with $\left(  m_{2}+\sigma^{2}\right)  ,$ i.e.%
\[
\left(  X_{t+t_{0}}-X_{t_{0}}\right)  ^{n}=\sum_{\theta_{n}\in\mathcal{I}_{n}%
}\Pi_{\theta_{n},t,\sigma}^{\left(  n\right)  }\mathcal{S}_{\theta_{n}%
,t,t_{0}}^{\prime\left(  H\right)  }+C_{t,\sigma}^{\left(  n\right)  },
\]
where $\mathcal{I}_{n}$ is defined in (\ref{setI}), $\Pi_{\theta_{n},t,\sigma
}^{\left(  n\right)  }$ and $C_{t,\sigma}^{\left(  n\right)  }$ are defined
above and $\mathcal{S}_{\left(  i_{1},i_{2},...,i_{j}\right)  ,t,t_{0}%
}^{\prime}$ is defined to be the integral:%
\begin{align*}
& \mathcal{S}_{\left(  i_{1},i_{2},...,i_{j}\right)  ,t,t_{0}}^{\prime\left(
H\right)  }\\
& ~=\sum_{k_{1}=1}^{i_{1}}\cdots\sum_{k_{j-1}=1}^{i_{j-1}}\sum_{k_{j}%
=1}^{i_{j}}b_{i_{1},k_{1}}\cdots b_{i_{j-1},k_{j-1}}b_{i_{j},k_{j}}\int
_{t_{0}}^{t+t_{0}}\int_{t_{0}}^{t_{1}-}\cdots\int_{t_{0}}^{t_{j-1}-}%
\mathrm{d}H_{t_{j}}^{\left(  k_{1}\right)  }\cdots\mathrm{d}H_{t_{2}}^{\left(
k_{j-1}\right)  }\mathrm{d}H_{t_{1}}^{\left(  k_{j}\right)  },
\end{align*}
$b_{n,k}$ is defined in Proposition \ref{PropositionYtoH}.
\end{theorem}

\begin{proof}
It follows directly from Theorems \ref{TheoremH} and \ref{newFormula}.\ 
\end{proof}

\noindent\textbf{Remark}\quad\label{RemarkJam}As noted in Section
\ref{SectionJam}, \cite{j05} derived an explicit formula for the chaotic
representation of $X_{t}^{k}$ in terms of the non-compensated power jump
processes, $X_{t}^{\left(  j\right)  }$, when $X_{t}$ is a semimartingale.
\ Our explicit formula gives the representation in terms of orthogonalized
compensated power jump processes, $H_{t}^{\left(  j\right)  }$. \ In the
following, we show that in the L\'{e}vy case, our formula complements
Jamshidian's one. \ We note the notation used by Jamshidian in Section
\ref{SectionJam}.\ \ If $X$ is a L\'{e}vy process, we can see that $\left[
X^{c}\right]  _{t}=\left[  X\right]  _{t}^{c}=\sigma^{2}t$ (where the
superscript $c $ stands for continuous part of the process) and hence $\left[
X\right]  _{t}^{\left(  2\right)  }=\sigma^{2}t+\sum_{s\leq t}\left(  \Delta
X_{s}\right)  ^{2}.$ \ With Jamshidian's notation, the $\sigma^{2}$ is
implicitly included in the $\left[  X\right]  _{t}^{\left(  2\right)  }$.
\smallskip\ \newline\cite{j05} defined $\mathcal{C}=\mathcal{C}^{\ast}%
\cap\mathcal{C}_{\ast}$, where $\mathcal{C}^{\ast}$ is the set of
semimartingales of finite moments with continuous compensators adapted to a
Brownian filtration, and $\mathcal{C}_{\ast}$ is the set of processes with
exponentially decreasing law. \ Jamshidian generalized the CRP from L\'{e}vy
processes to the set $\mathcal{C}.$ \ In proposition 8.2 of \cite{j05}, an
explicit formula for the chaotic representation with respect to the
non-compensated power jump processes for the semimartingales in $\mathcal{C}%
$\ when $t_{0}=0$ was derived. \ \cite{j05} defined the power jump processes
using the power brackets, see (\ref{JpowerJump}) and (\ref{Jcomp}). \ The
multi-indices were denoted by $I=\left(  i_{1},...,i_{p}\right)  \in
\mathbb{N}^{p},$ where $\mathbb{N}$ is the set of natural numbers, and for
integers $1\leq p\leq n,$
\begin{equation}
\mathbb{N}_{n}^{p}=\left\{  I=\left(  i_{1},...,i_{p}\right)  \in
\mathbb{N}^{p}:i_{1}+\cdots+i_{p}=n\right\}  ,\ \ \ p,n\in\mathbb{N}%
.\label{setN}%
\end{equation}
Note that from (\ref{setI}), $\mathcal{I}_{k}=\bigcup\limits_{n=1}^{k}%
\bigcup\limits_{p=1}^{n}\mathbb{N}_{n}^{p}.$ \ Proposition 8.2 of \cite{j05}
states that, for a semimartingale $X_{t}$ with $X_{0}=0,$ we have, for all
$n\in N$\
\begin{equation}
X_{t}^{n}=\sum_{p=1}^{n}\sum_{I\in\mathbb{N}_{n}^{p}}\frac{n!}{i_{1}!\cdots
i_{p}!}\int_{0}^{t}\int_{0}^{t_{1}-}\cdots\int_{0}^{t_{p-1}-}\mathrm{d}\left[
X\right]  _{t_{p}}^{\left(  i_{1}\right)  }\cdots\mathrm{d}\left[  X\right]
_{t_{2}}^{\left(  i_{p}-1\right)  }\mathrm{d}\left[  X\right]  _{t_{1}%
}^{\left(  i_{p}\right)  }.\label{j}%
\end{equation}
Since \cite{j05} only considered non-compensated processes, we substitute all
the $m_{j}$ in (\ref{C}) by zeros (since the compensators in the L\'{e}vy case
are $m_{j}t$), which makes $C_{t}^{\left(  k\right)  }=0$ for all $k\neq0.$
\ So $\Pi_{(i_{1},i_{2},...,i_{j}),t}^{\left(  k\right)  }$ is non-zero only
when $\sum_{p=1}^{j}i_{p}=k$, as defined in (\ref{setN}). \ Hence in the
L\'{e}vy case, Theorem \ref{newFormula} reduces to (\ref{j}).

\begin{corollary}
\label{EX}The expectation of $\left(  X_{t+t_{0}}-X_{t_{0}}\right)  ^{k}$ is
given by $C_{t,\sigma}^{\left(  n\right)  },$ which can be obtained by
replacing $m_{2}$ with $m_{2}+\sigma^{2}$ in $C_{t}^{\left(  k\right)  },$
given by equation (\ref{C2_crp}).
\end{corollary}

\begin{proof}
As the expectations of all the stochastic integrals are zero, this follows
directly from Theorem \ref{newFormula}. \ \ 
\end{proof}

\begin{corollary}
The expectation of $\left[  H_{t}^{\left(  1\right)  }\right]  ^{k}=\left[
\int_{0}^{t}\mathrm{d}H_{t_{1}}^{\left(  1\right)  }\right]  ^{k}$ can be
obtained by replacing $m_{2}$ with $m_{2}+\sigma^{2}$ and $m_{1}$ with $0$ in
$C_{t}^{\left(  k\right)  },$ given by Proposition \ref{TheoremC2}.
\end{corollary}

From Corollary \ref{EX}, $E\left[  X_{t}^{k}\right]  $ can be obtained by
replacing $m_{2}$ with $m_{2}+\sigma^{2}$ in $C_{t}^{\left(  k\right)  }.$
\ Since $H_{t}^{\left(  1\right)  }=X_{t}-m_{1}t$ and
\begin{equation}
\left[  X_{t}\right]  ^{k}=\left[  \int_{0}^{t}\mathrm{d}H_{t_{1}}^{\left(
1\right)  }+m_{1}t\right]  ^{k},\label{XH}%
\end{equation}
by putting $m_{1}=0$ in (\ref{XH}), we can conclude that the expectation of
$\left[  \int_{0}^{t}\mathrm{d}H_{t_{1}}^{\left(  1\right)  }\right]  ^{k}$
can be obtained by replacing $m_{2}$ with $m_{2}+\sigma^{2}$ and $m_{1}$ with
$0 $ in $C_{t}^{\left(  k\right)  }.$\newline

In the next section, we extend our results to chaos expansions in terms of the
Poisson random measure, with the use of the relationship between the two chaos
expansions derived by \cite{bdlop03}.

\section{The Chaos Expansion with respect to the Poisson random
measure\label{Sectionprm}}

\cite{i56} proved a chaos expansion for general L\'{e}vy processes in terms of
multiple integrals with respect to the compensated Poisson random measure.
Note that it is trivial to covert the representation to iterated integrals as
done by \cite{l04}. \ The \textit{compensated Poisson measure }is defined to
be $\tilde{N}\left(  \mathrm{d}t,\mathrm{d}x\right)  =N\left(  \mathrm{d}%
t,\mathrm{d}x\right)  -\nu\left(  \mathrm{d}x\right)  \mathrm{d}t,$ where
$\nu\left(  \mathrm{d}x\right)  $ is the L\'{e}vy measure of the underlying
L\'{e}vy process, $X$, and%
\[
N\left(  B\right)  =\#\left\{  t:\left(  t,\Delta X_{t}\right)  \in B\right\}
,\ \ \ B\in\mathcal{B}\left(  \left[  0,T\right]  \times\mathbb{R}_{0}\right)
,
\]
where $\mathcal{B}\left(  \left[  0,T\right]  \times\mathbb{R}_{0}\right)  $
is the \textit{Borel} $\sigma$\textit{-algebra} of $\left[  0,T\right]
\times\mathbb{R}_{0}$ and $\mathbb{R}_{0}=\mathbb{R-}\left\{  0\right\}  $, is
the jump measure of the process and its compensator is known to be%
\[
E\left[  N\left(  \mathrm{d}t,\mathrm{d}x\right)  \right]  =\nu\left(
\mathrm{d}x\right)  \mathrm{d}t.
\]
Let $f$ be a real function on $\left(  \left[  0,T\right]  \times
\mathbb{R}\right)  ^{n}$. \ Its \textit{symmetrization} $\tilde{f}$ with
respect to the variables $\left(  t_{1},x_{1}\right)  ,...,\left(  t_{n}%
,x_{n}\right)  $ is defined by%
\begin{equation}
\tilde{f}\left(  t_{1},x_{1},...,t_{n},x_{n}\right)  =\frac{1}{n!}%
\sum_{\mathbf{\pi}}f\left(  t_{\pi_{1}},x_{\pi_{1}},...,t_{\pi_{n}},x_{\pi
_{n}}\right)  ,\label{symmetrization}%
\end{equation}
where the sum is taken over all permutations $\mathbf{\pi}$ of $\left\{
1,...,n\right\}  .$ \ $f$ is said to be \textit{symmetric }if $f=\tilde{f}$.
\ The definition of symmetrization is used to represent the CRP in terms of
multiple integrals instead of iterative integrals.

\subsection{Pure jump case}

We first consider the representation for pure jump L\'{e}vy processes as in
\cite{l04}. \ Let $\tilde{L}_{2}\left(  \left(  \lambda\times\nu\right)
^{n}\right)  $ be the space of all square integrable symmetric functions on
$\left(  \left[  0,T\right]  \times\mathbb{R}\right)  ^{n}.$ \ In an iterative
integral such as (\ref{crp}), the time variables $t_{1},...,t_{n}$ are
monotonic. \ For ease of notation, we let%
\begin{equation}
G_{n}=\left\{  \left(  t_{1},x_{1},...,t_{n},x_{n}\right)  :0\leq t\leq
\cdots\leq t_{n}\leq T;x_{i}\in\mathbb{R},i=1,...,n\right\}  ,\label{gn}%
\end{equation}
and let $L_{2}\left(  G_{n}\right)  $ be the space of functions $g$ such that%
\[
\left\Vert g\right\Vert _{L_{2}\left(  G_{n}\right)  }^{2}=\int_{G_{n}}%
g^{2}\left(  t_{1},x_{1},...,t_{n},x_{n}\right)  \mathrm{d}t_{1}\nu\left(
\mathrm{d}x_{1}\right)  \cdots\mathrm{d}t_{n}\nu\left(  \mathrm{d}%
x_{n}\right)  <\infty.
\]
For $f\in L_{2}\left(  G_{n}\right)  ,$ let%
\[
J_{n}\left(  f\right)  =\int_{0}^{T}\int_{\mathbb{R}}\cdots\int_{0}^{t_{2}%
}\int_{\mathbb{R}}f\left(  t_{1},x_{1},...,t_{n},x_{n}\right)  \tilde
{N}\left(  \mathrm{d}t_{1},\mathrm{d}x_{1}\right)  \cdots\tilde{N}\left(
\mathrm{d}t_{n},\mathrm{d}x_{n}\right)  .
\]
For $f\in\tilde{L}_{2}\left(  \left(  \lambda\times\nu\right)  ^{n}\right)  ,
$ let%
\[
I_{n}\left(  f\right)  =\int_{\left(  \left[  0,T\right]  \times
\mathbb{R}\right)  ^{n}}f\left(  t_{1},x_{1},...,t_{n},x_{n}\right)  \tilde
{N}^{\otimes n}\left(  \mathrm{d}\mathbf{t},\mathrm{d}\mathbf{x}\right)
=n!J_{n}\left(  f\right)  .
\]

\begin{itemize}
\item \textit{Chaos Expansion for pure jump L\'{e}vy processes:}%
\textbf{\ \ }Let $F$ be a square integrable random variable adapted to the
underlying pure jump L\'{e}vy process, $X_{t}$. \ \cite{i56} proved that there
exists a unique sequence $\left\{  f_{n}\right\}  _{n=0}^{\infty}$ where
$f_{n}\in\tilde{L}_{2}\left(  \left[  0,T\right]  \times\mathbb{R}\right)
^{n}$ such that%
\begin{equation}
F=E\left(  F\right)  +\sum_{n=1}^{\infty}I_{n}\left(  f_{n}\right)
.\label{CEPRN}%
\end{equation}

\end{itemize}

\cite{bdlop03} derived relations between the two chaos expansions, that is,
between the expansion in terms of compensated power jump processes and the
expansion in terms of the Poisson random measure\textbf{. \ }\cite{bdlop03}
showed that the compensated power jump process defined in (\ref{CompPower})
satisfies the equation%
\begin{equation}
Y^{\left(  i\right)  }\left(  t\right)  =\int_{0}^{t}\int_{\mathbb{R}}%
x^{i}\tilde{N}\left(  \mathrm{d}s,\mathrm{d}x\right)  ,\text{ \ \ }0\leq t\leq
T,\ i=1,2,....\label{YN}%
\end{equation}
Hence, the CRP can be written as%
\begin{align}
F  & =E(F)+\sum_{j=1}^{\infty}\sum_{i_{1},...,i_{j}\geq1}\int_{0}^{\infty}%
\int_{0}^{{t_{1}}-}\cdots\int_{0}^{{t_{j-1}}-}f_{(i_{1},...,i_{j})}%
(t_{1},...,t_{j})\mathrm{d}Y_{t_{j}}^{(i_{j})}...\mathrm{d}Y_{t_{2}}^{(i_{2}%
)}\mathrm{d}Y_{t_{1}}^{(i_{1})}\label{line1}\\
& =E(F)+\sum_{j=1}^{\infty}\sum_{i_{1},...,i_{j}\geq1}\int_{0}^{\infty}%
\int_{\mathbb{R}}\int_{0}^{{t_{1}}-}\int_{\mathbb{R}}\cdots\int_{0}^{{t_{j-1}%
}-}\int_{\mathbb{R}}x_{j}^{i_{j}}\cdots x_{2}^{i_{2}}x_{1}^{i_{1}}\nonumber\\
& \times f_{(i_{1},...,i_{j})}(t_{1},...,t_{j})\tilde{N}\left(  \mathrm{d}%
t_{j},\mathrm{d}x_{j}\right)  \cdots\tilde{N}\left(  \mathrm{d}t_{2}%
,\mathrm{d}x_{2}\right)  \tilde{N}\left(  \mathrm{d}t_{1},\mathrm{d}%
x_{1}\right) \label{line3}\\
& =E(F)+\sum_{j=1}^{\infty}\int_{0}^{\infty}\int_{\mathbb{R}}\int_{0}^{{t_{1}%
}-}\int_{\mathbb{R}}\cdots\int_{0}^{{t_{j-1}}-}\int_{\mathbb{R}}g_{j}\left(
t_{1},x_{1},...,t_{j},x_{j}\right) \nonumber\\
& \times\tilde{N}\left(  \mathrm{d}t_{j},\mathrm{d}x_{j}\right)  \cdots
\tilde{N}\left(  \mathrm{d}t_{2},\mathrm{d}x_{2}\right)  \tilde{N}\left(
\mathrm{d}t_{1},\mathrm{d}x_{1}\right) \nonumber\\
& =E(F)+\sum_{j=1}^{\infty}J_{j}\left(  g_{j}\right)  =E(F)+\sum_{j=1}%
^{\infty}n!J_{j}\left(  \tilde{g}_{j}\right)  =E(F)+\sum_{j=1}^{\infty}%
I_{j}\left(  \tilde{g}_{j}\right)  ,\nonumber
\end{align}
where $\tilde{g}_{j}$ is the symmetrization (defined in (\ref{symmetrization}%
)) of the function $g_{j}$ given by%
\begin{align}
& g_{j}\left(  t_{1},x_{1},...,t_{j},x_{j}\right) \nonumber\\
& \ =\left\{
\begin{array}
[c]{ll}%
\sum_{i_{1},...,i_{j}\geq1}x_{1}^{i_{1}}\cdots x_{j}^{i_{j}}f_{(i_{1}%
,...,i_{j})}(t_{1},...,t_{j}), & \text{on }G_{n}\\
0 & \text{on }\left(  \left[  0,T\right]  \times\mathbb{R}\right)  ^{n}-G_{n}.
\end{array}
\right. \label{gj_prm}%
\end{align}
Therefore, by uniqueness, $\left\{  f_{n}\right\}  _{n=0}^{\infty}$ in
(\ref{CEPRN}) is given by%
\[
f_{n}=\tilde{g}_{n},\text{ \ \ \ }n=1,2,....
\]
This equation provides a simple relationship between the two expansions.
\ From Theorem \ref{newFormula}, of course,%
\begin{equation}
\left(  X_{t+t_{0}}-X_{t_{0}}\right)  ^{n}=\sum_{\theta_{n}\in\mathcal{I}_{n}%
}\Pi_{\theta_{n},t,\sigma}^{\left(  n\right)  }\mathcal{S}_{\theta_{n}%
,t,t_{0}}^{\prime}+C_{t,\sigma}^{\left(  n\right)  }.\label{newFormula2}%
\end{equation}
We can now use this relationship to derive a form for $\tilde{g}_{n}$ in terms
of $\mathcal{I}_{n},\ \Pi_{\theta_{n},t,\sigma}^{\left(  n\right)  }$ and
$C_{t,\sigma}^{\left(  n\right)  }.$ \ Let $\mathcal{K}_{l,s}=\left\{  \left(
i_{1},...,i_{l}\right)  |i_{j}\in\left\{  1,2,...,s\right\}  \text{ and}%
\sum_{j=1}^{l}i_{j}=s\right\}  .$ \ Since the length of a tuple must not be
greater than the sum of all the elements in the tuple (because an element must
be at least 1), $l\leq s.$ \ By definition, we have $\mathcal{I}%
_{n}=\mathop{\displaystyle \bigcup }\limits_{s=1}^{n}\mathop{\displaystyle
\bigcup }\limits_{l=1}^{s}\mathcal{K}_{l,s}.$ \ So we can write%
\[
\left(  X_{t+t_{0}}-X_{t_{0}}\right)  ^{n}=\sum_{l=1}^{n}\sum_{s=l}^{n}%
\sum_{\theta_{n}\in\mathcal{K}_{l,s}}\Pi_{\theta_{n},t,\sigma}^{\left(
n\right)  }\mathcal{S}_{\theta_{n},t,t_{0}}^{\prime}+C_{t,\sigma}^{\left(
n\right)  },
\]
where $\theta_{n}$ is the tuple $\left(  i_{1}^{\theta_{n}},...,i_{l}%
^{\theta_{n}}\right)  $ with $l$ elements which sum up to $s.$ \ Therefore, we
deduce that for $F=\left(  X_{t+t_{0}}-X_{t_{0}}\right)  ^{n}$ in
(\ref{line1}), $f_{(i_{1},...,i_{j})}(t_{1},...,t_{j})$ is given by
\begin{equation}
f_{(i_{1},...,i_{j})}(t_{1},...,t_{j})=\Pi_{\theta_{n},t,\sigma}^{\left(
n\right)  }.\label{fPI}%
\end{equation}
By (\ref{gj_prm}), we have then proved the following proposition.

\begin{proposition}
For any pure jump L\'{e}vy process $X_{t}$ satisfying condition
(\ref{condition}),
\[
\left(  X_{t+t_{0}}-X_{t_{0}}\right)  ^{n}=\sum_{l=1}^{n}I_{l}\left(
\tilde{g}_{l}^{\left(  n\right)  }\right)  +C_{t,\sigma}^{\left(  n\right)  },
\]
where $\tilde{g}_{l}^{\left(  n\right)  }$ is the symmetrization of the
function $g_{l}^{\left(  n\right)  }$ defined by%
\begin{align*}
& g_{l}^{\left(  n\right)  }\left(  t_{1},x_{1},...,t_{l},x_{l}\right) \\
& \ \ \ \ =\left\{
\begin{array}
[c]{ll}%
\sum_{s=l}^{n}\sum_{\theta_{n}\in\mathcal{K}_{l,s}}x_{1}^{i_{1}^{\theta_{n}}%
}\cdots x_{j}^{i_{j}^{\theta_{n}}}\Pi_{\theta_{n},t,\sigma}^{\left(  n\right)
}, & \text{on }G_{n}\\
0 & \text{on }\left(  \left[  0,T\right]  \times\mathbb{R}\right)  ^{n}-G_{n},
\end{array}
\right.
\end{align*}
where $C_{t,\sigma}^{\left(  n\right)  }$ and $\Pi_{\theta_{n},t,\sigma
}^{\left(  n\right)  }$ are defined in Definition \ref{DefinitionCPI}.
\end{proposition}%

\noindent
The following proposition gives a more straightforward representation.

\begin{proposition}
For any pure jump L\'{e}vy process $X_{t}$ satisfying condition
(\ref{condition}),%
\begin{align}
\left(  X_{t+t_{0}}-X_{t_{0}}\right)  ^{n}  & =\sum_{\theta_{n}\in
\mathcal{I}_{n}}\int_{0}^{\infty}\int_{\mathbb{R}}\int_{0}^{{t_{1}}-}%
\int_{\mathbb{R}}\cdots\int_{0}^{{t_{j-1}}-}\int_{\mathbb{R}}x_{j}%
^{i_{j}^{\theta_{n}}}\cdots x_{2}^{i_{2}^{\theta_{n}}}x_{1}^{i_{1}^{\theta
_{n}}}\nonumber\\
& \times\Pi_{\theta_{n},t,\sigma}^{\left(  n\right)  }\tilde{N}\left(
\mathrm{d}t_{j},\mathrm{d}x_{j}\right)  \cdots\tilde{N}\left(  \mathrm{d}%
t_{2},\mathrm{d}x_{2}\right)  \tilde{N}\left(  \mathrm{d}t_{1},\mathrm{d}%
x_{1}\right)  +C_{t,\sigma}^{\left(  n\right)  },\label{prm}%
\end{align}
where $C_{t,\sigma}^{\left(  n\right)  }$ and $\Pi_{\theta_{n},t,\sigma
}^{\left(  n\right)  }$ are defined in Definition \ref{DefinitionCPI}.
\end{proposition}

\begin{proof}
This follows directly by replacing $f_{(i_{1},...,i_{j})}(t_{1},...,t_{j})$ in
(\ref{line3}) by (\ref{fPI}).
\end{proof}

Note that both chaos expansions, that is, the expansion in terms of
compensated power jump processes and the expansion in terms of random measure,
depend on $\mathcal{I}_{n},$ $\Pi_{\theta_{n},t,\sigma}^{\left(  n\right)  }$
and $C_{t,\sigma}^{\left(  n\right)  }$. \ From (\ref{YN}), we note the
relationship between $Y^{\left(  i\right)  }\left(  t\right)  $ and $\tilde
{N}\left(  \mathrm{d}s,\mathrm{d}x\right)  .$ \ Because of the simple form of
this relationship, we can use Theorem \ref{Theorem3} to derive the explicit
representation of (\ref{prm}).

\subsection{General case}

We shall now discuss the general relationship between the two representations.
\ \cite{i56} proved the chaos expansion for general L\'{e}vy functionals.
\ \cite{bdlop03} gave the relationship between the chaos expansions in the
case with both a continuous (Wiener process) component and a pure jump
(Poisson random measure) component. \ In this general case, the stochastic
integrals are in terms of both Brownian motion, $W$, and the compensated
Poisson measure, $\tilde{N}\left(  \cdot,\cdot\right)  $. \ Hence, to unify
notation, \cite{bdlop03} defined:%
\begin{align*}
U_{1}=\left[  0,T\right]  \text{ \ \ }  & \text{and \ \ }U_{2}=\left[
0,T\right]  \times\mathbb{R}\\
\mathrm{d}Q_{1}\left(  \cdot\right)  =\mathrm{d}W\left(  \cdot\right)  \text{
\ \ }  & \text{and \ \ }Q_{2}\left(  \cdot\right)  =\tilde{N}\left(
\cdot,\cdot\right) \\
\mathrm{d}\left\langle Q_{1}\right\rangle =\mathrm{d}\lambda\text{ \ \ }  &
\text{and \ \ }\mathrm{d}\left\langle Q_{2}\right\rangle =\mathrm{d}%
\lambda\times\mathrm{d}\nu\\
\int_{U_{1}}g\left(  u^{\left(  1\right)  }\right)  Q_{1}\left(
\mathrm{d}u^{\left(  1\right)  }\right)   & =\int_{0}^{t}g\left(  s\right)
W\left(  \mathrm{d}s\right)  \text{ \ \ and}\\
\int_{U_{2}}g\left(  u^{\left(  2\right)  }\right)  Q_{2}\left(
\mathrm{d}u^{\left(  2\right)  }\right)   & =\int_{0}^{t}\int_{\mathbb{R}%
}g\left(  s,x\right)  \tilde{N}\left(  \mathrm{d}s,\mathrm{d}x\right)  ,
\end{align*}
where $\lambda\left(  \mathrm{d}t\right)  =\mathrm{d}t$ is the Lebesgue
measure on $\left[  0,T\right]  .$ \ Let $F$ be a square integrable random
variable adapted to the filtration generated by the underlying L\'{e}vy
process, $X$. \ \cite{bdlop03} proved that%
\begin{equation}
F=E\left[  F\right]  +\sum_{n=1}^{\infty}\sum_{j_{1},...,j_{n}=1,2}\tilde
{J}_{n}\left(  g_{n}^{\left(  j_{1},...,j_{n}\right)  }\right) \label{gprm}%
\end{equation}
for a unique sequence $g_{n}^{\left(  j_{1},...,j_{n}\right)  }$
($j_{1},...,j_{n}=1,2;\ n=1,2,...$) of deterministic functions in the
corresponding $L_{2}$-space, $L_{2}\left(  G_{n}\right)  =L_{2}\left(
G_{n},\otimes_{i=1}^{n}d\left\langle Q_{j_{i}}\right\rangle \right)  ,$ where
\[
G_{n}=\left\{  \left(  u_{1}^{\left(  j_{1}\right)  },...,u_{n}^{\left(
j_{n}\right)  }\right)  \in\Pi_{i=1}^{n}U_{j_{i}}:0\leq t_{1}\leq\cdots\leq
t_{n}\leq T\right\}
\]
with $u^{\left(  j_{i}\right)  }=t$ if $j_{i}=1$, and $u^{\left(
j_{i}\right)  }=\left(  t,x\right)  $ if $j_{i}=2,$ and
\begin{align*}
& \tilde{J}_{n}\left(  g_{n}^{\left(  j_{1},...,j_{n}\right)  }\right) \\
& \ =\int_{\Pi_{i=1}^{n}U_{j_{i}}}g_{n}^{\left(  j_{1},...,j_{n}\right)
}\left(  u_{1}^{\left(  j_{1}\right)  },...,u_{n}^{\left(  j_{n}\right)
}\right)  1_{G_{n}}\left(  u_{1}^{\left(  j_{1}\right)  },...,u_{n}^{\left(
j_{n}\right)  }\right)  Q_{j_{1}}\left(  \mathrm{d}u_{1}^{\left(
j_{1}\right)  }\right)  \cdots Q_{j_{n}}\left(  \mathrm{d}u_{n}^{\left(
j_{n}\right)  }\right)  .
\end{align*}

Similar to the pure jump case, we can derive the explicit formula for the
chaos expansion with respect to the Poisson random measure of a general
L\'{e}vy process, i.e. $\sigma\neq0.$ \ In this case, we have
\begin{align*}
Y^{\left(  1\right)  }\left(  t\right)   & =\sigma\int_{0}^{t}\mathrm{d}%
W\left(  ds\right)  +\int_{0}^{t}\int_{\mathbb{R}}x\tilde{N}\left(
\mathrm{d}s,\mathrm{d}x\right) \\
Y^{\left(  i\right)  }\left(  t\right)   & =\int_{0}^{t}\int_{\mathbb{R}}%
x^{i}\tilde{N}\left(  \mathrm{d}s,\mathrm{d}x\right)  ,\text{ \ \ }0\leq t\leq
T,\ i=2,3,....
\end{align*}
To derive the relation between the two chaos expansions, we introduce the
following notation:%
\begin{align*}
R^{\left(  1\right)  }\left(  \mathrm{d}s,\mathrm{d}x\right)   &
=\sigma\mathrm{d}W\left(  \mathrm{d}s\right)  +\int_{\mathbb{R}}x\tilde
{N}\left(  \mathrm{d}s,\mathrm{d}x\right) \\
R^{\left(  i\right)  }\left(  \mathrm{d}s,\mathrm{d}x\right)   &
=\int_{\mathbb{R}}x^{i}\tilde{N}\left(  \mathrm{d}s,\mathrm{d}x\right)
,\text{ \ \ }i=2,3,....
\end{align*}
Hence the CRP with respect to the power jump processes can be written as%
\begin{align*}
F  & =E(F)+\sum_{j=1}^{\infty}\sum_{i_{1},...,i_{j}\geq1}\int_{0}^{\infty}%
\int_{0}^{{t_{1}}-}\cdots\int_{0}^{{t_{j-1}}-}f_{(i_{1},...,i_{j})}%
(t_{1},...,t_{j})\mathrm{d}Y_{t_{j}}^{(i_{j})}...\mathrm{d}Y_{t_{2}}^{(i_{2}%
)}\mathrm{d}Y_{t_{1}}^{(i_{1})}\\
& =E(F)+\sum_{j=1}^{\infty}\sum_{i_{1},...,i_{j}\geq1}\int_{0}^{\infty}%
\int_{0}^{{t_{1}}-}\cdots\int_{0}^{{t_{j-1}}-}f_{(i_{1},...,i_{j})}%
(t_{1},...,t_{j})\\
& \times R^{(i_{j})}\left(  \mathrm{d}t_{j},\mathrm{d}x_{j}\right)
...R^{(i_{2})}\left(  \mathrm{d}t_{2},\mathrm{d}x\right)  R^{(i_{1})}\left(
\mathrm{d}t_{1},\mathrm{d}x\right)  .
\end{align*}
From Theorem \ref{newFormula},%
\begin{align*}
\left(  X_{t+t_{0}}-X_{t_{0}}\right)  ^{n}  & =\sum_{\theta_{n}\in
\mathcal{I}_{n}}\Pi_{\theta_{n},t,\sigma}^{\left(  n\right)  }\mathcal{S}%
_{\theta_{n},t,t_{0}}^{\prime}+C_{t,\sigma}^{\left(  n\right)  }\\
& =\sum_{\theta_{n}\in\mathcal{I}_{n}}\int_{0}^{\infty}\int_{0}^{{t_{1}}%
-}\cdots\int_{0}^{{t_{j-1}}-}\Pi_{\theta_{n},t,\sigma}^{\left(  n\right)  }\\
& \times R^{(i_{j}^{\theta_{n}})}\left(  \mathrm{d}t_{j},\mathrm{d}%
x_{j}\right)  ...R^{(i_{2}^{\theta_{n}})}\left(  \mathrm{d}t_{2}%
,\mathrm{d}x\right)  R^{(i_{1}^{\theta_{n}})}\left(  \mathrm{d}t_{1}%
,\mathrm{d}x\right)  +C_{t,\sigma}^{\left(  n\right)  }.
\end{align*}
We have then proved the following proposition.

\begin{proposition}
\label{gefprm}For any L\'{e}vy process $X_{t}$ satisfying condition
(\ref{condition}),%
\begin{align*}
\left(  X_{t+t_{0}}-X_{t_{0}}\right)  ^{n}  & =\sum_{\theta_{n}\in
\mathcal{I}_{n}}\int_{0}^{\infty}\int_{0}^{{t_{1}}-}\cdots\int_{0}^{{t_{j-1}%
}-}\Pi_{\theta_{n},t,\sigma}^{\left(  n\right)  }\\
& \times R^{(i_{j}^{\theta_{n}})}\left(  \mathrm{d}t_{j},\mathrm{d}%
x_{j}\right)  ...R^{(i_{2}^{\theta_{n}})}\left(  \mathrm{d}t_{2}%
,\mathrm{d}x\right)  R^{(i_{1}^{\theta_{n}})}\left(  \mathrm{d}t_{1}%
,\mathrm{d}x\right)  +C_{t,\sigma}^{\left(  n\right)  },
\end{align*}
where $C_{t,\sigma}^{\left(  n\right)  }$ and $\Pi_{\theta_{n},t,\sigma
}^{\left(  n\right)  }$ are defined in Definition \ref{DefinitionCPI}.
\end{proposition}

\section{The explicit chaos expansions for a common kind of L\'{e}vy
functionals\label{SectionGeneral}}

Note that we have only found the explicit representations for powers of
increments of L\'{e}vy processes. \ In this section, we explain how the
explicit formulae for a common kind of L\'{e}vy functionals might be obtained
using multivariate Taylor expansion.

Assume that a real function $g,$ possessing derivatives of all orders, is such
that%
\begin{equation}
F=g\left(  X_{t_{1}},X_{t_{2}}-X_{t_{1}},...,X_{t_{n}}-X_{t_{n-1}}\right)
,\label{g}%
\end{equation}
where the indices $0\leq t_{1}<t_{2}<\cdots<t_{n}$ are known and $n$ is
finite. \ By expressing $F$ in terms of power of increments of $X$, we can use
our explicit formula to obtain the CRP of $F$. \ For example, in financial
applications, $g$ corresponds to all pricing functions of contingent claims
which depend on the underlying asset at a finite number of time points.
\ Suppose $\left\{  X_{t},0\leq t\leq T\right\}  $ is the background driving
L\'{e}vy process and time is now $t=t_{n}$. \ Suppose the underlying asset,
$\left\{  S_{t},0\leq t\leq T\right\}  $, is given by the exponential-L\'{e}vy
model, see \cite[Chapter 8.4]{ct03}:%
\[
S_{t}=S_{0}\exp\left(  X_{t}\right)  ,
\]
where $S_{0}$ is the initial value of the underlying asset at time $t=0$.
$\ $Then, for example, we can represent $F$ as the pricing functions of a
number of contingent claims listed in Table 1.

In (\ref{g}), let $x_{1}=X_{t_{1}},x_{2}=X_{t_{2}}-X_{t_{1}},...,x_{n}%
=X_{t_{n}}-X_{t_{n-1}}$. \ \ If $g$ is not a linear combination of powers of
$x_{i}$, we need to use the multivariate Taylor series, see \cite{jj88}, about
the points $x_{i}=0,i=1,...,n$ to obtain such a representation:%
\begin{equation}
g\left(  x_{1},...,x_{n}\right)  =\sum_{j=0}^{\infty}\left\{  \frac{1}%
{j!}\left[  \sum_{k=1}^{n}x_{k}\frac{\partial}{\partial x_{k}^{\prime}%
}\right]  ^{j}g\left(  x_{1}^{\prime},...,x_{n}^{\prime}\right)  \right\}
_{x_{1}^{\prime}=0,...,x_{n}^{\prime}=0}.\label{MultiTaylor}%
\end{equation}
Note that this representation exists when $g$ possesses derivatives of all
orders at zero. \ To show typical elements in this representation, we note the
special case of $n=2$:%
\begin{align*}
g\left(  x_{1},x_{2}\right)   & =\sum_{j=0}^{\infty}\left\{  \frac{1}%
{j!}\left[  x_{1}\frac{\partial}{\partial x_{1}^{\prime}}+x_{2}\frac{\partial
}{\partial x_{2}^{\prime}}\right]  ^{j}g\left(  x_{1}^{\prime},x_{2}^{\prime
}\right)  \right\}  _{x_{1}^{\prime}=0,x_{2}^{\prime}=0}\\
& =g\left(  0,0\right)  +\left[  x_{1}\left.  \frac{\partial g}{\partial
x_{1}^{\prime}}\right\vert _{x_{1}^{\prime}=0,x_{2}^{\prime}=0}+x_{2}\left.
\frac{\partial g}{\partial x_{2}^{\prime}}\right\vert _{x_{1}^{\prime}%
=0,x_{2}^{\prime}=0}\right] \\
& +\frac{1}{2!}\left[  x_{1}^{2}\left.  \frac{\partial^{2}g}{\partial
x_{1}^{\prime2}}\right\vert _{x_{1}^{\prime}=0,x_{2}^{\prime}=0}+2x_{1}%
x_{2}\left.  \frac{\partial^{2}g}{\partial x_{1}^{\prime}\partial
x_{2}^{\prime}}\right\vert _{x_{1}^{\prime}=0,x_{2}^{\prime}=0}+x_{2}%
^{2}\left.  \frac{\partial^{2}g}{\partial x_{2}^{\prime2}}\right\vert
_{x_{1}^{\prime}=0,x_{2}^{\prime}=0}\right]  +\cdots.
\end{align*}

\begin{center}%
\begin{tabular}
[c]{|ll|}\hline
\textbf{Name} & \multicolumn{1}{|l|}{\textbf{Formula}}\\\hline
\multicolumn{1}{|l|}{Forward and future contracts on a} & $F_{t}=S_{t}%
\exp\left(  r\left(  T-t\right)  \right)  =S_{0}\exp\left(  X_{t}+r\left(
T-t\right)  \right)  ,$\\
\multicolumn{1}{|l|}{security providing no income} & where $r$ is the risk
free interest rate and $T$ is the\\
\multicolumn{1}{|l|}{} & maturity of the contract.\\\hline
\multicolumn{1}{|l|}{Forward and future contracts on a} & $F_{t}=\left(
S_{t}-I\right)  \exp\left(  r\left(  T-t\right)  \right)  $\\
\multicolumn{1}{|l|}{security providing a known cash} & $\ \ =\left(
S_{0}\exp\left(  X_{t}\right)  -I\right)  \exp\left(  r\left(  T-t\right)
\right)  ,$\\
\multicolumn{1}{|l|}{income} & where $I$ is the present value of the
perfectly\\
\multicolumn{1}{|l|}{} & predictable income on $S$.\\\hline
\multicolumn{1}{|l|}{Forward and future contracts on a} & $F_{t}=S_{t}%
\exp\left(  \left(  r-r_{f}\right)  \left(  T-t\right)  \right)  $\\
\multicolumn{1}{|l|}{foreign currency} & $\ \ =S_{0}\exp\left(  X_{t}+\left(
r-r_{f}\right)  \left(  T-t\right)  \right)  ,\ $where $r_{f}$ is the\\
\multicolumn{1}{|l|}{} & risk free interest rate of the foreign
currency.\\\hline
\multicolumn{1}{|l|}{Forward and future contracts on} & $F_{t}=\left(
S_{t}+U\right)  \exp\left(  r\left(  T-t\right)  \right)  $\\
\multicolumn{1}{|l|}{commodity} & $\ \ =\left(  S_{0}\exp\left(  X_{t}\right)
+U\right)  \exp\left(  r\left(  T-t\right)  \right)  ,$\\
\multicolumn{1}{|l|}{} & where $U$ is the present value of all storage
costs.\\\hline
\multicolumn{1}{|l|}{European call options} & $F\left(  t,S_{t}\right)
=\exp\left(  -r\left(  T-t\right)  \right)  E_{Q}\left[  \left(
S_{T}-K\right)  ^{+}|\mathcal{F}_{t}\right]  ,$\\
\multicolumn{1}{|l|}{} & where $K$ is the strike, $T$ is the maturity, $Q$ is
the\\
\multicolumn{1}{|l|}{} & risk neutral measure and $\mathcal{F}_{t}$ is the
filtration of $S$.\\\hline
\multicolumn{1}{|l|}{Up-and-out barrier call options} & $F\left(
t,S_{t}\right)  =\exp\left(  -r\left(  T-t\right)  \right)  E_{Q}\left[
\left(  S_{T}-K\right)  ^{+}1_{\left\{  M_{T}^{S}<H\right\}  }\right]  ,$\\
\multicolumn{1}{|l|}{} & where $H$ is the barrier and\\
\multicolumn{1}{|l|}{} & $M_{t}^{S}=\sup\left\{  S_{u};0\leq u\leq t\right\}
,\text{ }0\leq t\leq T.$\\\hline
\multicolumn{1}{|l|}{Up-and-in barrier call options} & $F\left(
t,S_{t}\right)  =\exp\left(  -r\left(  T-t\right)  \right)  E_{Q}\left[
\left(  S_{T}-K\right)  ^{+}1_{\left\{  M_{T}^{S}\geq H\right\}  }\right]
.$\\\hline
\multicolumn{1}{|l|}{Down-and-out barrier call options} & $F\left(
t,S_{t}\right)  =\exp\left(  -r\left(  T-t\right)  \right)  E_{Q}\left[
\left(  S_{T}-K\right)  ^{+}1_{\left\{  m_{T}^{S}>H\right\}  }\right]  ,$\\
\multicolumn{1}{|l|}{} & where $m_{t}^{S}=\inf\left\{  S_{u};0\leq u\leq
t\right\}  ,\text{ }0\leq t\leq T.$\\\hline
\multicolumn{1}{|l|}{Down-and-in barrier call options} &
\multicolumn{1}{|l|}{$F\left(  t,S_{t}\right)  =\exp\left(  -r\left(
T-t\right)  \right)  E_{Q}\left[  \left(  S_{T}-K\right)  ^{+}1_{\left\{
m_{T}^{S}\leq H\right\}  }\right]  .$}\\\hline
\multicolumn{1}{|l|}{Lookback options with a floating strike} &
\multicolumn{1}{|l|}{$F\left(  t,S_{t}\right)  =\exp\left(  -r\left(
T-t\right)  \right)  E_{Q}\left[  M_{T}^{S}-S_{T}\right]  .$}\\\hline
\multicolumn{1}{|l|}{Lookback options with a fixed strike} &
\multicolumn{1}{|l|}{$F\left(  t,S_{t}\right)  =\exp\left(  -r\left(
T-t\right)  \right)  E_{Q}\left[  \left(  M_{T}^{S}-K\right)  ^{+}\right]  .$%
}\\\hline
\multicolumn{1}{|l|}{Asian call options} & \multicolumn{1}{|l|}{$F\left(
t,S_{t}\right)  =\frac{\exp\left(  -r\left(  T-t\right)  \right)  }{n}%
E_{Q}\left[  \left.  \left(  \sum_{k=1}^{n}S_{t_{k}}-nK\right)  ^{+}%
\right\vert \mathcal{F}_{t}\right]  .$}\\\hline
\end{tabular}
. \ 
\end{center}

Table 1: \ The contingent claims and their pricing formulae to which Taylor's
expansion can be applied at some values of $S_{t}.$

Let $g_{j_{1},j_{2},...,j_{l}}^{\left(  l\right)  }\left(  \mathbf{0}\right)
=\frac{1}{l!}\left.  \frac{\partial^{l}g}{\partial x_{j_{1}}^{\prime}\partial
x_{j_{2}}^{\prime}\cdots\partial x_{j_{l}}^{\prime}}\right\vert _{x_{1}%
^{\prime}=0,...,x_{n}^{\prime}=0}.$ $\ $As in \cite[Lemma 2]{cns05}, we assume
that%
\begin{equation}
\sum_{l=2}^{\infty}\sum_{j_{1},...,j_{l}\in\left\{  1,...,n\right\}
}\left\vert g_{j_{1},j_{2},...,j_{l}}^{\left(  l\right)  }\left(
\mathbf{0}\right)  \right\vert R^{l}<\infty,\label{taylorCon}%
\end{equation}
for all $R>0.$ \ The multivariate Taylor series (\ref{MultiTaylor}) expresses
$F$ in terms of sum of products of powers of increments of $X.$ \ From Theorem
\ref{newFormulaH}, we can substitute $x_{i},\ i=1,2,...$with the iterated
integrals with respect to the orthogonal martingales$.$ \ 

For all $F\in L^{2}(\Omega,\mathcal{F})$ having the form (\ref{g}), then
\begin{align}
F  & =\sum_{j=0}^{\infty}\left\{  \frac{1}{j!}\left[  \sum_{k=1}^{n}\left(
X_{t_{k}}-X_{t_{k-1}}\right)  \frac{\partial}{\partial x_{k}^{\prime}}\right]
^{j}g\left(  x_{1}^{\prime},...,x_{n}^{\prime}\right)  \right\}
_{x_{1}^{\prime}=0,...,x_{n}^{\prime}=0}\nonumber\\
& =g\left(  0,0,...,0\right)  +\sum_{j=1}^{n}\left(  X_{t_{j}}-X_{t_{j-1}%
}\right)  g_{j}^{\left(  1\right)  }\left(  \mathbf{0}\right)  +\frac{1}%
{2!}\sum_{j=1}^{n}\left(  X_{t_{j}}-X_{t_{j-1}}\right)  ^{2}g_{j,j}^{\left(
2\right)  }\left(  \mathbf{0}\right) \nonumber\\
& +\frac{1}{2!}\sum_{j_{1}=1}^{n}\sum_{j_{2}=1}^{n}1_{\left\{  j_{1}\neq
j_{2}\right\}  }\left(  X_{t_{j_{1}}}-X_{t_{j_{1}-1}}\right)  \left(
X_{t_{j_{2}}}-X_{t_{j_{2}-1}}\right)  g_{j_{1},j_{2}}^{\left(  2\right)
}\left(  \mathbf{0}\right) \nonumber\\
& +\frac{1}{3!}\sum_{j=1}^{n}\left(  X_{t_{j}}-X_{t_{j-1}}\right)
^{3}g_{j,j,j}^{\left(  3\right)  }\left(  \mathbf{0}\right) \nonumber\\
& +\frac{3}{3!}\sum_{j_{1}=1}^{n}\sum_{j_{2}=1}^{n}1_{\left\{  j_{1}\neq
j_{2}\right\}  }\left(  X_{t_{j_{1}}}-X_{t_{j_{1}-1}}\right)  ^{2}\left(
X_{t_{j_{2}}}-X_{t_{j_{2}-1}}\right)  g_{j_{1},j_{1},j_{2}}^{\left(  3\right)
}\left(  \mathbf{0}\right) \nonumber\\
& +\frac{1}{3!}\sum_{j_{1}=1}^{n}\sum_{j_{2}=1}^{n}\sum_{j_{3}=1}%
^{n}1_{\left\{  j_{1}\neq j_{2}\neq j_{3}\right\}  }\left(  X_{t_{j_{1}}%
}-X_{t_{j_{1}-1}}\right)  \left(  X_{t_{j_{2}}}-X_{t_{j_{2}-1}}\right)
\nonumber\\
& \times\left(  X_{t_{j_{3}}}-X_{t_{j_{3}-1}}\right)  g_{j_{1},j_{2},j_{3}%
}^{\left(  3\right)  }\left(  \mathbf{0}\right)  +\cdots,\label{multi}%
\end{align}
where $\left(  X_{t_{i}}-X_{t_{i-1}}\right)  ^{n}$'s are given by Theorem
\ref{newFormulaH} and we assume $X_{t_{0}}=0.$ \ The sums converge for every
$\omega\in\Omega$ because of (\ref{taylorCon}).

Since $0\leq t_{1}<t_{2}<\cdots<t_{n}$, the product of two iterated integrals
with non-overlapping limits results in an iterated integral: if $i\leq
j-1,\ u,v\in\left\{  1,2,3,...\right\}  $ and$~\phi_{i},\phi_{j}$ are the
predictable integrands,
\begin{align*}
\int_{t_{i-1}}^{t_{i}}\phi_{i}\ \mathrm{d}H_{s_{1}}^{\left(  u\right)  }%
\times\int_{t_{j-1}}^{t_{j}}\phi_{j}\ \mathrm{d}H_{r_{1}}^{\left(  v\right)
}  & =\int_{t_{j-1}}^{t_{j}}\int_{t_{i-1}}^{t_{i}}\phi_{i}\phi_{j}%
\ \mathrm{d}H_{s_{1}}^{\left(  u\right)  }\mathrm{d}H_{r_{1}}^{\left(
v\right)  }\\
& =\int_{0}^{t_{j}}\int_{0}^{t_{i}}1_{\left\{  s_{1}>t_{i-1}\right\}
}1_{\left\{  r_{1}>t_{j-1}\right\}  }\phi_{i}\phi_{j}\ \mathrm{d}H_{s_{1}%
}^{\left(  u\right)  }\mathrm{d}H_{r_{1}}^{\left(  v\right)  }\\
& =\int_{0}^{t_{j}}\int_{0}^{r_{1}}1_{\left\{  t_{i}>s_{1}>t_{i-1}\right\}
}1_{\left\{  r_{1}>t_{j-1}\right\}  }\phi_{i}\phi_{j}\ \mathrm{d}H_{s_{1}%
}^{\left(  u\right)  }\mathrm{d}H_{r_{1}}^{\left(  v\right)  },
\end{align*}
since $r_{1}>t_{j-1}\geq t_{i},$ giving an iterated integral. \ Hence, we get
a chaos expansion of $F$ in terms of iterated integrals with respect to
orthogonalized compensated power jump processes.

Note that in some applications, it is only necessary to apply Taylor's Theorem
directly to $F$ to obtain a PRP representation. \ While the approach given in
this section gives the CRP of $F$, each $\left(  X_{t_{j_{n}}}-X_{t_{j_{n}-1}%
}\right)  ^{m}$ consists of an infinite sum and therefore (\ref{multi}) is
composed of two levels of infinite sums. \ \cite{y07} applied Taylor's Theorem
directly to obtain the PRP of European and exotic option prices for hedging.

\section{Simulations using the explicit formula\label{SectionSimulation}}

To verify the theoretical results given in Section \ref{SectionPJ}, we
simulate the underlying L\'{e}vy processes and compare the values of $\left(
X_{t+t_{0}}-X_{t_{0}}\right)  ^{n}$ with the value given by its chaos
expansion. In simulations we apply the stochastic Euler scheme for the
stochastic differential equations (SDEs) of general L\'{e}vy processes. \ The
rate of convergence of this scheme for L\'{e}vy processes was discussed by
\cite{pt97}. \ For an up-to-date introduction to numerical solutions of SDEs,
see for example \cite{hk02}, \cite{h01}, \cite{k02} and \cite{kp99}. \ 

The processes considered are Gamma process and a combination of Wiener and
Gamma processes. \ For illustration, we run simulations for $k=4$ and $k=9$ in
the pure jump case and $k=5$ and $k=8$ for the combined case. \ The plots
produced are shown in Figures 1, 3, 5, 7 in Appendix \ref{AppendixPlots}
respectively. \ In the second and fourth simulations, we set $t_{0}=0.0099$
and $t_{0}=0.0019$ respectively. \ These simulations substantiate our explicit
formula for the CRP for $t_{0}\geq0$. \ We see that processes generated using
the CRP and those generated directly from the Gamma process jump at the same
time points. \ The differences between the two are plotted in Figures 2, 4, 6,
8 accordingly. \ Note that the axis of Figures 2, 4, 6, 8 are in much smaller
scales than those in Figures 1, 3, 5, 7. \ We deduce that the difference is
due to approximation errors of the stochastic Euler scheme. \ The errors
decrease with the step size $\Delta.\ \ $In each of the Figures 1, 3, 5, 7,
independent realizations of the Gamma and Wiener processes are used.

\section{Conclusion\label{SectionConclusion}}

L\'{e}vy processes were introduced in mathematical finance to improve the
performance of some of the financial models which are based on using Brownian
motion as the underlying process and to model stylized features observed in
financial processes.
The derivation of an explicit formula for the CRP has been of the focus of
considerable study, for previous work, see \cite{lsuv02}, \cite{bdlop03},
\cite{l04} and \cite{esv05}. \ In this paper, we have derived a computational
explicit formula for the construction of CRP of the powers of increments of
L\'{e}vy processes in terms of orthogonal compensated power jump processes and
its CRP in terms of Poisson random measures. \ \cite{j05} extended the CRP in
terms of power jump processes to a large class of semimartingales and we have
shown that in the L\'{e}vy case, our formula complements the one given by
\cite{j05}. \ Our explicit formula shows that the integrands of the stochastic
integrals in the CRP of the powers of increments of L\'{e}vy processes do not
depend on the integrating variables nor the starting time, which makes the
construction and simulation of the CRP much easier. \ The coefficients of the
CRP depend on $m_{i}$'s which represent the moments of the process with
respect to its L\'{e}vy measure. \ In this paper, we consider only L\'{e}vy
processes and their compensators are always of the form $m_{i}t.$ \ Using the
same calculation, it is trivial to extend the representation to
semimartingales which stochastic compensators have known representations.
\ The CRP of the pricing functions for some common financial derivatives can
be found by expressing the pricing functions in terms of powers of increments
of the underlying L\'{e}vy process using Taylor's expansion.

\bigskip

\setcounter{section}{0} \renewcommand{\thesection}{\Alph{section}}
\renewcommand{\thesubsection}{\Alph{section}.\arabic{subsection}} \numberwithin{equation}{section}%

\small

\begin{center}
{\large {\textbf{APPENDICES}}}
\end{center}

\setcounter{equation}{0}\renewcommand{\theequation}{A.\arabic{equation}}

\section{A note on the Nualart and Schoutens
representation\label{SectionNSTypo}}

\cite{ns00} derived the basic result for representing $\left(  X_{t+t_{0}%
}-X_{t_{0}}\right)  ^{k}$ when $t_{0}=0$ and $k=2.$ \ In the proof of the CRP,
\cite{ns00} made use of Proposition 2 in their paper, given in Section
\ref{SectionLevy} and the following equation derived from the Ito formula
(equation (5) in \cite{ns00}):%
\begin{align}
& \left(  X_{t+t_{0}}-X_{t_{0}}\right)  ^{k}\nonumber\\
& \ =\frac{\sigma^{2}}{2}k\left(  k-1\right)  \left(  \left(  X_{t+t_{0}%
}-X_{t_{0}}\right)  ^{k-2}t-\int_{0}^{t}s\,\ \mathrm{d}\left(  X_{s+t_{0}%
}-X_{t_{0}}\right)  ^{k-2}\right) \label{NS5_1}\\
& \ \ \ \ \ +\sum_{j=1}^{k}\binom{k}{j}\int_{t_{0}}^{t+t_{0}}\left(
X_{s-}-X_{t_{0}}\right)  ^{k-j}\mathrm{d}Y_{s}^{\left(  j\right)  }+\sum
_{j=1}^{k-1}\binom{k}{j}m_{j}t\left(  X_{t+t_{0}}-X_{t_{0}}\right)
^{k-j}\label{NS5}\\
& \ \ \ \ \ -\sum_{j=1}^{k-1}\binom{k}{j}m_{j}\int_{t_{0}}^{t+t_{0}%
}s\,\ \mathrm{d}\left(  X_{s}-X_{t_{0}}\right)  ^{k-j}+m_{k}t.\label{NS5_3}%
\end{align}
There is a small inaccuracy in this equation and we provide the corrected one
necessary for the derivation of the explicit formula. \ The second term in
(\ref{NS5}) should be
\[
\sum_{j=1}^{k-1}\binom{k}{j}m_{j}\left(  t+t_{0}\right)  \left(  X_{t+t_{0}%
}-X_{t_{0}}\right)  ^{k-j}%
\]
rather than $\sum_{j=1}^{k-1}\binom{k}{j}m_{j}t\left(  X_{t+t_{0}}-X_{t_{0}%
}\right)  ^{k-j}$. \ The error propagates from equation (4) in \cite{ns00}.
\ By integration by parts, $\sum_{j=1}^{k}\binom{k}{j}m_{j}\int_{t_{0}%
}^{t+t_{0}}\left(  X_{s-}-X_{t_{0}}\right)  ^{k-j}\mathrm{d}s$ should give%
\[
\sum_{j=1}^{k-1}\binom{k}{j}m_{j}\left(  t+t_{0}\right)  \left(  X_{t+t_{0}%
}-X_{t_{0}}\right)  ^{k-j}-\sum_{j=1}^{k-1}\binom{k}{j}m_{j}\int_{t_{0}%
}^{t+t_{0}}s\ \mathrm{d}\left(  X_{s}-X_{t_{0}}\right)  ^{k-j}+m_{k}t
\]
rather than the term
\[
\sum_{j=1}^{k-1}\binom{k}{j}m_{j}t\left(  X_{t+t_{0}}-X_{t_{0}}\right)
^{k-j}-\sum_{j=1}^{k-1}\binom{k}{j}m_{j}\int_{t_{0}}^{t+t_{0}}s\ \mathrm{d}%
\left(  X_{s}-X_{t_{0}}\right)  ^{k-j}+m_{k}t
\]
stated in \cite[p.114]{ns00}. \ Omitting $t_{0}$ makes the constant term of
the representation not equal to the expectation of $\left(  X_{t+t_{0}%
}-X_{t_{0}}\right)  ^{k}$ since it depends on $t_{0}$. \ Equation (5) in
\cite{ns00} should in fact be:%
\begin{align}
& \left(  X_{t+t_{0}}-X_{t_{0}}\right)  ^{k}\nonumber\\
& \ =\frac{\sigma^{2}}{2}k\left(  k-1\right)  \left(  \left(  X_{t+t_{0}%
}-X_{t_{0}}\right)  ^{k-2}t-\int_{0}^{t}s\ \mathrm{d}\left(  X_{s+t_{0}%
}-X_{t_{0}}\right)  ^{k-2}\right) \label{Corrected1}\\
& \ \ \ \ \ +\sum_{j=1}^{k}\binom{k}{j}\int_{t_{0}}^{t+t_{0}}\left(
X_{s-}-X_{t_{0}}\right)  ^{k-j}\mathrm{d}Y_{s}^{\left(  j\right)  }+\sum
_{j=1}^{k-1}\binom{k}{j}m_{j}\left(  t+t_{0}\right)  \left(  X_{t+t_{0}%
}-X_{t_{0}}\right)  ^{k-j}\label{Corrected2}\\
& \ \ \ \ -\sum_{j=1}^{k-1}\binom{k}{j}m_{j}\int_{t_{0}}^{t+t_{0}%
}s\ \mathrm{d}\left(  X_{s}-X_{t_{0}}\right)  ^{k-j}+m_{k}t.\label{Corrected}%
\end{align}
As an illustration of this representation, we derive $\left(  G_{t+t_{0}%
}-G_{t_{0}}\right)  ^{2}$ using (\ref{NS5_1})-(\ref{NS5_3}) to inspect the
constant terms. \ Note that $G_{t}$ is a L\'{e}vy process with $\sigma^{2}=0,$
so the terms in (\ref{NS5_1}) are equal to zero. \ We have%
\begin{align*}
\left(  G_{t+t_{0}}-G_{t_{0}}\right)  ^{2}  & =2\int_{t_{0}}^{t+t_{0}}\left(
G_{t_{1}-}-G_{t_{0}}\right)  \mathrm{d}\hat{G}_{t_{1}}^{\left(  1\right)
}+\int_{t_{0}}^{t+t_{0}}\mathrm{d}\hat{G}_{t_{1}}^{\left(  2\right)  }%
+2m_{1}t\left(  G_{t+t_{0}}-G_{t_{0}}\right) \\
& -2m_{1}\int_{t_{0}}^{t+t_{0}}t_{1}\ \mathrm{d}\left(  G_{t_{1}}-G_{t_{0}%
}\right)  +m_{2}t\\
& =2\int_{t_{0}}^{t+t_{0}}\left[  \left(  \hat{G}_{t_{1}-}^{\left(  1\right)
}-\hat{G}_{t_{0}}^{\left(  1\right)  }\right)  +m_{1}\left(  t_{1}%
-t_{0}\right)  \right]  \mathrm{d}\hat{G}_{t_{1}}^{\left(  1\right)  }%
+\int_{t_{0}}^{t+t_{0}}\mathrm{d}\hat{G}_{t_{1}}^{\left(  2\right)  }\\
& +2m_{1}t\left[  \left(  \hat{G}_{t+t_{0}}^{\left(  1\right)  }-\hat
{G}_{t_{0}}^{\left(  1\right)  }\right)  +m_{1}t\right]  -2m_{1}\int_{t_{0}%
}^{t+t_{0}}t_{1}\ \mathrm{d}\left[  \hat{G}_{t_{1}}^{\left(  1\right)  }%
+m_{1}t_{1}\right]  +m_{2}t\\
& =2\int_{t_{0}}^{t+t_{0}}\int_{t_{0}}^{t_{1}-}\mathrm{d}\hat{G}_{t_{2}%
}^{\left(  1\right)  }\mathrm{d}\hat{G}_{t_{1}}^{\left(  1\right)  }%
+\int_{t_{0}}^{t+t_{0}}\mathrm{d}\hat{G}_{t_{1}}^{\left(  2\right)  }%
+2m_{1}\left(  t-t_{0}\right)  \int_{t_{0}}^{t+t_{0}}\mathrm{d}\hat{G}_{t_{1}%
}^{\left(  1\right)  }\\
& +m_{1}^{2}t^{2}+m_{2}t-2m_{1}^{2}tt_{0}.
\end{align*}
The expectation of $2\int_{t_{0}}^{t+t_{0}}\int_{t_{0}}^{t_{1}-}$%
\textrm{d}$\hat{G}_{t_{2}}^{\left(  1\right)  }$\textrm{d}$\hat{G}_{t_{1}%
}^{\left(  1\right)  }+\int_{t_{0}}^{t+t_{0}}$\textrm{d}$\hat{G}_{t_{1}%
}^{\left(  2\right)  }+2m_{1}\left(  t-t_{0}\right)  \int_{t_{0}}^{t+t_{0}}%
$\textrm{d}$\hat{G}_{t_{1}}^{\left(  1\right)  }$ is zero since the
compensated processes $\hat{G}_{t}^{\left(  1\right)  }$ and $\hat{G}%
_{t}^{\left(  2\right)  }$ have zero means. \ We see that $m_{1}^{2}%
t^{2}+m_{2}t-2m_{1}^{2}tt_{0}$ depends on $t_{0}$ which in fact cannot be the
expectation of $\left(  G_{t+t_{0}}-G_{t_{0}}\right)  ^{2}$ since the
increments of $G_{t}$ are stationary. \ Starting from (\ref{Corrected1}%
)-(\ref{Corrected}), we can find that
\begin{equation}
\left(  G_{t+t_{0}}-G_{t_{0}}\right)  ^{2}=2\int_{t_{0}}^{t+t_{0}}\int_{t_{0}%
}^{t_{1}-}\mathrm{d}\hat{G}_{t_{2}}^{\left(  1\right)  }\mathrm{d}\hat
{G}_{t_{1}}^{\left(  1\right)  }+2m_{1}t\int_{t_{0}}^{t+t_{0}}\mathrm{d}%
\hat{G}_{t_{1}}^{\left(  1\right)  }+\int_{t_{0}}^{t+t_{0}}\mathrm{d}\hat
{G}_{t_{1}}^{\left(  2\right)  }+m_{1}^{2}t^{2}+m_{2}t.\label{g2pj}%
\end{equation}

\setcounter{equation}{0}\renewcommand{\theequation}{B.\arabic{equation}}

\section{Proof of Proposition \ref{TheoremC2}\label{AppendixTheoremC2}}

We prove this result using strong induction. \ Clearly, the proposition is
true for $k=1$ and $2.$ \ Assume the proposition is true for $k=n$, where $n$
is an integer $\geq1.$ \ Then for $k=n+1,$ firstly we prove that the sum of
the indices of all the $m_{q}$'s appear in each of the terms of $C_{t}%
^{\left(  n+1\right)  }$ (given by Proposition \ref{Theorem1}) are equal to
$n+1.$\ \ We have%
\begin{equation}
C_{t}^{\left(  n+1\right)  }=\sum_{j=1}^{n}\binom{n+1}{j}m_{j}tC_{t}^{\left(
n+1-j\right)  }-\sum_{j=1}^{n}\binom{n+1}{j}m_{j}\int_{0}^{t}t_{1}%
\ \mathrm{d}C_{t_{1}}^{\left(  n+1-j\right)  }+m_{n+1}t.\label{3.1.15}%
\end{equation}
By the induction step, the tuples of the indices of all the $m_{q}$'s
appearing in each of the terms of $C_{t}^{\left(  n+1-j\right)  }$ are
elements of $\mathcal{L}_{n+1-j}$ defined in (\ref{Lk}). \ Since we have
$m_{j}C_{t}^{\left(  n+1-j\right)  }$ appearing in the first term of
(\ref{3.1.15}), $m_{j}C_{t_{1}}^{\left(  n+1-j\right)  }$ in the second term
and $m_{n+1}$ in the last term, it is clear that the tuples of the indices of
all the $m_{q}$'s appearing in each of the terms of $C_{t}^{\left(
n+1\right)  }$ are elements of $\mathcal{L}_{n+1}.$ \ Now from (\ref{3.1.15}),
the first term can be proven to be%
\begin{align*}
& \sum_{j=1}^{n}\binom{n+1}{j}m_{j}tC_{t}^{\left(  n+1-j\right)  }\\
& ~=\sum_{j=1}^{n}\sum_{\phi_{n+1-j}=\left(  i_{1}^{(n+1-j)},i_{2}%
^{(n+1-j)},...,i_{l}^{(n+1-j)}\right)  \in\mathcal{L}_{n+1-j}}\frac{1}%
{l!}\left(  i_{1}^{(n+1-j)},i_{2}^{(n+1-j)},...,i_{l}^{(n+1-j)},j\right)  !\\
& ~~~~~~\times\left(  p_{1}^{\phi_{n+1-j}},p_{2}^{\phi_{n+1-j}},...,p_{n+1-j}%
^{\phi_{n+1-j}}\right)  !\left[  \prod\limits_{q\in\phi_{n+1-j}\cup\left\{
j\right\}  }m_{q}\right]  t^{l+1}%
\end{align*}
and the second term can be shown to be%
\begin{align*}
& -\sum_{j=1}^{n}\binom{n+1}{j}m_{j}\int_{0}^{t}t_{1}\ \mathrm{d}C_{t_{1}%
}^{\left(  n+1-j\right)  }\\
& ~=-\sum_{j=1}^{n}\sum_{\phi_{n+1-j}=\left(  i_{1}^{(n+1-j)},i_{2}%
^{(n+1-j)},...,i_{l}^{(n+1-j)}\right)  \in\mathcal{L}_{n+1-j}}\frac{1}%
{l!}\left(  i_{1}^{(n+1-j)},i_{2}^{(n+1-j)},...,i_{l}^{(n+1-j)},j\right)  !\\
& ~~~~~~\times\left(  p_{1}^{\phi_{n+1-j}},p_{2}^{\phi_{n+1-j}},...,p_{n+1-j}%
^{\phi_{n+1-j}}\right)  !\left[  \prod\limits_{q\in\phi_{n+1-j}\cup\left\{
j\right\}  }m_{q}\right]  \frac{l}{l+1}t^{l+1}.
\end{align*}
Hence,%
\begin{align*}
C_{t}^{\left(  n+1\right)  }  & =\sum_{j=1}^{n}\sum_{\phi_{n+1-j}=\left(
i_{1}^{(n+1-j)},i_{2}^{(n+1-j)},...,i_{l}^{(n+1-j)}\right)  \in\mathcal{L}%
_{n+1-j}}\frac{1}{l!}\left(  i_{1}^{(n+1-j)},i_{2}^{(n+1-j)},...,i_{l}%
^{(n+1-j)},j\right)  !\\
& \times\left(  p_{1}^{\phi_{n+1-j}},p_{2}^{\phi_{n+1-j}},...,p_{n+1-j}%
^{\phi_{n+1-j}}\right)  !\left[  \prod\limits_{q\in\phi_{n+1-j}\cup\left\{
j\right\}  }m_{q}\right]  t^{l+1}\frac{1}{l+1}+m_{n+1}t.
\end{align*}
Next we are going to prove that
\begin{align}
& \sum_{\phi_{n+1}=\left(  i_{1}^{(n+1)},i_{2}^{(n+1)},...,i_{l+1}%
^{(n+1)}\right)  \in\mathcal{L}_{n+1}}\frac{1}{\left(  l+1\right)  !}\left(
i_{1}^{(n+1)},i_{2}^{(n+1)},...,i_{l+1}^{(n+1)}\right)  !\nonumber\\
& \times\left(  p_{1}^{\phi_{n+1}},p_{2}^{\phi_{n+1}},...,p_{n+1}^{\phi_{n+1}%
}\right)  !\left[  \prod\limits_{q\in\phi_{n+1}}m_{q}\right]  t^{l+1}%
\nonumber\\
& ~=\sum_{j=1}^{n}\sum_{\phi_{n+1-j}=\left(  i_{1}^{(n+1-j)},i_{2}%
^{(n+1-j)},...,i_{l}^{(n+1-j)}\right)  \in\mathcal{L}_{n+1-j}}\frac{1}{\left(
l+1\right)  !}\left(  i_{1}^{(n+1-j)},i_{2}^{(n+1-j)},...,i_{l}^{(n+1-j)}%
,j\right)  !\nonumber\\
& ~~~~~~\times\left(  p_{1}^{\phi^{n+1-j}},p_{2}^{\phi^{n+1-j}},...,p_{n+1-j}%
^{\phi^{n+1-j}}\right)  !\left[  \prod\limits_{q\in\phi_{n+1-j}\cup\left\{
j\right\}  }m_{q}\right]  t^{l+1}+m_{n+1}t.\label{goal}%
\end{align}
On the R.H.S., we are adding a $j$ to each tuple $\left(  i_{1}^{(n+1-j)}%
,i_{2}^{(n+1-j)},...,i_{l}^{(n+1-j)}\right)  $ such that $\sum_{q=1}^{l}%
i_{q}^{(n+1-j)}+j=n+1$. \ Suppose $\phi_{n+1}=\left(  i_{1}^{(n+1)}%
,i_{2}^{(n+1)},...,i_{l+1}^{(n+1)}\right)  $ has one extra element compared to
the tuple $\left(  i_{1}^{(n+1-j)},i_{2}^{(n+1-j)},...,i_{l}^{(n+1-j)}\right)
$ and otherwise they are the same. \ Since $\sum_{q=1}^{l+1}i_{q}%
^{(n+1)}=n+1,$ to obtain $\left(  i_{1}^{(n+1)},i_{2}^{(n+1)},...,i_{l+1}%
^{(n+1)}\right)  $ from $\left(  i_{1}^{(n+1-j)},i_{2}^{(n+1-j)}%
,...,i_{l}^{(n+1-j)}\right)  $, we are adding an element $j$ to the latter
such that the sum of the tuple is equal to $n+1$. \ Suppose there are $r$
distinct value(s) in $\left(  i_{1}^{(n+1)},i_{2}^{(n+1)},...,i_{l+1}%
^{(n+1)}\right)  .$ \ Let $x_{1},x_{2},...,x_{r}$ be the distinct values in
$\left(  i_{1}^{(n+1)},i_{2}^{(n+1)},...,i_{l+1}^{(n+1)}\right)  $ and let
$f_{i},i=1,...,r$ be the number of times $x_{i}$ appears in $\left(
i_{1}^{(n+1)},i_{2}^{(n+1)},...,i_{l+1}^{(n+1)}\right)  .$ \ Note that
$\sum_{q=1}^{r}f_{q}$ is equal to the length of the tuple $\left(
i_{1}^{(n+1)},i_{2}^{(n+1)},...,i_{l+1}^{(n+1)}\right)  ,$ that is,
$\sum_{q=1}^{r}f_{q}=l+1.$ \ Since $\left(  i_{1}^{(n+1)},i_{2}^{(n+1)}%
,...,i_{l+1}^{(n+1)}\right)  $ can be obtained by adding an element $j$ to a
tuple $\left(  i_{1}^{(n+1-j)},i_{2}^{(n+1-j)},...,i_{l}^{(n+1-j)}\right)  $
whose elements add up to $n+1-j,$ $j$ can take one of the $r$ distinct
value(s): $x_{1},x_{2},...,x_{r}.$ \ For example, suppose $j=x_{i},$ then the
corresponding term on the R.H.S. of (\ref{goal}) is
\[
\frac{\left(  n+1\right)  !}{\left(  x_{1}!\right)  ^{f_{1}}\left(
x_{2}!\right)  ^{f_{2}}\cdots\left(  x_{i}!\right)  ^{f_{i}-1}\cdots\left(
x_{r}!\right)  ^{f_{r}}x_{i}!f_{1}!\cdots\left(  f_{i}-1\right)  !\cdots
f_{r}!}\left[  \prod\limits_{q\in\phi_{n+1}}m_{q}\right]  t^{l+1}\frac{1}%
{\sum_{q=1}^{r}f_{q}}.
\]
Summing up all the possible $j\in\left\{  x_{1},x_{2},...,x_{r}\right\}  ,$%
\begin{align*}
& \sum_{i=1}^{r}\frac{\left(  n+1\right)  !}{\left(  x_{1}!\right)  ^{f_{1}%
}\left(  x_{2}!\right)  ^{f_{2}}\cdots\left(  x_{i}!\right)  ^{f_{i}-1}%
\cdots\left(  x_{r}!\right)  ^{f_{r}}x_{i}!f_{1}!\cdots\left(  f_{i}-1\right)
!\cdots f_{r}!}\left[  \prod\limits_{q\in\phi_{n+1}}m_{q}\right]  t^{l+1}%
\frac{1}{\sum_{q=1}^{r}f_{q}}\\
& ~=\frac{\left(  n+1\right)  !}{\left(  x_{1}!\right)  ^{f_{1}}\cdots\left(
x_{r}!\right)  ^{f_{r}}f_{1}!\cdots f_{r}!}\frac{1}{\sum_{q=1}^{r}f_{q}%
}\left[  \prod\limits_{q\in\phi_{n+1}}m_{q}\right]  t^{l+1}\sum_{i=1}^{r}%
f_{i}\\
& ~=\frac{\left(  n+1\right)  !}{\left(  x_{1}!\right)  ^{f_{1}}\cdots\left(
x_{r}!\right)  ^{f_{r}}f_{1}!\cdots f_{r}!}\left[  \prod\limits_{q\in
\phi_{n+1}}m_{q}\right]  t^{l+1}.
\end{align*}
For the case $\phi_{n+1}=\left(  i_{1}^{(n+1)}\right)  ,$ it is clear that the
L.H.S. of (\ref{goal}) is equal to $m_{n+1}t.\ \ $Hence, by applying the same
argument to each possible tuple $\left(  i_{1}^{(n+1)},i_{2}^{(n+1)}%
,...,i_{l+1}^{(n+1)}\right)  \in\mathcal{L}_{n+1}$, we have proven
(\ref{goal}) and therefore%
\begin{align*}
C_{t}^{\left(  n+1\right)  }  & =\sum_{\phi_{n+1}=\left(  i_{1}^{(n+1)}%
,i_{2}^{(n+1)},...,i_{l}^{(n+1)}\right)  \in\mathcal{L}_{n+1}}\frac{1}%
{l!}\left(  i_{1}^{(n+1)},i_{2}^{(n+1)},...,i_{l}^{(n+1)}\right)  !\\
& \times\left(  p_{1}^{\phi_{n+1}},p_{2}^{\phi_{n+1}},...,p_{n+1}^{\phi_{n+1}%
}\right)  !\left[  \prod\limits_{q\in\phi_{n+1}}m_{q}\right]  t^{l}.
\end{align*}

\setcounter{equation}{0}\renewcommand{\theequation}{C.\arabic{equation}}

\section{Proof of Theorem \ref{Theorem3}\label{AppendixB}}

We prove the result using strong induction. \ Firstly we need to consider
$\mathcal{I}_{k}$ defined by equation (\ref{setI}). \ We need to know what
tuples are in $\mathcal{I}_{k+1}$ but not in $\mathcal{I}_{k}$, and these
correspond to those elements adding up exactly to $k+1$. \ Let $\mathcal{J}%
_{k+1}$ be the collection of these tuples, that is, $\mathcal{I}_{k+1}%
\equiv\mathcal{I}_{k}\cup\mathcal{J}_{k+1}.$ \ We have
\[
\mathcal{J}_{k+1}=\left\{  \left(  i_{1},i_{2},...,i_{j}\right)
|\ j\in\left\{  1,2,...,k+1\right\}  ,\ i_{p}\in\left\{  1,2,...,k+1\right\}
\text{ and }\sum_{p=1}^{j}i_{p}=k+1\right\}  .
\]
To construct $\mathcal{J}_{k+1}$ from $\mathcal{I}_{k}$, we can simply add an
element to the end of each tuple in $\mathcal{I}_{k}$ so that the elements of
each new tuple add up exactly to $k+1$, and finally including the tuple
$\left(  k+1\right)  $ in $\mathcal{J}_{k+1}$.\ 

We are going to prove by strong induction that $\left(  G_{t+t_{0}}-G_{t_{0}%
}\right)  ^{k}=\sum_{\theta_{k}\in\mathcal{I}_{k}}\Pi_{\theta_{k},t}^{\left(
k\right)  }\mathcal{S}_{\theta_{k},t,t_{0}}+C_{t}^{\left(  k\right)  }$ for
any non-negative integer $k$. \ For $k=0$, clearly both sides equal 1. \ For
$k=1$ and $2$, it can be checked easily that the proposition is true. \ Assume
the proposition is true for $k=0,1,2,...,n$, where $n$ is a positive integer.
\ Note that it is sufficient to prove the representation for $G_{t}^{n+1}$
only since we can always let $G_{t+t_{0}}-G_{t_{0}}=F_{t},$ which is also a
L\'{e}vy process and we have, the $i$-th power jump process of $\left\{
F_{t},t\geq0\right\}  ,$ $F_{t}^{\left(  i\right)  }=G_{t+t_{0}}^{\left(
i\right)  }-G_{t_{0}}^{\left(  i\right)  }$ for $i=1,2,3,...$. \ Since both
$\left\{  F_{t},t\geq0\right\}  $ and $\left\{  G_{t},t\geq0\right\}  $ are
created by the same infinitely divisible distribution, the compensators for
their $i$-th power jump processes are both equal to $m_{i}t.$ \ Hence, we have
the $i$-th compensated power jump process of $\left\{  F_{t},t\geq0\right\}
,$
\begin{equation}
\hat{F}_{t}^{\left(  i\right)  }=\hat{G}_{t+t_{0}}^{\left(  i\right)  }%
-\hat{G}_{t_{0}}^{\left(  i\right)  }.\label{FG}%
\end{equation}
For $k=n+1$, by (\ref{Corrected1})-(\ref{Corrected}),
\begin{align}
G_{t}^{n+1}  & =\sum_{j=1}^{n+1}\binom{n+1}{j}\int_{0}^{t}G_{t_{1}-}%
^{n+1-j}\ \mathrm{d}\hat{G}_{t_{1}}^{\left(  j\right)  }+\sum_{j=1}^{n}%
\binom{n+1}{j}m_{j}tG_{t}^{n+1-j}\nonumber\\
& -\sum_{j=1}^{n}\binom{n+1}{j}m_{j}\int_{0}^{t}t_{1}\ \mathrm{d}G_{t_{1}%
}^{n+1-j}+m_{n+1}t.\label{CorrectedGn+1}%
\end{align}
Firstly, we want to prove that all the stochastic integrals in $G_{t}^{n+1}$
is of the form $\mathcal{S}_{\theta_{n+1},t,0},\ $where $\theta_{n+1}%
\in\mathcal{I}_{n+1}.$ \ From (\ref{CorrectedGn+1}), it is clear that the
first term is the only term introducing new stochastic integrals which are not
in $\mathcal{I}_{n}.$ \ The general term of the stochastic integrals in the
first term is%
\begin{equation}
\int_{0}^{t}G_{t_{1}-}^{n+1-j}\ \mathrm{d}\hat{G}_{t_{1}}^{\left(  j\right)
},\text{\ \ \ }j=1,2,...,n+1.\label{SI1}%
\end{equation}
By assumption,%
\[
G_{t_{1}-}^{n+1-j}=\sum_{\theta_{n+1-j}\in\mathcal{I}_{n+1-j}}\Pi
_{\theta_{n+1-j},t_{1}}^{\left(  n+1-j\right)  }\mathcal{S}_{\theta
_{n+1-j},t_{1},0}+C_{t_{1}}^{\left(  n+1-j\right)  },\ j=1,2,...,n+1.
\]
When $j=1$ in (\ref{SI1}), we have $\int_{0}^{t}G_{t_{1}-}^{n}$\textrm{d}%
$\hat{G}_{t_{1}}^{\left(  1\right)  },$ meaning that we are adding a 1 to the
end of all tuples in $\mathcal{I}_{n}$. \ Since by definition%
\[
\mathcal{I}_{n}=\left\{  \left(  i_{1},i_{2},...,i_{j}\right)  |\ j\in\left\{
1,2,...,n\right\}  ,\ i_{p}\in\left\{  1,2,...,n\right\}  \text{ and }%
\sum_{p=1}^{j}i_{p}\leq n\right\}  ,
\]
we know that the sum of the elements of the new tuples we get from adding a 1
to the end of each tuple of $\mathcal{I}_{n}$ is less than or equal to $n+1 $.
\ Similarly, when $j=2$, we have $\int_{0}^{t}G_{t_{1}-}^{n-1}$\textrm{d}%
$\hat{G}_{t_{1}}^{\left(  2\right)  },$ meaning that we are adding a 2 to the
end of all tuples in $\mathcal{I}_{n-1}$ and since by definition%
\[
\mathcal{I}_{n-1}=\left\{  \left(  i_{1},i_{2},...,i_{j}\right)
|\ j\in\left\{  1,2,...,n-1\right\}  ,\ i_{p}\in\left\{  1,2,...,n-1\right\}
\text{ and }\sum_{p=1}^{j}i_{p}\leq n-1\right\}  ,
\]
we know that the sum of the elements of the new tuples we get from adding a 2
to the end of each tuple of $\mathcal{I}_{n-1}$ is less than or equal to
$n+1$.\ \ We can continue the same argument until $j=n$. \ When $j=n+1$, we
have $\int_{0}^{t}$\textrm{d}$\hat{G}_{t_{1}}^{\left(  n+1\right)  }.$ \ Since
$\mathcal{I}_{n}\supset\mathcal{I}_{n-1}\supset...\supset\mathcal{I}%
_{2}\supset\mathcal{I}_{1},$ the above way of introducing new stochastic
integrals is the same as adding an element to the end of each tuple in
$\mathcal{I}_{n}$ so that the elements of each new tuple add up exactly to
$n+1$. \ Hence all the elements in $\mathcal{J}_{n+1}$ have been created and
since $\mathcal{I}_{n+1}\equiv\mathcal{I}_{n}\cup\mathcal{J}_{n+1},$ we have
proved that all the stochastic integrals in $G_{t}^{n+1}$ have the form
$\mathcal{S}_{\theta_{n+1},t,0},$ where $\theta_{n+1}\in\mathcal{I}_{n+1}. $

By definition, $C_{t}^{\left(  n+1\right)  }$ is the term in $G_{t}^{n+1}$ not
containing any stochastic integral. \ Hence it is correct to write
$C_{t}^{\left(  n+1\right)  }$ as the final term.

Finally, we want to consider the coefficients of the stochastic integrals,
that is, we are going to prove Proposition \ref{Theorem2}. \ By assumption of
the induction step, we have%
\begin{align}
G_{t}^{n+1}  & =\sum_{j=1}^{n}\binom{n+1}{j}\sum_{\theta_{n+1-j}\in
\mathcal{I}_{n+1-j}}\int_{0}^{t}\Pi_{\theta_{n+1-j},t_{1}}^{\left(
n+1-j\right)  }\mathcal{S}_{\theta_{n+1-j},t_{1},0}\ \mathrm{d}\hat{G}_{t_{1}%
}^{\left(  j\right)  }\nonumber\\
& +\sum_{j=1}^{n}\binom{n+1}{j}m_{j}t\sum_{\theta_{n+1-j}\in\mathcal{I}%
_{n+1-j}}\Pi_{\theta_{n+1-j},t}^{\left(  n+1-j\right)  }\mathcal{S}%
_{\theta_{n+1-j},t,0}\nonumber\\
& -\sum_{j=1}^{n}\binom{n+1}{j}m_{j}\sum_{\theta_{n+1-j}\in\mathcal{I}%
_{n+1-j}}\int_{0}^{t}t_{1}\ \mathrm{d}\left[  \Pi_{\theta_{n+1-j},t_{1}%
}^{\left(  n+1-j\right)  }\mathcal{S}_{\theta_{n+1-j},t_{1},0}\right]
\nonumber\\
& +\sum_{j=1}^{n}\binom{n+1}{j}\int_{0}^{t}C_{t_{1}}^{\left(  n+1-j\right)
}\ \mathrm{d}\hat{G}_{t_{1}}^{\left(  j\right)  }+\sum_{j=1}^{n}\binom{n+1}%
{j}m_{j}tC_{t}^{\left(  n+1-j\right)  }\nonumber\\
& -\sum_{j=1}^{n}\binom{n+1}{j}m_{j}\int_{0}^{t}t_{1}\ \mathrm{d}\left[
C_{t_{1}}^{\left(  n+1-j\right)  }\right]  +\int_{0}^{t}\mathrm{d}\hat
{G}_{t_{1}}^{\left(  n+1\right)  }+m_{n+1}t\nonumber\\
& =L_{1}+L_{2}+L_{3}+L_{4}+L_{5}+L_{6}+\int_{0}^{t}\mathrm{d}\hat{G}_{t_{1}%
}^{\left(  n+1\right)  }+m_{n+1}t.
\end{align}
Let $\mathcal{K}_{l,s}=\left\{  \left(  i_{1},...,i_{l}\right)  |i_{j}%
\in\left\{  1,2,...,s\right\}  \text{ and }\sum_{j=1}^{l}i_{j}=s\right\}  .$
\ Since the length of a tuple must not be greater than the sum of all the
elements in the tuple (because an element must be at least 1), $l\leq s.$ \ By
definition, we have $\mathcal{I}_{n}%
=\mathop{\displaystyle \bigcup }\limits_{s=1}^{n}%
\mathop{\displaystyle \bigcup }\limits_{l=1}^{s}\mathcal{K}_{l,s}.$ \ For any
$\theta_{l,s}\in\mathcal{K}_{l,s}$, let $\theta_{l,s}=\left(  i_{1}%
^{\theta_{l,s}},i_{2}^{\theta_{l,s}},...,i_{l}^{\theta_{l,s}}\right)  .$ \ It
is obvious from Proposition \ref{Theorem1} that $C_{t}^{\left(  k\right)  }$
has the form $C_{t}^{\left(  k\right)  }=q_{0}^{\left(  k\right)  }%
+q_{1}^{\left(  k\right)  }t+q_{2}^{\left(  k\right)  }t^{2}+...+q_{k}%
^{\left(  k\right)  }t^{k}.$ \ Note that $q_{0}^{\left(  k\right)  }$ is
non-zero only when $k=0$. \ When $k=0$, by definition $C_{t}^{\left(
k\right)  }=1$, so we have $q_{0}^{\left(  0\right)  }=1$. \ We need to find
out the recursive relationships between the $q_{r}^{\left(  k\right)  }$'s.
\ From (\ref{C}), for $k>1,$%
\begin{align*}
q_{1}^{\left(  k\right)  }t+q_{2}^{\left(  k\right)  }t^{2}+\cdots
+q_{k}^{\left(  k\right)  }t^{k}  & =\sum_{j=1}^{k-1}\binom{k}{j}m_{j}t\left[
q_{1}^{\left(  k-j\right)  }t+q_{2}^{\left(  k-j\right)  }t^{2}+\cdots
+q_{k-j}^{\left(  k-j\right)  }t^{k-j}\right] \\
& -\sum_{j=1}^{k-1}\binom{k}{j}m_{j}\int_{0}^{t}t_{1}\ \mathrm{d}\left[
q_{1}^{\left(  k-j\right)  }t_{1}+q_{2}^{\left(  k-j\right)  }t_{1}^{2}%
+\cdots+q_{k-j}^{\left(  k-j\right)  }t_{1}^{k-j}\right]  +m_{k}t\\
& =m_{k}t+\sum_{j=1}^{k-1}\binom{k}{j}m_{j}\left[  q_{1}^{\left(  k-j\right)
}t^{2}+q_{2}^{\left(  k-j\right)  }t^{3}+\cdots+q_{k-j}^{\left(  k-j\right)
}t^{k-j+1}\right] \\
& -\sum_{j=1}^{k-1}\binom{k}{j}m_{j}\left[  \frac{1}{2}q_{1}^{\left(
k-j\right)  }t^{2}+\frac{2}{3}q_{2}^{\left(  k-j\right)  }t^{3}+\cdots
+\frac{k-j}{k-j+1}q_{k-j}^{\left(  k-j\right)  }t^{k-j+1}\right]  .
\end{align*}
By comparing the coefficients of $t$, $q_{1}^{\left(  k\right)  }=m_{k}.$ \ By
comparing the coefficients of $t^{r}$, $r=2,...,k$,%
\begin{equation}
q_{r}^{\left(  k\right)  }=\sum_{j=1}^{k+1-r}\binom{k}{j}m_{j}q_{r-1}^{\left(
k-j\right)  }-\sum_{j=1}^{k+1-r}\binom{k}{j}m_{j}\frac{r-1}{r}q_{r-1}^{\left(
k-j\right)  }=\frac{1}{r}\sum_{j=1}^{k+1-r}\binom{k}{j}m_{j}q_{r-1}^{\left(
k-j\right)  }.\label{q}%
\end{equation}
To ease notation, we let
\begin{align*}
\mathbf{F}_{1}  & =\frac{\left(  n+1\right)  !}{i_{1}^{\theta_{l,s}}%
!i_{2}^{\theta_{l,s}}!\cdots i_{l}^{\theta_{l,s}}!j!},~\mathbf{F}_{2}%
=\frac{\mathbf{F}_{1}}{\left(  n+1-j-s\right)  !},~\mathbf{G}_{i}^{j}=\hat
{G}_{t_{i}}^{\left(  i_{j}^{\theta_{l,s}}\right)  },\ \mathbf{I}_{1}=\int
_{0}^{t_{1}-}\cdots\int_{0}^{t_{l}-}\mathrm{d}\mathbf{G}_{l+1}^{1}%
\cdots\mathrm{d}\mathbf{G}_{2}^{l},\\
\mathbf{I}_{2}  & =\int_{0}^{t_{1}-}\cdots\int_{0}^{t_{l-1}-}\mathrm{d}%
\mathbf{G}_{l}^{1}\cdots\mathrm{d}\mathbf{G}_{2}^{l-1}\ ,\ \mathbf{I}_{3}%
=\int_{0}^{t}\mathbf{I}_{2}\ \mathrm{d}\mathbf{G}_{1}^{l},\ \mathbf{\tilde{q}%
}_{i}=q_{i}^{\left(  n+1-j-s\right)  }.
\end{align*}
Note that it is only for simplicity in writing out the equations. \ When doing
calculation, we should always use the long but clear notation. \ So, we have%
\begin{align*}
L_{1}  & =\sum_{j=1}^{n}\sum_{s=1}^{n+1-j}\sum_{l=1}^{s}\sum_{\theta_{l,s}%
\in\mathcal{K}_{l,s}}\int_{0}^{t}\mathbf{F}_{2}\sum_{w=0}^{n+1-j-s}%
\mathbf{\tilde{q}}_{w}t_{1}^{w}\mathbf{I}_{1}\ \mathrm{d}\hat{G}_{t_{1}%
}^{\left(  j\right)  }.\\
L_{2}  & =\sum_{j=1}^{n}m_{j}\sum_{s=1}^{n+1-j}\left\{  1_{\left\{
s=n+1-j\right\}  }t\sum_{l=1}^{s}\sum_{\theta_{l,s}\in\mathcal{K}_{l,s}%
}\mathbf{F}_{1}\mathbf{I}_{3}+1_{\left\{  s\leq n-j\right\}  }\sum_{l=1}%
^{s}\sum_{\theta_{l,s}\in\mathcal{K}_{l,s}}\mathbf{F}_{2}\sum_{w=1}%
^{n+1-j-s}\mathbf{\tilde{q}}_{w}t^{w+1}\mathbf{I}_{3}\right\} \\
L_{3}  & =-\sum_{j=1}^{n}m_{j}\sum_{s=1}^{n+1-j}\sum_{l=1}^{s}\sum
_{\theta_{l,s}\in\mathcal{K}_{l,s}}\mathbf{F}_{2}\left\{  \mathbf{\tilde{q}%
}_{0}\int_{0}^{t}t_{1}\mathbf{I}_{2}\ \mathrm{d}\mathbf{G}_{1}^{l}+\frac{1}%
{2}\mathbf{\tilde{q}}_{1}t^{2}\mathbf{I}_{3}+\frac{1}{2}\mathbf{\tilde{q}}%
_{1}\int_{0}^{t}t_{1}^{2}\mathbf{I}_{2}\ \mathrm{d}\mathbf{G}_{1}^{l}\right.
\\
& \left.  +\sum_{w=2}^{n+1-j-s}\mathbf{\tilde{q}}_{w}\left\{  \frac{w}%
{w+1}t^{w+1}\mathbf{I}_{3}+\frac{1}{w+1}\int_{0}^{t}t_{1}^{w+1}\mathbf{I}%
_{2}\ \mathrm{d}\mathbf{G}_{1}^{l}\right\}  \right\}  .\\
L_{4}  & =\sum_{j=1}^{n}\binom{n+1}{j}\int_{0}^{t}\sum_{w=1}^{n+1-j}%
q_{w}^{\left(  n+1-j\right)  }t_{1}^{w}\ \mathrm{d}\hat{G}_{t_{1}}^{\left(
j\right)  }.\\
L_{5}  & =\sum_{j=1}^{n}\binom{n+1}{j}m_{j}\sum_{w=1}^{n+1-j}q_{w}^{\left(
n+1-j\right)  }t^{w+1}.\\
L_{6}  & =-\sum_{j=1}^{n}\binom{n+1}{j}m_{j}\sum_{w=1}^{n+1-j}q_{w}^{\left(
n+1-j\right)  }\left\{  \frac{w}{w+1}t^{w+1}\right\}  .
\end{align*}
Next, consider $L_{1}$ and $L_{3}.$ \ Let $u,v\in\left\{  1,2,...,n-1\right\}
$ and $u+v\leq n.$ \ In $L_{1}$, when $j=u,$ $s=v$ (hence $s\leq n-j$),
\begin{align}
L_{1}  & =\sum_{l=1}^{v}\sum_{\theta_{l,v}\in\mathcal{K}_{l,v}}\int_{0}%
^{t}\frac{\left(  n+1\right)  !}{i_{1}^{\theta_{l,v}}!i_{2}^{\theta_{l,v}%
}!\cdots i_{l}^{\theta_{l,v}}!u!}\frac{1}{\left(  n+1-u-v\right)  !}\left\{
m_{n+1-u-v}t_{1}+\sum_{w=2}^{n+1-u-v}q_{w}^{\left(  n+1-u-v\right)  }t_{1}%
^{w}\right\} \nonumber\\
& \times\int_{0}^{t_{1}-}\int_{0}^{t_{2}-}\cdots\int_{0}^{t_{l}-}%
\mathrm{d}\hat{G}_{t_{l+1}}^{\left(  i_{1}^{\theta_{l,v}}\right)  }%
\cdots\mathrm{d}\hat{G}_{t_{2}}^{\left(  i_{l}^{\theta_{l,v}}\right)
}\mathrm{d}\hat{G}_{1}^{\left(  u\right)  }.
\end{align}
Since $s=v,$ $l\in\left\{  1,2,...,v\right\}  ,$ we have by definition
$\left(  i_{1}^{\theta_{l,v}},i_{2}^{\theta_{l,v}},...,i_{l}^{\theta_{l,v}%
}\right)  \in\mathcal{J}_{v}.$ \ In $L_{3}$, when $j=n+1-u-v$ (hence
$j\in\left\{  1,...,n-1\right\}  $), $s=u+v$ (hence $s=n+1-j$) and
$i_{l}^{\theta_{l,s}}=u$ (hence $i_{l}^{\theta_{l,s}}<s$),%
\begin{align*}
L_{3}  & =-m_{n+1-u-v}\sum_{l=1}^{u+v}\sum_{\theta_{l,u+v}\in\mathcal{K}%
_{l,u+v}}\frac{\left(  n+1\right)  !}{i_{1}^{\theta_{l,u+v}}!i_{2}%
^{\theta_{l,u+v}}!\cdots i_{l-1}^{\theta_{l,u+v}}!u!\left(  n+1-u-v\right)
!}\\
& \times\int_{0}^{t}t_{1}\int_{0}^{t_{1}-}\int_{0}^{t_{2}-}\cdots\int
_{0}^{t_{l-1}-}\mathrm{d}\hat{G}_{t_{l}}^{\left(  i_{1}^{\theta_{l,u+v}%
}\right)  }\cdots\mathrm{d}\hat{G}_{t_{2}}^{\left(  i_{l-1}^{\theta_{l,u+v}%
}\right)  }\mathrm{d}\hat{G}_{t_{1}}^{\left(  u\right)  }.
\end{align*}
Since $s=u+v$ and $i_{l}^{\theta_{l,s}}=u$, $\sum_{p=1}^{l-1}i_{p}%
^{\theta_{l,u+v}}=v,$ we have by definition $\left(  i_{1}^{\theta_{l,u+v}%
},i_{2}^{\theta_{l,u+v}},...,i_{l-1}^{\theta_{l,u+v}}\right)  \in
\mathcal{J}_{v}.$ \ Hence the terms%
\begin{align*}
& \sum_{l=1}^{v}\sum_{\theta_{l,v}\in\mathcal{K}_{l,v}}\int_{0}^{t}%
\frac{\left(  n+1\right)  !}{i_{1}^{\theta_{l,v}}!i_{2}^{\theta_{l,v}}!\cdots
i_{l}^{\theta_{l,v}}!u!}\frac{1}{\left(  n+1-u-v\right)  !}m_{n+1-u-v}t_{1}\\
& \times\int_{0}^{t_{1}-}\int_{0}^{t_{2}-}\cdots\int_{0}^{t_{l}-}%
\mathrm{d}\hat{G}_{t_{l+1}}^{\left(  i_{1}^{\theta_{l,v}}\right)  }%
\cdots\mathrm{d}\hat{G}_{t_{2}}^{\left(  i_{l}^{\theta_{l,v}}\right)
}\mathrm{d}\hat{G}_{1}^{\left(  u\right)  }%
\end{align*}
in $L_{1}$ and%
\begin{align*}
& -m_{n+1-u-v}\sum_{l=1}^{u+v}\sum_{\theta_{l,u+v}\in\mathcal{K}_{l,u+v}}%
\frac{\left(  n+1\right)  !}{i_{1}^{\theta_{l,u+v}}!i_{2}^{\theta_{l,u+v}%
}!\cdots i_{l-1}^{\theta_{l,u+v}}!u!\left(  n+1-u-v\right)  !}\\
& \times\int_{0}^{t}t_{1}\int_{0}^{t_{1}-}\int_{0}^{t_{2}-}\cdots\int
_{0}^{t_{l-1}-}\mathrm{d}\hat{G}_{t_{l}}^{\left(  i_{1}^{\theta_{l,u+v}%
}\right)  }\cdots\mathrm{d}\hat{G}_{t_{2}}^{\left(  i_{l-1}^{\theta_{l,u+v}%
}\right)  }\mathrm{d}\hat{G}_{t_{1}}^{\left(  u\right)  }%
\end{align*}
cancel each other. \ So we now have%
\begin{align*}
L_{1}  & =\sum_{j=1}^{n}\left\{  1_{\left\{  j\leq n-1\right\}  }\sum
_{s=1}^{n+1-j}\left\{  1_{\left\{  s\leq n-j\right\}  }\sum_{l=1}^{s}%
\sum_{\theta_{l,s}\in\mathcal{K}_{l,s}}\int_{0}^{t}\mathbf{F}_{2}\sum
_{w=2}^{n+1-j-s}\mathbf{\tilde{q}}_{w}t_{1}^{w}\mathbf{I}_{1}\ \mathrm{d}%
\hat{G}_{t_{1}}^{\left(  j\right)  }\right.  \right. \\
& \left.  \left.  +1_{\left\{  s=n+1-j\right\}  }\sum_{l=1}^{s}\sum
_{\theta_{l,s}\in\mathcal{K}_{l,s}}\int_{0}^{t}\ \mathbf{F}_{1}\mathbf{I}%
_{1}\ \mathrm{d}\hat{G}_{t_{1}}^{\left(  j\right)  }\right\}  +1_{\left\{
j=n\right\}  }\left(  n+1\right)  \int_{0}^{t}\int_{0}^{t_{1}-}\mathrm{d}%
\hat{G}_{t_{2}}^{\left(  1\right)  }\mathrm{d}\hat{G}_{t_{1}}^{\left(
n\right)  }\right\}  .
\end{align*}
Since $q_{0}^{\left(  k\right)  }=0$ for $k>0$,
\begin{align*}
L_{3}  & =-\sum_{j=1}^{n}m_{j}\left\{  1_{\left\{  j\leq n-1\right\}  }%
\sum_{s=1}^{n+1-j}\left\{  1_{\left\{  s\leq n-j\right\}  }\sum_{l=1}^{s}%
\sum_{\theta_{l,s}\in\mathcal{K}_{l,s}}\mathbf{F}_{2}\left\{  \frac{1}%
{2}\mathbf{\tilde{q}}_{1}t^{2}\mathbf{I}_{3}+\frac{1}{2}\mathbf{\tilde{q}}%
_{1}\int_{0}^{t}t_{1}^{2}\mathbf{I}_{2}\ \mathrm{d}\mathbf{G}_{1}^{l}\right.
\right.  \right. \\
& \left.  +\sum_{w=2}^{n+1-j-s}\mathbf{\tilde{q}}_{w}\left\{  \frac{w}%
{w+1}t^{w+1}\mathbf{I}_{3}+\frac{1}{w+1}\int_{0}^{t}t_{1}^{w+1}\mathbf{I}%
_{2}\ \mathrm{d}\mathbf{G}_{1}^{l}\right\}  \right\} \\
& \left.  \left.  +1_{\left\{  s=n+1-j\right\}  }\frac{\left(  n+1\right)
!}{\left(  n+1-j\right)  !j!}\int_{0}^{t}t_{1}\ \mathrm{d}\hat{G}_{t_{1}%
}^{\left(  n+1-j\right)  }\right\}  +1_{\left\{  j=n\right\}  }m_{n}\left(
n+1\right)  \int_{0}^{t}t_{1}\ \mathrm{d}\hat{G}_{t_{1}}^{\left(  1\right)
}\right\}  .
\end{align*}
Next, consider $L_{3}$ and $L_{4}.$ \ Let $u\in\{1,2,...,n\}.$ \ In $L_{3}$,
when $j=n+1-u$ (hence $j\in\left\{  1,...,n\right\}  $)$,\ s=u$ (hence
$s=n+1-j$) and $i_{l}^{\theta_{l,s}}=u$ (hence $i_{l}^{\theta_{l,s}}=s$), we
have%
\[
L_{3}=m_{n+1-u}\frac{\left(  n+1\right)  !}{u!\left(  n+1-u\right)  !}\int
_{0}^{t}t_{1}\ \mathrm{d}\hat{G}_{t_{1}}^{\left(  u\right)  }.
\]
In $L_{4}$, when $j=u$, we have%
\begin{equation}
L_{4}=\binom{n+1}{u}\int_{0}^{t}\left\{  m_{n+1-u}t_{1}+\sum_{w=2}%
^{n+1-u}q_{w}^{\left(  n+1-u\right)  }t_{1}^{w}\right\}  \mathrm{d}\hat
{G}_{t_{1}}^{\left(  u\right)  }.\nonumber
\end{equation}
Hence the terms%
\[
m_{n+1-u}\frac{\left(  n+1\right)  !}{u!\left(  n+1-u\right)  !}\int_{0}%
^{t}t_{1}\ \mathrm{d}\hat{G}_{t_{1}}^{\left(  u\right)  }%
\]
in $L_{3}$ and
\[
\binom{n+1}{u}\int_{0}^{t}m_{n+1-u}t_{1}\ \mathrm{d}\hat{G}_{t_{1}}^{\left(
u\right)  }%
\]
cancel each other. \ In $L_{3}$, since the terms where $\left(  j=n\right)  $
and $\left(  j\leq n-1,\ s=n+1-j\right)  $ get cancelled, we can sum $j$ from
1 to $n-1$ and sum $s$ from 1 to $n-j$.%
\begin{align*}
L_{3}  & =-\sum_{j=1}^{n-1}m_{j}\sum_{s=1}^{n-j}\sum_{l=1}^{s}\sum
_{\theta_{l,s}\in\mathcal{K}_{l,s}}\mathbf{F}_{2}\sum_{w=1}^{n+1-j-s}%
\mathbf{\tilde{q}}_{w}\left\{  \frac{w}{w+1}t^{w+1}\mathbf{I}_{3}+\frac
{1}{w+1}\int_{0}^{t}t_{1}^{w+1}\mathbf{I}_{2}\ \mathrm{d}\mathbf{G}_{1}%
^{l}\right\}  .\\
L_{4}  & =\sum_{j=1}^{n-1}\left\{  \binom{n+1}{j}\int_{0}^{t}\sum
_{w=2}^{n+1-j}\frac{1}{w}\left[  \sum_{z=1}^{n+2-j-w}\binom{n+1-j}{z}%
m_{z}q_{w-1}^{\left(  n+1-j-z\right)  }\right]  t_{1}^{w}\ \mathrm{d}\hat
{G}_{t_{1}}^{\left(  j\right)  }\right\}  .
\end{align*}
by (\ref{q}).

Consider $L_{4}$ and $L_{3}$. \ Let $u\in\left\{  1,2,...,n-1\right\}  $,
$v\in\left\{  1,2,...,n-1\right\}  ,\ u+v\leq n,\ x\in\left\{
1,2,...,v\right\}  $ and hence$\ x+u\leq n.$ \ In $L_{4}$, when
$j=u,\ w=n+2-u-v$ (hence $w\in\left\{  2,...,n+1-j\right\}  $)$,\ z=x$ (hence
$z\in\left\{  1,...,n+2-j-w\right\}  $)$,$%
\begin{equation}
L_{4}=\binom{n+1}{u}\int_{0}^{t}\frac{1}{n+2-u-v}\binom{n+1-u}{x}%
m_{x}q_{n+1-u-v}^{\left(  n+1-u-x\right)  }t_{1}^{n+2-u-v}\ \mathrm{d}\hat
{G}_{t_{1}}^{\left(  u\right)  }.\label{L4}%
\end{equation}
In $L_{3}$, when $j=x$ (hence $j\in\left\{  1,...,n-1\right\}  $), $s=u$
(hence $s\leq n-j$)$,\ i_{l}^{\theta_{l,s}}=u$ (hence $i_{l}^{\theta_{l,s}}%
=s$)$,\ w=n+1-u-v$ (hence $w\leq n+1-s-j$)$,$
\begin{align*}
L_{3}  & =m_{x}\frac{\left(  n+1\right)  !}{u!x!}\frac{1}{\left(
n+1-x-u\right)  !}q_{n+1-u-v}^{\left(  n+1-u-x\right)  }\\
& \times\left\{  \frac{n+1-u-v}{n+2-u-v}t^{n+2-u-v}\int_{0}^{t}\mathrm{d}%
\hat{G}_{t_{1}}^{\left(  u\right)  }+\frac{1}{n+2-u-v}\int_{0}^{t}%
t_{1}^{n+2-u-v}\ \mathrm{d}\hat{G}_{t_{1}}^{\left(  u\right)  }\right\}
\end{align*}
where the $2^{nd}$ term cancels (\ref{L4}). \ So now we have%
\begin{align*}
L_{4}  & =0.\\
L_{3}  & =-\sum_{j=1}^{n-1}m_{j}\sum_{s=1}^{n-j}\left\{  1_{\left\{
s=1\right\}  }\frac{\left(  n+1\right)  !}{j!}\frac{1}{\left(  n-j\right)
!}\sum_{w=1}^{n-j}q_{w}^{\left(  n-j\right)  }\frac{w}{w+1}t^{w+1}\int_{0}%
^{t}\mathrm{d}\hat{G}_{t_{1}}^{\left(  1\right)  }\right. \\
& +1_{\left\{  2\leq s\leq n-j\right\}  }\sum_{l=1}^{s}\sum_{\theta_{l,s}%
\in\mathcal{K}_{l,s}}\left\{  1_{\left\{  i_{l}^{\theta_{l,s}}<s\right\}
}\mathbf{F}_{2}\sum_{w=1}^{n+1-j-s}\mathbf{\tilde{q}}_{w}\left\{  \frac
{w}{w+1}t^{w+1}\mathbf{I}_{3}+\frac{1}{w+1}\int_{0}^{t}t_{1}^{w+1}%
\mathbf{I}_{2}\ \mathrm{d}\mathbf{G}_{1}^{l}\right\}  \right. \\
& \left.  \left.  +1_{\left\{  i_{l}^{\theta_{l,s}}=s\right\}  }\frac{\left(
n+1\right)  !}{s!j!}\frac{1}{\left(  n+1-j-s\right)  !}\sum_{w=1}%
^{n+1-j-s}\mathbf{\tilde{q}}_{w}\frac{w}{w+1}t^{w+1}\int_{0}^{t}\mathrm{d}%
\hat{G}_{t_{1}}^{\left(  s\right)  }\right\}  \right\}  .
\end{align*}
Next, consider $L_{1}$ and $L_{3}$. \ By the equation for $q_{w}^{\left(
n+1-j-s\right)  }$ given in (\ref{q}), we have%
\begin{align*}
L_{1}  & =\sum_{j=1}^{n}\left\{  1_{\left\{  j\leq n-1\right\}  }\sum
_{s=1}^{n+1-j}\left\{  1_{\left\{  s\leq n-1-j\right\}  }\sum_{l=1}^{s}%
\sum_{\theta_{l,s}\in\mathcal{K}_{l,s}}\int_{0}^{t}\mathbf{F}_{2}\sum
_{w=2}^{n+1-j-s}\frac{1}{w}\right.  \right. \\
& \times\sum_{z=1}^{n+2-j-s-w}\binom{n+1-j-s}{z}m_{z}q_{w-1}^{\left(
n+1-j-s-z\right)  }t_{1}^{w}\mathbf{I}_{1}\ \mathrm{d}\hat{G}_{t_{1}}^{\left(
j\right)  }\\
& \left.  \left.  +1_{\left\{  s=n+1-j\right\}  }\sum_{l=1}^{s}\sum
_{\theta_{l,s}\in\mathcal{K}_{l,s}}\int_{0}^{t}\mathbf{F}_{1}\mathbf{I}%
_{1}\ \mathrm{d}\hat{G}_{t_{1}}^{\left(  j\right)  }\right\}  +1_{\left\{
j=n\right\}  }\left(  n+1\right)  \int_{0}^{t}\int_{0}^{t_{1}-}\mathrm{d}%
\hat{G}_{t_{2}}^{\left(  1\right)  }\mathrm{d}\hat{G}_{t_{1}}^{\left(
n\right)  }\right\}  .
\end{align*}
Let $u\in\left\{  1,2,...,n-2\right\}  ,\ v\in\left\{  1,2,...,n-2\right\}
,\ u+v\leq n-1,\ x\in\left\{  1,2,...,v\right\}  ,\ \beta\in\left\{
1,2,...,v+1-x\right\}  .$ \ In $L_{1}$, when $j=u,\ s=n-u-v$ (hence
$s\in\left\{  1,...,n-1-j\right\}  $)$,\ w=x+1$ (hence $w\in\left\{
2,...,n+1-j-s\right\}  $)$,\ z=\beta$ (hence $z\in\left\{
1,...,n+2-j-s-w\right\}  $),%

\begin{align}
L_{1}  & =\sum_{l=1}^{n-u-v}\sum_{\theta_{l,n-u-v}\in\mathcal{K}_{l,n-u-v}%
}\int_{0}^{t}\frac{\left(  n+1\right)  !}{i_{1}^{\theta_{l,n-u-v}}%
!i_{2}^{\theta_{l,n-u-v}}!\cdots i_{l}^{\theta_{l,n-u-v}}!u!}\frac{1}{\left(
v+1\right)  !}\frac{1}{x+1}\binom{v+1}{\beta}\nonumber\\
& \times m_{\beta}q_{x}^{\left(  v+1-\beta\right)  }t_{1}^{x+1}\int_{0}%
^{t_{1}-}\cdots\int_{0}^{t_{l}-}\mathrm{d}\hat{G}_{t_{l+1}}^{\left(
i_{1}^{\theta_{l,n-u-v}}\right)  }\cdots\mathrm{d}\hat{G}_{t_{2}}^{\left(
i_{l}^{\theta_{l,n-u-v}}\right)  }\mathrm{d}\hat{G}_{t_{1}}^{\left(  u\right)
}.\label{L1_3}%
\end{align}
By definition, since $s=n-u-v$ and $l\in\left\{  1,2,...,n-u-v\right\}  ,$
$\left(  i_{1}^{\theta_{l,n-u-v}},i_{2}^{\theta_{l,n-u-v}},...,i_{l}%
^{\theta_{l,n-u-v}}\right)  \in\mathcal{J}_{n-u-v}.$ \ In $L_{3}$, when
$j=\beta$ (hence $j\in\left\{  1,...,n-2\right\}  $)$,\ s=n-v$ (hence
$s\in\left\{  2,...,n-j\right\}  $)$,\ i_{l}^{\theta_{l,s}}=u$ (hence
$i_{l}^{\theta_{l,s}}<s$)$,\ w=x$ (hence $w\in\left\{  1,...,n+1-j-s\right\}
$)$,$%
\begin{align*}
L_{3}  & =-m_{\beta}\sum_{l=1}^{n-v}\sum_{\theta_{l,n-v}\in\mathcal{K}%
_{l,n-v}}\frac{\left(  n+1\right)  !}{i_{1}^{\theta_{l,n-v}}!i_{2}%
^{\theta_{l,n-v}}!\cdots i_{l-1}^{\theta_{l,n-v}}!u!\beta!}\frac{1}{\left(
v+1-\beta\right)  !}\\
& \times q_{x}^{\left(  v+1-\beta\right)  }\left\{  \frac{x}{x+1}t^{x+1}%
\int_{0}^{t}\int_{0}^{t_{1}-}\cdots\int_{0}^{t_{l-1}-}\mathrm{d}\hat{G}%
_{t_{l}}^{\left(  i_{1}^{\theta_{l,n-v}}\right)  }\cdots\mathrm{d}\hat{G}%
_{t2}^{\left(  i_{l-1}^{\theta_{l,n-v}}\right)  }\mathrm{d}\hat{G}_{t_{1}%
}^{\left(  u\right)  }\right. \\
& \left.  +\frac{1}{x+1}\int_{0}^{t}t_{1}^{x+1}\int_{0}^{t_{1}-}\cdots\int
_{0}^{t_{l-1}-}\mathrm{d}\hat{G}_{t_{l}}^{\left(  i_{1}^{\theta_{l,n-v}%
}\right)  }\cdots\mathrm{d}\hat{G}_{t2}^{\left(  i_{l-1}^{\theta_{l,n-v}%
}\right)  }\mathrm{d}\hat{G}_{t_{1}}^{\left(  u\right)  }\right\}  .
\end{align*}
The final term in $L_{3}$%
\begin{align*}
& -m_{\beta}\sum_{l=1}^{n-v}\sum_{\theta_{l,n-v}\in\mathcal{K}_{l,n-v}}%
\frac{\left(  n+1\right)  !}{i_{1}^{\theta_{l,n-v}}!i_{2}^{\theta_{l,n-v}%
}!\cdots i_{l-1}^{\theta_{l,n-v}}!u!\beta!}\frac{1}{\left(  v+1-\beta\right)
!}\\
& \ \ \times q_{x}^{\left(  v+1-\beta\right)  }\frac{1}{x+1}\int_{0}^{t}%
t_{1}^{x+1}\int_{0}^{t_{1}-}\cdots\int_{0}^{t_{l-1}-}\mathrm{d}\hat{G}_{t_{l}%
}^{\left(  i_{1}^{\theta_{l,n-v}}\right)  }\cdots\mathrm{d}\hat{G}%
_{t2}^{\left(  i_{l-1}^{\theta_{l,n-v}}\right)  }\mathrm{d}\hat{G}_{t_{1}%
}^{\left(  u\right)  }%
\end{align*}
clearly cancels (\ref{L1_3}) in $L_{1}.$ \ So now we can write%
\begin{align*}
L_{1}  & =\sum_{j=1}^{n}1_{\left\{  j\leq n\right\}  }\sum_{s=1}%
^{n+1-j}1_{\left\{  s=n+1-j\right\}  }\sum_{l=1}^{s}\sum_{\theta_{l,s}%
\in\mathcal{K}_{l,s}}\int_{0}^{t}\ \mathbf{F}_{1}\mathbf{I}_{1}\ \mathrm{d}%
\hat{G}_{t_{1}}^{\left(  j\right)  }.\\
L_{3}  & =-\sum_{j=1}^{n-1}m_{j}\sum_{s=1}^{n-j}\left\{  1_{\left\{
s=1\right\}  }\frac{\left(  n+1\right)  !}{j!}\frac{1}{\left(  n-j\right)
!}\sum_{w=1}^{n-j}q_{w}^{\left(  n-j\right)  }\frac{w}{w+1}t^{w+1}\int_{0}%
^{t}\mathrm{d}\hat{G}_{t_{1}}^{\left(  1\right)  }\right. \\
& +1_{\left\{  2\leq s\leq n-j\right\}  }\sum_{l=1}^{s}\sum_{\theta_{l,s}%
\in\mathcal{K}_{l,s}}\left\{  1_{\left\{  i_{l}^{\theta_{l,s}}<s\right\}
}\mathbf{F}_{2}\sum_{w=1}^{n+1-j-s}\mathbf{\tilde{q}}_{w}\frac{w}{w+1}%
t^{w+1}\mathbf{I}_{3}\right. \\
& \left.  \left.  +1_{\left\{  i_{l}^{\theta_{l,s}}=s\right\}  }\frac{\left(
n+1\right)  !}{s!j!}\frac{1}{\left(  n+1-j-s\right)  !}\sum_{w=1}%
^{n+1-j-s}\mathbf{\tilde{q}}_{w}\frac{w}{w+1}t^{w+1}\int_{0}^{t}\mathrm{d}%
\hat{G}_{t_{1}}^{\left(  s\right)  }\right\}  \right\}  .
\end{align*}
We can now simplify it as%
\[
L_{3}=-\sum_{j=1}^{n-1}m_{j}\sum_{s=1}^{n-j}\sum_{l=1}^{s}\sum_{\theta
_{l,s}\in\mathcal{K}_{l,s}}\mathbf{F}_{2}\sum_{w=1}^{n+1-j-s}\mathbf{\tilde
{q}}_{w}\frac{w}{w+1}t^{w+1}\ \mathbf{I}_{3}.
\]
All together, we have%
\begin{align*}
L_{1}  & =\sum_{j=1}^{n}1_{\left\{  j\leq n\right\}  }\left\{  \sum
_{s=1}^{n+1-j}1_{\left\{  s=n+1-j\right\}  }\sum_{l=1}^{s}\sum_{\theta
_{l,s}\in\mathcal{K}_{l,s}}\int_{0}^{t}\mathbf{F}_{1}\mathbf{I}_{1}%
\ \mathrm{d}\hat{G}_{t_{1}}^{\left(  j\right)  }\right\}  .\\
L_{2}  & =\sum_{j=1}^{n}m_{j}\sum_{s=1}^{n+1-j}\left\{  1_{\left\{
s=n+1-j\right\}  }\sum_{l=1}^{s}\sum_{\theta_{l,s}\in\mathcal{K}_{l,s}%
}t\mathbf{F}_{1}\mathbf{I}_{3}+1_{\left\{  s\leq n-j\right\}  }\sum_{l=1}%
^{s}\sum_{\theta_{l,s}\in\mathcal{K}_{l,s}}\mathbf{F}_{2}\sum_{w=1}%
^{n+1-j-s}\ \mathbf{\tilde{q}}_{w}t^{w+1}\mathbf{I}_{3}\right\}  .\\
L_{3}  & =-\sum_{j=1}^{n-1}m_{j}\sum_{s=1}^{n-j}\sum_{l=1}^{s}\sum
_{\theta_{l,s}\in\mathcal{K}_{l,s}}\mathbf{F}_{2}\sum_{w=1}^{n+1-j-s}%
\mathbf{\tilde{q}}_{w}\frac{w}{w+1}t^{w+1}\mathbf{I}_{3}.\\
L_{4}  & =0.\\
L_{5}  & =\sum_{j=1}^{n}\binom{n+1}{j}m_{j}\sum_{w=1}^{n+1-j}q_{w}^{\left(
n+1-j\right)  }t^{w+1}.\\
L_{6}  & =-\sum_{j=1}^{n}\binom{n+1}{j}m_{j}\sum_{w=1}^{n+1-j}q_{w}^{\left(
n+1-j\right)  }\frac{w}{w+1}t^{w+1}.
\end{align*}
Since at the beginning of the proof, we have already showed that the
stochastic integrals of $G_{t}^{n+1}$ are of the form $\mathcal{S}%
_{\theta_{n+1},t,0}$ where $\theta_{n+1,t}\in\mathcal{I}_{n+1}.$ \ We are now
going to show that the coefficient of each $\mathcal{S}_{\theta_{n+1},t,0}$ is
$\Pi_{\theta_{n+1,t}}^{\left(  n+1\right)  }.$ \ Consider $\int_{0}^{t_{1}%
-}\int_{0}^{t_{2}-}\cdots\int_{0}^{t_{l}-}\mathrm{d}\hat{G}_{t_{l+1}}^{\left(
i_{1}^{\theta_{l,s}}\right)  }\cdots\mathrm{d}\hat{G}_{t_{2}}^{\left(
i_{l}^{\theta_{l,s}}\right)  }\mathrm{d}\hat{G}_{t_{1}}^{\left(  j\right)  }$
where $\theta_{l,s}\in\mathcal{K}_{l,s},\ j\in\left\{  1,2,...,n\right\}
,\ s=n+1-j.$ \ This stochastic integral only appears in $L_{1}$. \ And its
coefficient is $\frac{\left(  n+1\right)  !}{i_{1}^{\theta_{l,s}}%
!i_{2}^{\theta_{l,s}}!\cdots i_{l}^{\theta_{l,s}}!j!}$. \ And from (\ref{PI}),
since $n+1-s-j=n+1-\left(  n+1-j\right)  -j=0,$
\[
\Pi_{\left(  i_{1}^{\theta_{l,s}},i_{2}^{\theta_{l,s}},...,i_{l}^{\theta
_{l,s}},j\right)  }^{\left(  n+1\right)  }=\frac{\left(  n+1\right)  !}%
{i_{1}^{\theta_{l,s}}!i_{2}^{\theta_{l,s}}!\cdots i_{l}^{\theta_{l,s}}%
!j!0!}C_{t}^{\left(  0\right)  }=\frac{\left(  n+1\right)  !}{i_{1}%
^{\theta_{l,s}}!i_{2}^{\theta_{l,s}}!\cdots i_{l}^{\theta_{l,s}}!j!}%
\]
since $C_{t}^{\left(  0\right)  }=1$ by definition (\ref{C}). \ Hence we have
proved that the coefficient is given by $\Pi_{\left(  i_{1}^{\theta_{l,s}%
},i_{2}^{\theta_{l,s}},...,i_{l}^{\theta_{l,s}},j\right)  }^{\left(
n+1\right)  }.$ \ Next, we change the summation sign of $j$ and $s$ in $L_{2}%
$.%
\begin{align*}
L_{2}  & =\sum_{s=1}^{n}\left\{  \sum_{l=1}^{s}\sum_{\theta_{l,s}%
\in\mathcal{K}_{l,s}}m_{n+1-s}t\frac{\left(  n+1\right)  !}{i_{1}%
^{\theta_{l,s}}!i_{2}^{\theta_{l,s}}!\cdots i_{l}^{\theta_{l,s}}!\left(
n+1-s\right)  !}\ \mathbf{I}_{3}\right. \\
& \left.  +\sum_{j=1}^{n-s}\sum_{l=1}^{s}\sum_{\theta_{l,s}\in\mathcal{K}%
_{l,s}}m_{j}\mathbf{F}_{2}\sum_{w=1}^{n+1-j-s}\mathbf{\tilde{q}}_{w}%
t^{w+1}\mathbf{I}_{3}\right\}  .
\end{align*}
Similarly, by changing the summation sign of $j$ and $w$, we have%
\begin{align*}
L_{2}  & =\sum_{s=1}^{n}\left\{  \sum_{l=1}^{s}\sum_{\theta_{l,s}%
\in\mathcal{K}_{l,s}}m_{n+1-s}t\frac{\left(  n+1\right)  !}{i_{1}%
^{\theta_{l,s}}!i_{2}^{\theta_{l,s}}!\cdots i_{l}^{\theta_{l,s}}!\left(
n+1-s\right)  !}\mathbf{I}_{3}\right. \\
& \left.  +\sum_{w=1}^{n-s}\sum_{j=1}^{n+1-w-s}\sum_{l=1}^{s}\sum
_{\theta_{l,s}\in\mathcal{K}_{l,s}}m_{j}\mathbf{F}_{2}\mathbf{\tilde{q}}%
_{w}t^{w+1}\mathbf{I}_{3}\right\}  .
\end{align*}
By (\ref{q}), $\frac{1}{w+1}\sum_{j=1}^{n+1-w-s}\frac{\left(  n+1-s\right)
!}{j!\left(  n+1-j-s\right)  !}m_{j}q_{w}^{\left(  n+1-s-j\right)  }%
=q_{w+1}^{\left(  n+1-s\right)  },$ so we have%
\begin{align*}
L_{2}  & =\sum_{s=1}^{n}\left\{  \sum_{l=1}^{s}\sum_{\theta_{l,s}%
\in\mathcal{K}_{l,s}}m_{n+1-s}t\frac{\left(  n+1\right)  !}{i_{1}%
^{\theta_{l,s}}!i_{2}^{\theta_{l,s}}!\cdots i_{l}^{\theta_{l,s}}!\left(
n+1-s\right)  !}\mathbf{I}_{3}\right. \\
& \left.  +\sum_{w=1}^{n-s}\sum_{l=1}^{s}\sum_{\theta_{l,s}\in\mathcal{K}%
_{l,s}}\frac{\left(  n+1\right)  !}{i_{1}^{\theta_{l,s}}!i_{2}^{\theta_{l,s}%
}!\cdots i_{l}^{\theta_{l,s}}!}\left(  w+1\right)  \frac{1}{\left(
n+1-s\right)  !}q_{w+1}^{\left(  n+1-s\right)  }t^{w+1}\mathbf{I}_{3}\right\}
.
\end{align*}
Changing $\sum_{w=1}^{n-s}$ to $\sum_{w=2}^{n+1-s}$, we have%
\begin{align*}
L_{2}  & =\sum_{s=1}^{n}\left\{  \sum_{l=1}^{s}\sum_{\theta_{l,s}%
\in\mathcal{K}_{l,s}}m_{n+1-s}t\frac{\left(  n+1\right)  !}{i_{1}%
^{\theta_{l,s}}!i_{2}^{\theta_{l,s}}!\cdots i_{l}^{\theta_{l,s}}!\left(
n+1-s\right)  !}\mathbf{I}_{3}\right. \\
& \left.  +\sum_{w=2}^{n+1-s}\sum_{l=1}^{s}\sum_{\theta_{l,s}\in
\mathcal{K}_{l,s}}\frac{\left(  n+1\right)  !}{i_{1}^{\theta_{l,s}}%
!i_{2}^{\theta_{l,s}}!\cdots i_{l}^{\theta_{l,s}}!}\frac{w}{\left(
n+1-s\right)  !}q_{w}^{\left(  n+1-s\right)  }t^{w}\mathbf{I}_{3}\right\}  .
\end{align*}
Similarly,%
\[
L_{3}=-\sum_{s=1}^{n-1}\sum_{w=1}^{n-s}\sum_{j=1}^{n+1-w-s}m_{j}\sum_{l=1}%
^{s}\sum_{\theta_{l,s}\in\mathcal{K}_{l,s}}\mathbf{F}_{2}\mathbf{\tilde{q}%
}_{w}\frac{w}{w+1}t^{w+1}\mathbf{I}_{3}.
\]
By (\ref{q}), $\frac{1}{w+1}\sum_{j=1}^{n+1-w-s}\frac{\left(  n+1-s\right)
!}{j!\left(  n+1-j-s\right)  !}m_{j}q_{w}^{\left(  n+1-s-j\right)  }%
=q_{w+1}^{\left(  n+1-s\right)  },$ so we have%

\[
L_{3}=-\sum_{s=1}^{n-1}\sum_{w=2}^{n+1-s}\sum_{l=1}^{s}\sum_{\theta_{l,s}%
\in\mathcal{K}_{l,s}}\frac{\left(  n+1\right)  !}{i_{1}^{\theta_{l,s}}%
!i_{2}^{\theta_{l,s}}!\cdots i_{l}^{\theta_{l,s}}!}\frac{w-1}{\left(
n+1-s\right)  !}q_{w}^{\left(  n+1-s\right)  }t^{w}\mathbf{I}_{3}%
\]
For $s=1$, the stochastic integral $\int_{0}^{t}\mathrm{d}\hat{G}_{t_{1}%
}^{\left(  1\right)  }$ appears in both $L_{2}$ and $L_{3}$. \ Its coefficient
is given by%
\begin{align*}
& \sum_{w=2}^{n}\left(  n+1\right)  wq_{w}^{\left(  n\right)  }t^{w}%
+m_{n}\left(  n+1\right)  t-\left(  n+1\right)  \sum_{w=2}^{n}\left(
w-1\right)  q_{w}^{\left(  n\right)  }t^{w}\\
& \ =\left(  n+1\right)  \left[  m_{n}t+\sum_{w=2}^{n}q_{w}^{\left(  n\right)
}t^{w}\right]  =\left(  n+1\right)  C_{t}^{\left(  n\right)  }.
\end{align*}
By (\ref{PI}),
\[
\Pi_{\left(  1\right)  }^{\left(  n+1\right)  }=\frac{\left(  n+1\right)
!}{\left(  n+1-1\right)  !}C_{t}^{\left(  n+1-1\right)  }=\left(  n+1\right)
C_{t}^{\left(  n\right)  }.
\]
For $s\in\left\{  2,3,...,n-1\right\}  $, the coefficients of the stochastic
integral
\[
\int_{0}^{t}\int_{0}^{t_{1}-}\cdots\int_{0}^{t_{l-1}-}\mathrm{d}\hat{G}%
_{t_{l}}^{\left(  i_{1}^{\theta_{l,s}}\right)  }\cdots\mathrm{d}\hat{G}%
_{t_{2}}^{\left(  i_{l-1}^{\theta_{l,s}}\right)  }\mathrm{d}\hat{G}_{t_{1}%
}^{\left(  i_{l}^{\theta_{l,s}}\right)  }%
\]
is given by%
\begin{align*}
& m_{n+1-s}t\frac{\left(  n+1\right)  !}{i_{1}^{\theta_{l,s}}!i_{2}%
^{\theta_{l,s}}!\cdots i_{l}^{\theta_{l,s}}!\left(  n+1-s\right)  !}%
+\sum_{w=2}^{n+1-s}\frac{\left(  n+1\right)  !}{i_{1}^{\theta_{l,s}}%
!i_{2}^{\theta_{l,s}}!\cdots i_{l}^{\theta_{l,s}}!}\frac{w}{\left(
n+1-s\right)  !}q_{w}^{\left(  n+1-s\right)  }t^{w}\\
& -\sum_{w=2}^{n+1-s}\frac{\left(  n+1\right)  !}{i_{1}^{\theta_{l,s}}%
!i_{2}^{\theta_{l,s}}!\cdots i_{l}^{\theta_{l,s}}!}\frac{\left(  w-1\right)
}{\left(  n+1-s\right)  !}q_{w}^{\left(  n+1-s\right)  }t^{w}\\
& \ =\sum_{w=1}^{n+1-s}\frac{\left(  n+1\right)  !}{i_{1}^{\theta_{l,s}}%
!i_{2}^{\theta_{l,s}}!\cdots i_{l}^{\theta_{l,s}}!}\frac{1}{\left(
n+1-s\right)  !}q_{w}^{\left(  n+1-s\right)  }t^{w}\\
& \ =\frac{\left(  n+1\right)  !}{i_{1}^{\theta_{l,s}}!i_{2}^{\theta_{l,s}%
}!\cdots i_{l}^{\theta_{l,s}}!}\frac{1}{\left(  n+1-s\right)  !}C_{t}^{\left(
n+1-s\right)  }\ =\Pi_{\left(  i_{1}^{\theta_{l,s}},i_{2}^{\theta_{l,s}%
},...,i_{l}^{\theta_{l,s}}\right)  }^{\left(  n+1\right)  }%
\end{align*}
by (\ref{PI}). \ For $s=n$, the stochastic integral appears in $L_{2}$ only
and its coefficient is given by%
\[
m_{1}t\frac{\left(  n+1\right)  !}{i_{1}^{\theta_{l,s}}!i_{2}^{\theta_{l,s}%
}!\cdots i_{l}^{\theta_{l,s}}!}=\frac{\left(  n+1\right)  !}{i_{1}%
^{\theta_{l,s}}!i_{2}^{\theta_{l,s}}!\cdots i_{l}^{\theta_{l,s}}!}%
C_{t}^{\left(  1\right)  }=\Pi_{\left(  i_{1}^{\theta_{l,s}},i_{2}%
^{\theta_{l,s}},...,i_{l}^{\theta_{l,s}}\right)  }^{\left(  n+1\right)  }.
\]
The stochastic integral $\int_{0}^{t}\mathrm{d}\hat{G}_{t_{1}}^{\left(
n+1\right)  }$ appears only once in $G_{t}^{n+1}$ and its coefficient is equal
to one. \ By (\ref{PI}),%
\[
\Pi_{\left(  n+1\right)  }^{\left(  n+1\right)  }=\frac{\left(  n+1\right)
!}{\left(  n+1\right)  !}C_{t}^{\left(  0\right)  }=1.
\]
Finally, we have to show that $L_{5}+L_{6}+m_{n+1}t=C_{t}^{\left(  n+1\right)
}.$ \ By (\ref{q}), \ $\frac{1}{w+1}\sum_{j=1}^{n+1-w}\binom{n+1}{j}m_{j}%
q_{w}^{\left(  n+1-j\right)  }=q_{w+1}^{\left(  n+1\right)  },$
\begin{align*}
L_{5}  & =\sum_{j=1}^{n}\binom{n+1}{j}m_{j}\sum_{w=1}^{n+1-j}q_{w}^{\left(
n+1-j\right)  }t^{w+1}=\sum_{w=1}^{n}\sum_{j=1}^{n+1-w}\binom{n+1}{j}%
m_{j}q_{w}^{\left(  n+1-j\right)  }t^{w+1}=\sum_{w=1}^{n}\left(  w+1\right)
q_{w+1}^{\left(  n+1\right)  }t^{w+1}.\\
L_{6}  & =-\sum_{w=1}^{n}\sum_{j=1}^{n+1-w}\binom{n+1}{j}m_{j}q_{w}^{\left(
n+1-j\right)  }\frac{w}{w+1}t^{w+1}=-\sum_{w=1}^{n}wq_{w+1}^{\left(
n+1\right)  }t^{w+1}.
\end{align*}
Hence
\[
L_{5}+L_{6}+m_{n+1}t=\sum_{w=1}^{n}q_{w+1}^{\left(  n+1\right)  }%
t^{w+1}+m_{n+1}t=\sum_{w=2}^{n+1}q_{w}^{\left(  n+1\right)  }t^{w}%
+m_{n+1}t=\sum_{w=1}^{n+1}q_{w}^{\left(  n+1\right)  }t^{w}=C_{n+1}^{\left(
k\right)  }.
\]
Thus, we have proved that%
\[
G_{t}^{n+1}=\sum_{\theta_{n+1}\in\mathcal{I}_{n+1}}\Pi_{\theta_{n+1}%
,t}^{\left(  n+1\right)  }\mathcal{S}_{\theta_{n+1},t,0}+C_{t}^{\left(
n+1\right)  }.
\]
As explained in (\ref{FG}), since $F_{t}=G_{t+t_{0}}-G_{t_{0}}$ is also a
L\'{e}vy process, we can write%
\[
F_{t}^{n+1}=\sum_{s=1}^{n+1}\sum_{l=1}^{s}\sum_{\theta_{l,s}\in\mathcal{K}%
_{l,s}}\Pi_{\theta_{l,s},t}^{\left(  n+1\right)  }\int_{0}^{t}\int_{0}%
^{t_{1}-}\cdots\int_{0}^{t_{l-1}-}\mathrm{d}\hat{F}_{t_{l}}^{\left(
i_{1}^{\theta_{l,s}}\right)  }\cdots\mathrm{d}\hat{F}_{t_{2}}^{\left(
i_{l-1}^{\theta_{l,s}}\right)  }\mathrm{d}\hat{F}_{t_{1}}^{\left(
i_{l}^{\theta_{l,s}}\right)  }%
\]
and since $\mathrm{d}\hat{F}_{t}^{\left(  i\right)  }=\mathrm{d}\left(
\hat{G}_{t+t_{0}}^{\left(  i\right)  }-\hat{G}_{t_{0}}^{\left(  i\right)
}\right)  =\mathrm{d}\hat{G}_{t+t_{0}}^{\left(  i\right)  },$ by changing of
variables, we have%
\begin{align*}
\left(  G_{t+t_{0}}-G_{t_{0}}\right)  ^{n+1}  & =\sum_{s=1}^{n+1}\sum
_{l=1}^{s}\sum_{\theta_{l,s}\in\mathcal{K}_{l,s}}\Pi_{\theta_{l,s},t}^{\left(
n+1\right)  }\int_{t_{0}}^{t+t_{0}}\int_{t_{0}}^{t_{1}-}\cdots\int_{t_{0}%
}^{t_{l-1}-}\mathrm{d}\hat{G}_{t_{l}}^{\left(  i_{1}^{\theta_{l,s}}\right)
}\cdots\mathrm{d}\hat{G}_{t_{2}}^{\left(  i_{l-1}^{\theta_{l,s}}\right)
}\mathrm{d}\hat{G}_{t_{1}}^{\left(  i_{l}^{\theta_{l,s}}\right)  }\\
& =\sum_{\theta_{n+1}\in\mathcal{I}_{n+1}}\Pi_{\theta_{n+1},t}^{\left(
n+1\right)  }\mathcal{S}_{\theta_{n+1},t,t_{0}}+C_{t}^{\left(  n+1\right)  }.
\end{align*}
Therefore, by the principle of strong induction,
\[
\left(  G_{t+t_{0}}-G_{t_{0}}\right)  ^{k}=\sum_{\theta_{k}\in\mathcal{I}_{k}%
}\Pi_{\theta_{k},t}^{\left(  k\right)  }\mathcal{S}_{\theta_{k},t,t_{0}}%
+C_{t}^{\left(  k\right)  }%
\]
for all non-negative integers $k$.

\setcounter{equation}{0}\renewcommand{\theequation}{D.\arabic{equation}}

\section{Proof of Proposition \ref{PropositionYtoH}\label{AppendixYtoH}}

We prove by induction. \ Assume the proposition is true for all $k\geq n.$
\ Now, consider $n+1,$%
\begin{align*}
Y^{\left(  n+1\right)  }  & =H^{\left(  n+1\right)  }-\sum_{l=1}^{n}%
a_{n+1,l}Y^{\left(  l\right)  }=H^{\left(  n+1\right)  }-\sum_{l=1}%
^{n}a_{n+1,l}\left\{  H^{\left(  l\right)  }+\sum_{k=1}^{l-1}b_{l,k}H^{\left(
k\right)  }\right\} \\
& =H^{\left(  n+1\right)  }-\sum_{l=1}^{n}a_{n+1,l}\sum_{k=1}^{l}%
b_{l,k}H^{\left(  k\right)  }=H^{\left(  n+1\right)  }+\sum_{k=1}^{n}%
b_{n+1,k}H^{\left(  k\right)  },
\end{align*}
which completes the proof.

\setcounter{equation}{0}\renewcommand{\theequation}{E.\arabic{equation}}

\section{Plots\label{AppendixPlots}}%

\begin{center}%
%

{\parbox[b]{3.1228in}{\begin{center}
\includegraphics[
height=2.2044in,
width=3.1228in
]%
{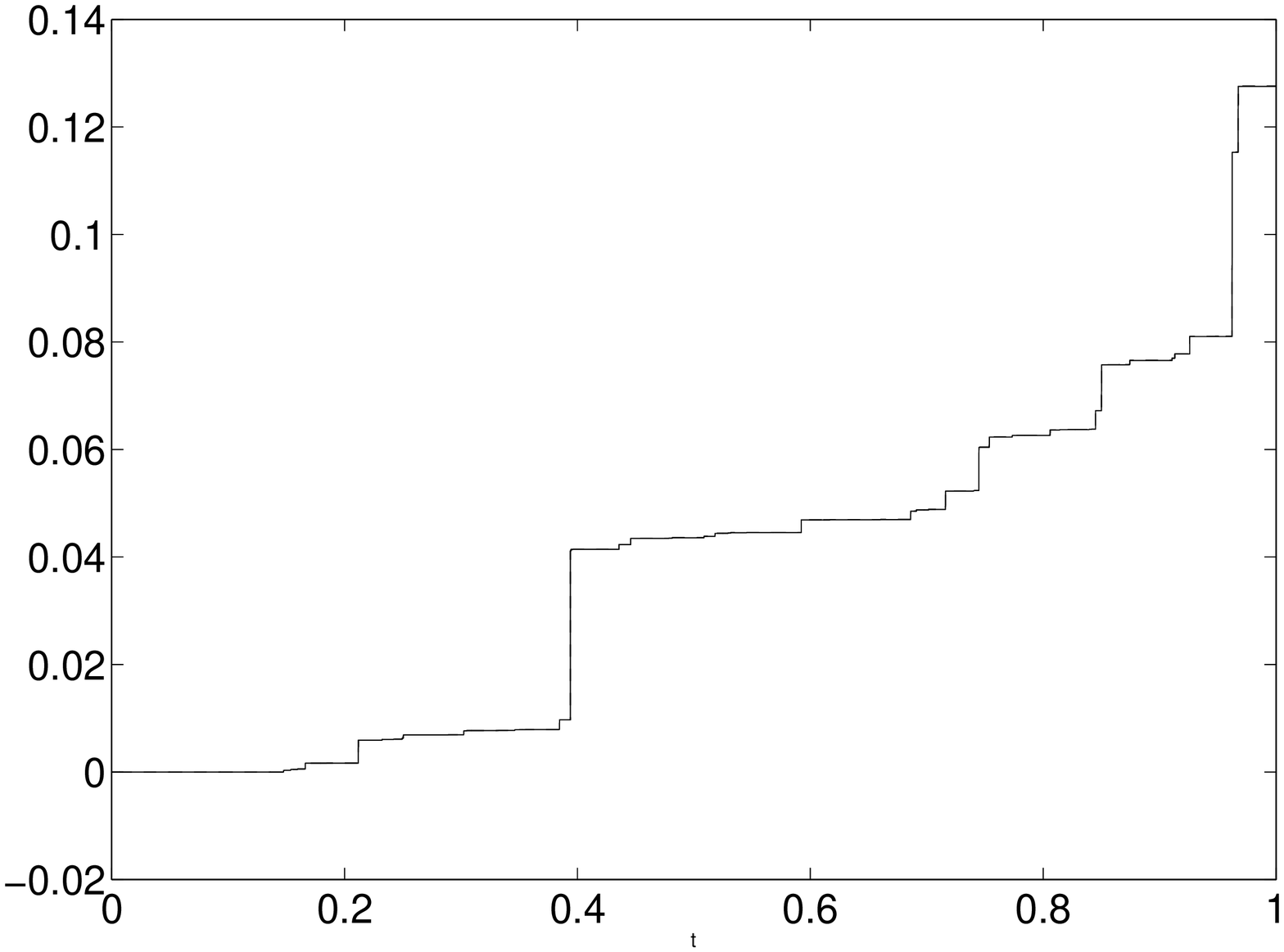}%
\\
Figure 1: $G_{t}^{4}$ generated using CRP and directly from the Gamma process.
\end{center}}}%
{\parbox[b]{3.1228in}{\begin{center}
\includegraphics[
height=2.2044in,
width=3.1228in
]%
{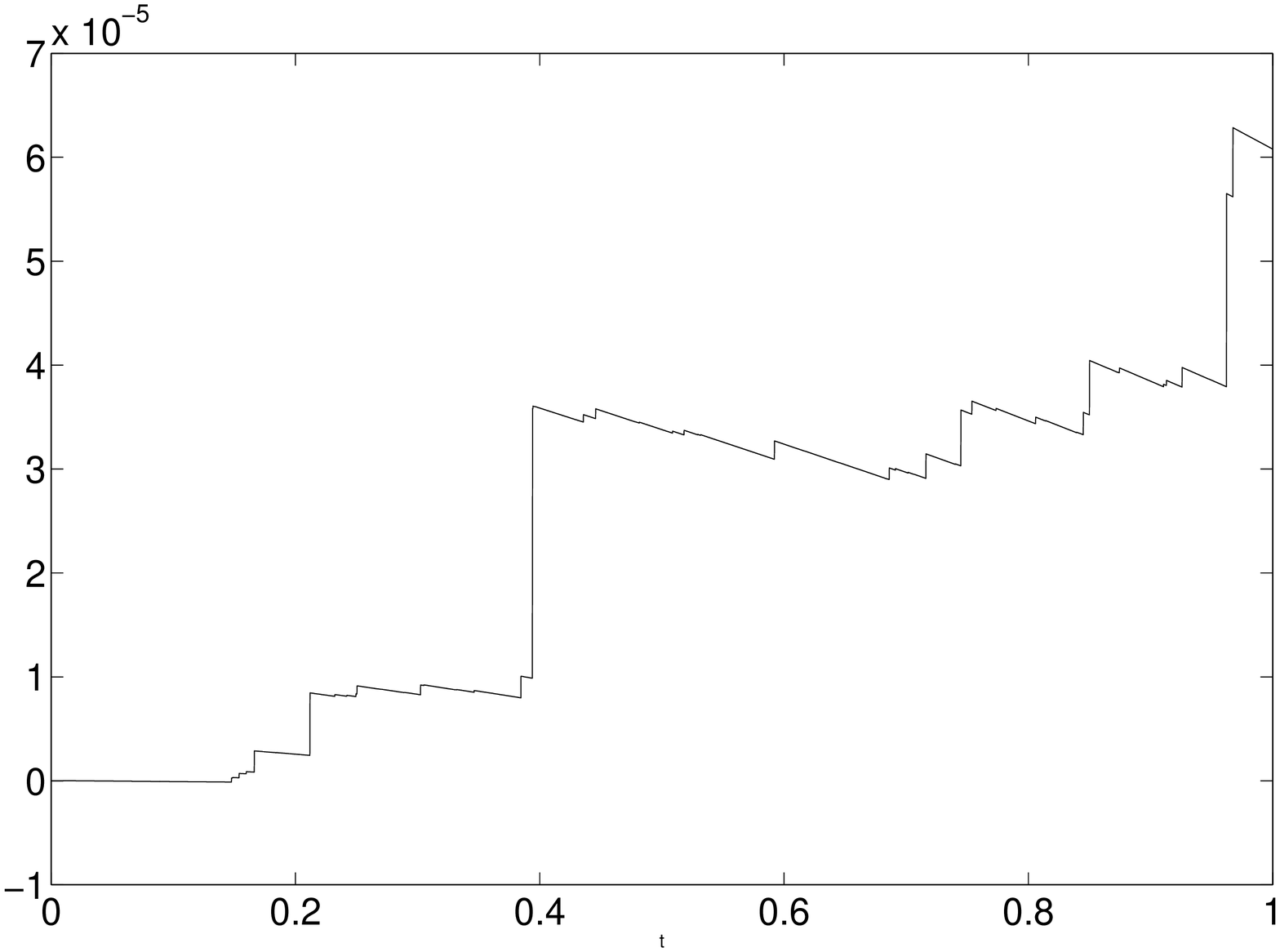}%
\\
Figure 2: The difference of the two series in Figure 1.
\end{center}}}%
%

{\parbox[b]{3.122in}{\begin{center}
\includegraphics[
height=2.2044in,
width=3.122in
]%
{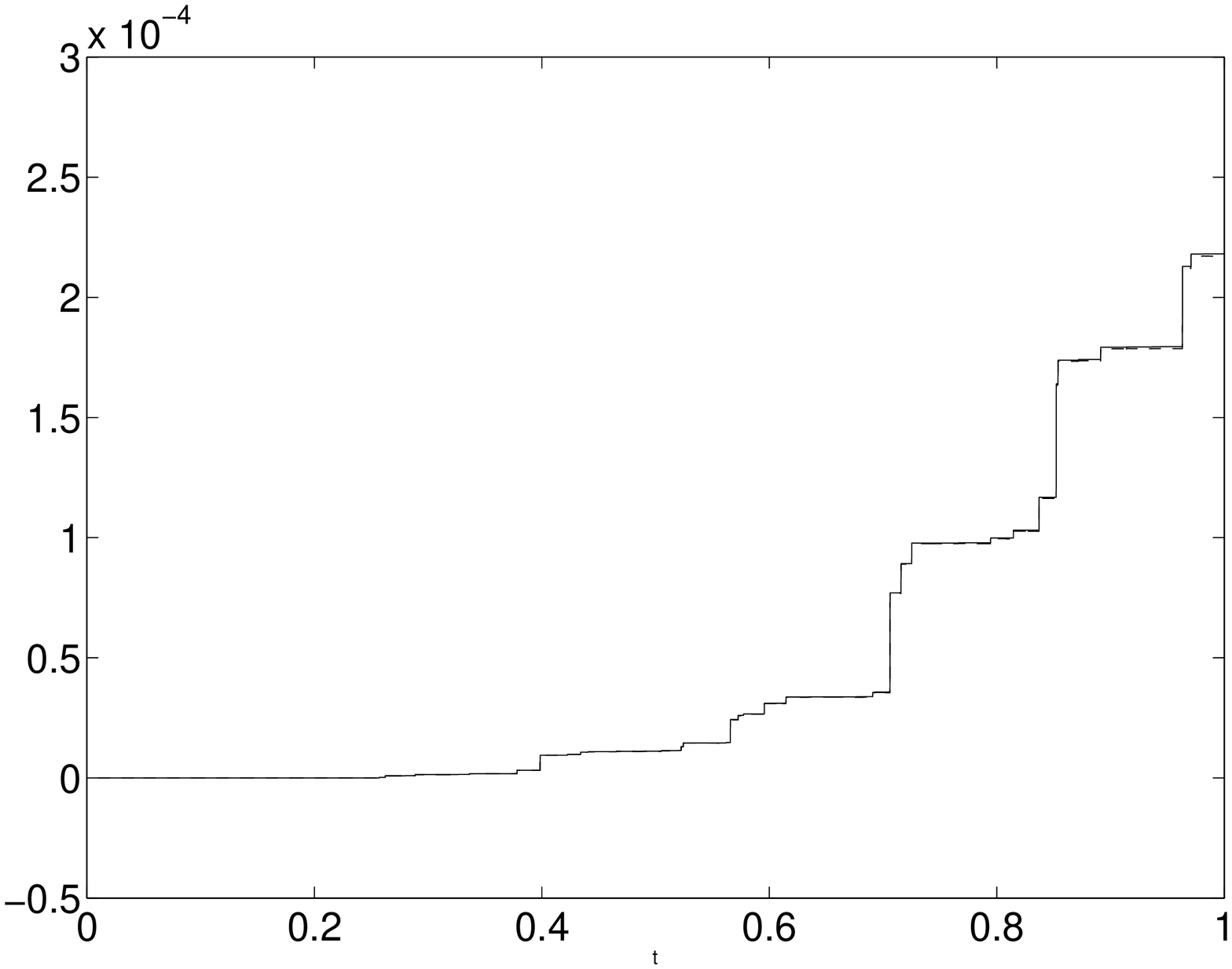}%
\\
Figure 3: $\left(  G_{t+t_{0}}-G_{t_{0}}\right)  ^{9}$ generated using CRP and
directly from the Gamma process.
\end{center}}}%
{\parbox[b]{3.122in}{\begin{center}
\includegraphics[
height=2.2044in,
width=3.122in
]%
{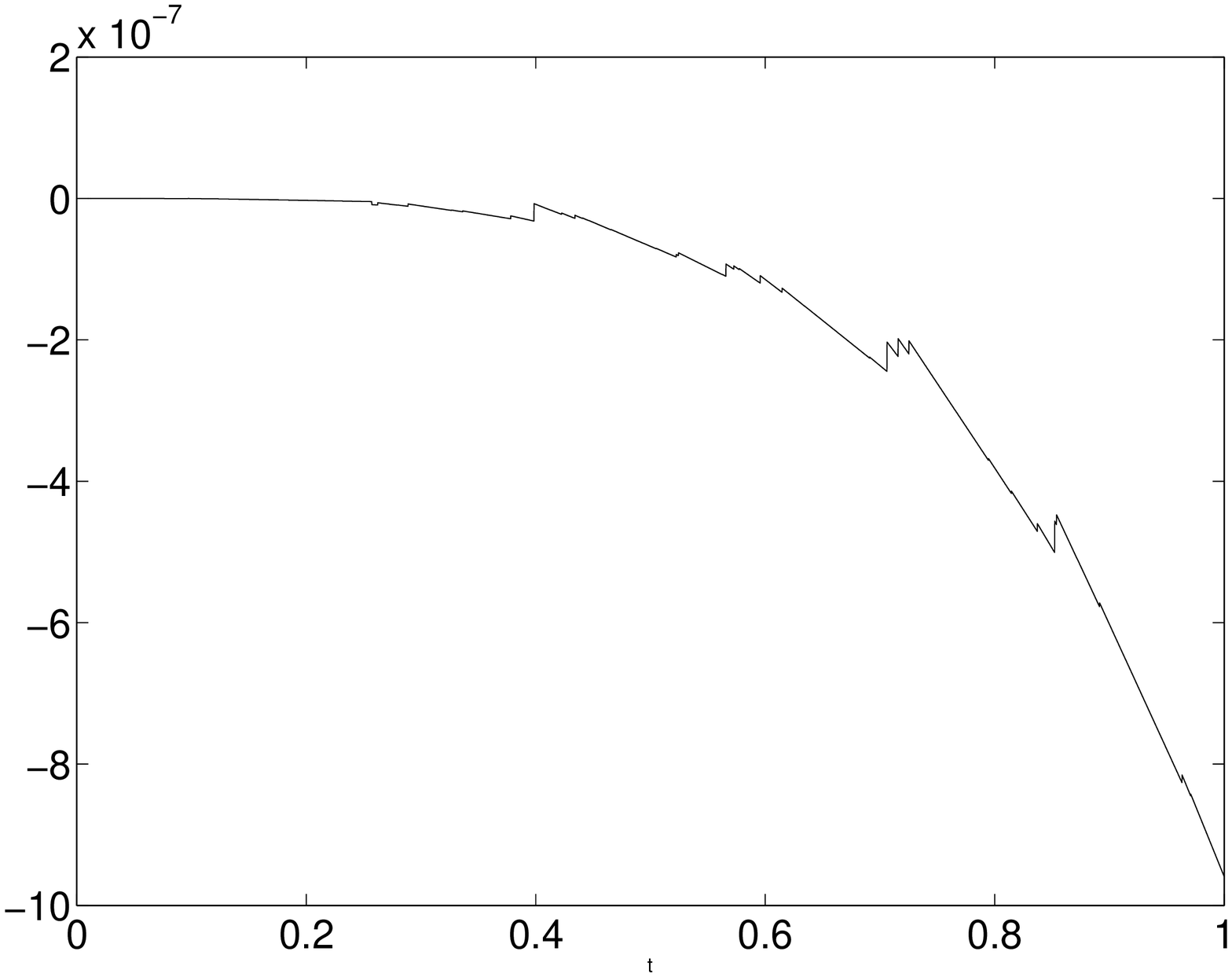}%
\\
Figure 4: The difference of the two series in Figure 3.
\end{center}}}%
%

{\parbox[b]{3.122in}{\begin{center}
\includegraphics[
height=2.2044in,
width=3.122in
]%
{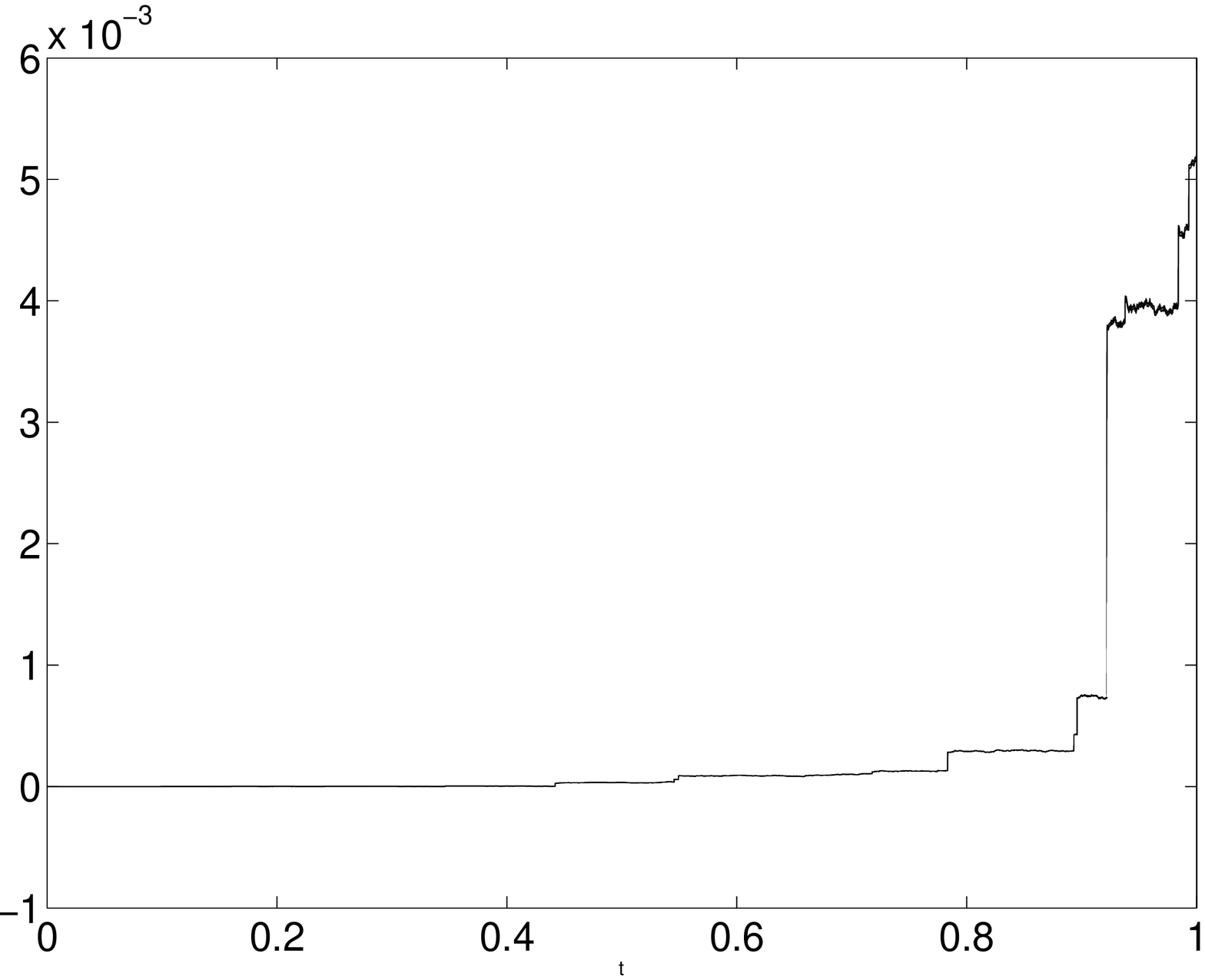}%
\\
Figure 5: $X_{t}^{5}$ generated using CRP and directly from the Wiener and
Gamma processes.
\end{center}}}%
{\parbox[b]{3.1228in}{\begin{center}
\includegraphics[
height=2.2044in,
width=3.1228in
]%
{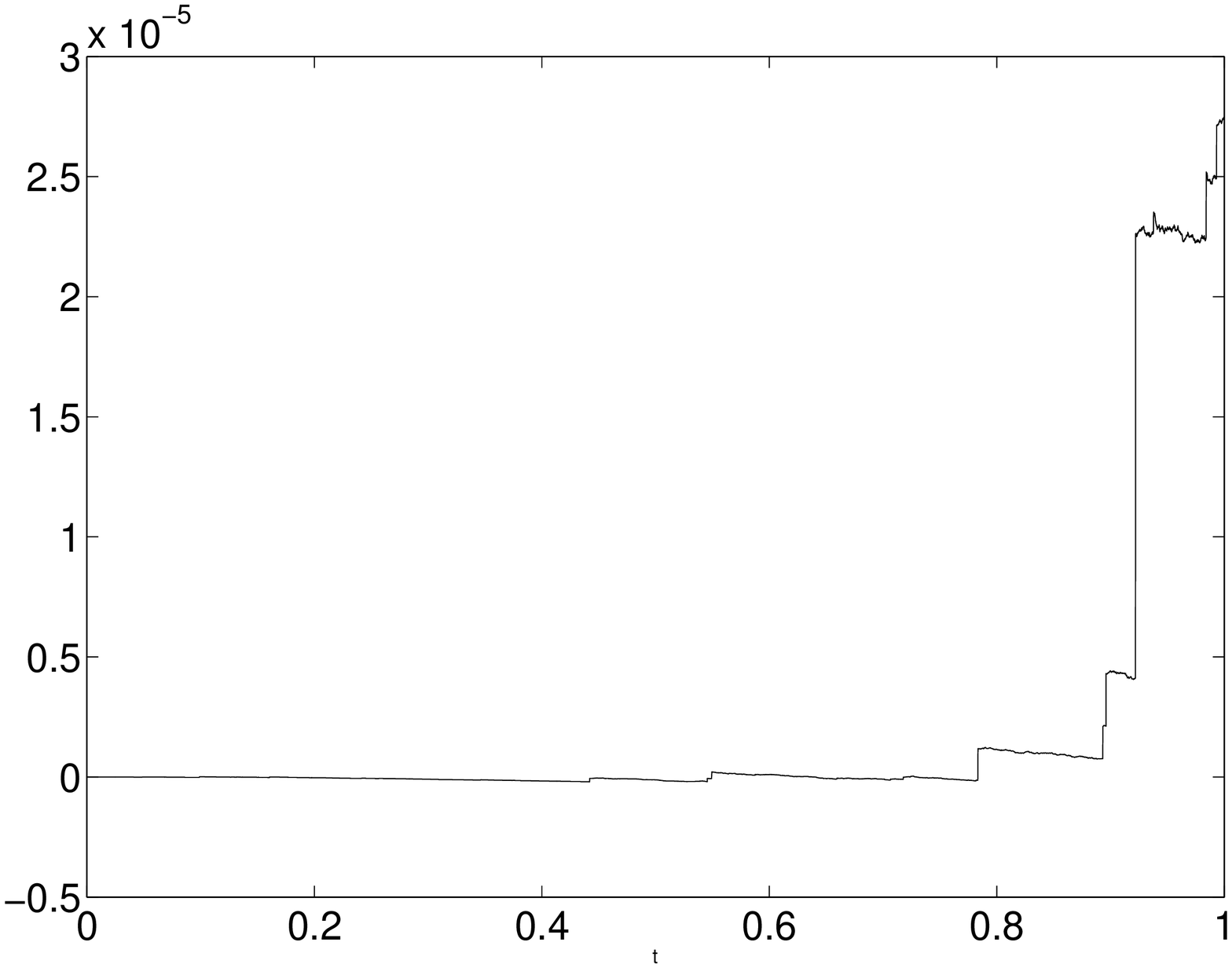}%
\\
Figure 6: The difference of the two series in Figure 5.
\end{center}}}%
%

{\parbox[b]{3.122in}{\begin{center}
\includegraphics[
height=2.2044in,
width=3.122in
]%
{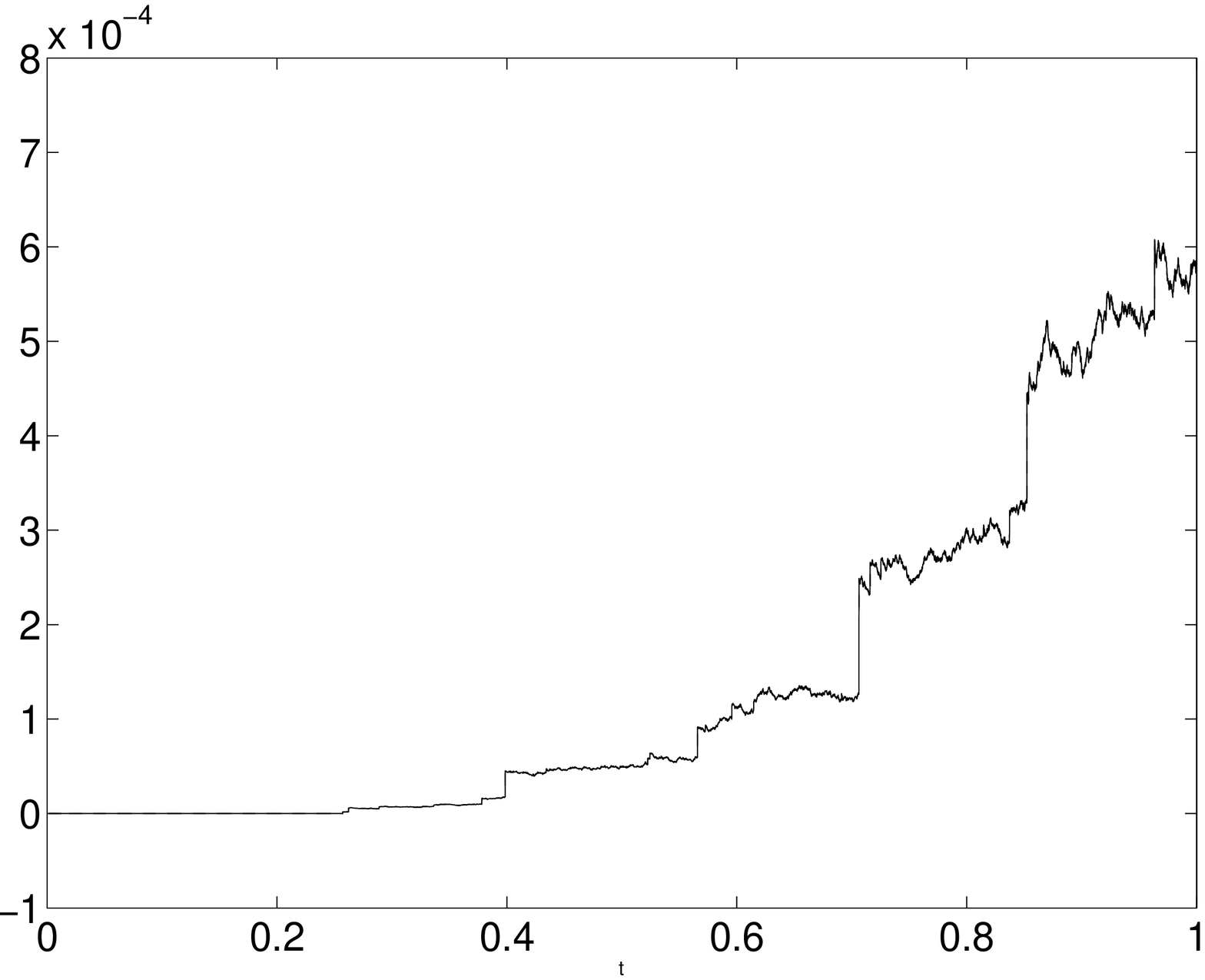}%
\\
Figure 7: $\left(  X_{t+t_{0}}-X_{t_{0}}\right)  ^{8}$ generated using CRP and
directly from the Wiener and Gamma processes.
\end{center}}}%
{\parbox[b]{3.1228in}{\begin{center}
\includegraphics[
height=2.2044in,
width=3.1228in
]%
{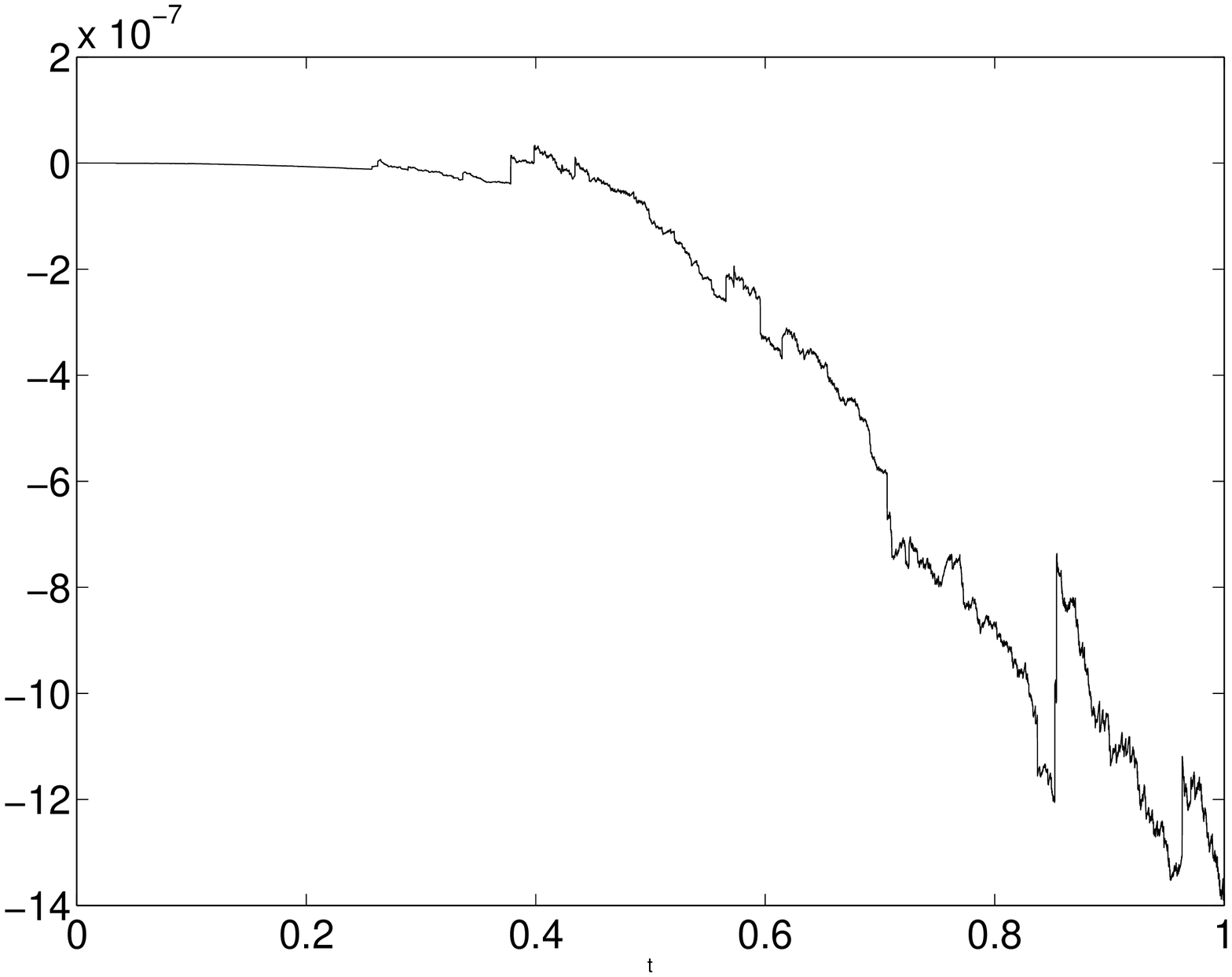}%
\\
Figure 8: The difference of the two series in Figure 7.
\end{center}}}%
%

\end{center}%

Figures 1-8: Solid line is generated using the CRP and the dotted line is
generated by the Wiener and Gamma processes. \ Time step $=\frac{1}%
{10000},\ a=10,\ b=20.$ \ In Figure 3, $t_{0}=0.0099$; in Figure 5,
$\sigma=0.01;$ in Figure 7, $t_{0}=0.0019$ and $\sigma=0.02.$

\bibliographystyle{authordate4}
\bibliography{MyReferences}

\bigskip
\end{document}